\documentclass[a4paper]{article}[12pt]
\usepackage{amsmath,cite,tikz,multirow,amssymb,mathrsfs,latexsym,bm,inputenc,longtable}

\begin{document}
\newcommand{\qed}{\hphantom{.}\hfill $\Box$\medbreak}
\newcommand{\proof}{\noindent{\bf Proof \ }}
\newtheorem{Theorem}{Theorem}[section]
\newtheorem{Lemma}[Theorem]{Lemma}
\newtheorem{Corollary}[Theorem]{Corollary}
\newtheorem{Remark}[Theorem]{Remark}
\newtheorem{Example}[Theorem]{Example}
\newtheorem{Definition}[Theorem]{Definition}
\newtheorem{Construction}[Theorem]{Construction}
\newtheorem{Proposition}[Theorem]{Proposition}

\title{\large{\bf Triangle decompositions of $\lambda K_v-\lambda K_w-\lambda K_u$ \footnote{Supported by NSFC under Grant 11971053 (Y. Chang), and NSFC under Grant 11871095 (T. Feng).}}}

\author{\small Yueting Li, Yanxun Chang, Tao Feng \\
\footnotesize  Department of Mathematics, Beijing Jiaotong University, Beijing 100044, P. R. China\\
\footnotesize  yuetli@bjtu.edu.cn, yxchang@bjtu.edu.cn, tfeng@bjtu.edu.cn}

\date{}
\maketitle

\noindent {\bf Abstract:} Denote by $\lambda K_v$ the complete graph of order $v$ with multiplicity $\lambda$. Let $\lambda K_v-\lambda K_w-\lambda K_u$ be the graph obtained from $\lambda K_v$ by the removal of the edges of two vertex disjoint complete multi-subgraphs with multiplicity $ \lambda $ of orders $ w $ and $ u $, respectively. When $\lambda$ is odd, it is shown that there exists a triangle decomposition of $\lambda K_v-\lambda K_w-\lambda K_u$ if and only if $v\geq w+u+\max\{u,w\}$, $ \lambda \left({v\choose 2}-{u\choose 2}-{w\choose 2}\right) \equiv 0 \pmod 3$ and $\lambda (v-w) \equiv \lambda (v-u) \equiv \lambda (v-1) \equiv 0 \pmod 2$. When $\lambda$ is even, it is shown that for large enough $v$, the elementary necessary conditions for the existence of a triangle decomposition of $\lambda K_v-\lambda K_w-\lambda K_u$ are also sufficient.

\vskip12pt


\noindent {\bf Keywords}$\colon$triangle decomposition; incomplete triple system; group divisible design


\section{Introduction}

All graphs (multigraphs) considered in this paper are loopless. A {\em $\triangle$-decomposition} of a multigraph $ G $ is a set $ \mathcal{T} $ of triangles in $ G $ such that each edge of $G$ occurs in exactly one triangle of $ \mathcal{T} $. Denote by $\lambda K_v$ the complete graph of order $v$ with multiplicity $\lambda$. Hanani \cite{h61} in 1961 proved that there exists a $ \triangle $-decomposition of $ \lambda K_v  $ if and only if $ \lambda\equiv0\pmod{\gcd(v-2,6)}$ and $ v\neq2 $.

Let $\lambda K_{v}-\lambda K_{w}$ be the graph obtained from $\lambda K_v$ by the removal of the edges of a complete multi-subgraph with multiplicity $ \lambda $ of order $ w $. In 1979, Stern \cite{s} established necessary and sufficient conditions for the existence of a $\triangle$-decomposition of $ \lambda K_{v}-\lambda K_{w} $.

\begin{Theorem} \label{ITS-one hole} {\rm\cite{s}}
There exists a $\triangle$-decomposition of $\lambda K_v-\lambda K_w$ if and only if $(1)$ $w=v$, or $(2)$ $w=0$, $v\neq 2$ and $\lambda \equiv 0 \pmod{{\rm gcd}(v-2,6)}$, or $(3)$ $0<w<v$, $v\geq 2w+1$, $\lambda({v\choose 2}-{w\choose 2}) \equiv 0 \pmod{3}$ and $\lambda(v-1)\equiv \lambda(v-w)\equiv0 \pmod{2}$.
\end{Theorem}

Furthermore, the graph obtained from $\lambda K_v$ by the removal of the edges of two vertex disjoint complete multi-subgraphs with multiplicity $ \lambda $ of orders $ w $ and $ u $, respectively, is denoted by $ \lambda K_{v}-\lambda K_{w}-\lambda K_{u} $. When $\lambda=1$, Bryant and Horsley \cite{bh} proved the following theorem.

\begin{Theorem} \label{lambda=1} {\rm\cite{bh}}
Let $v$, $w$ and $u$ be positive integers. There exists a $\triangle$-decomposition of $ K_v- K_w- K_u$ if and only if $v$, $u$ and $ w $ are odd, $ {v\choose 2}-{u\choose 2}-{w\choose 2} \equiv 0 \pmod 3$, and $v\geq w+u+\max\{u,w\}$.
\end{Theorem}

This paper is devoted to examining the existence of $\triangle$-decompositions of $\lambda K_v-\lambda K_w-\lambda K_u$. This question is equivalent to the existence of a class of incomplete triple systems (ITSs).

A {\em $ \Diamond $-ITS$_\lambda (v;w,u;z)$} is a quadruple $(V;W,U;\mathcal{B}) $, where $V$ is a $v$-set, $W, U\subseteq V$ (called {\em holes}) with $|W|=w$ and $|U|=u$ for which $ |W\cap U|=z $, and $\mathcal{B}$ is a collection of 3-subsets (called {\em triples} or {\em blocks}) of $V$ in which every pair of distinct elements from $V$ not contained within $U$ or $W$ appears in exactly $ \lambda $ triples in $\mathcal{B}$, while pairs within $U$ or $W$ appear in no triple. When $ z=0 $, it is simply written as a $\Diamond$-ITS$_\lambda(v;w,u)$. Clearly, the existence of $\Diamond$-ITS$_\lambda(v;w,u)$s is equivalent to the existence of $\triangle$-decompositions of $ \lambda K_{v}-\lambda K_{w}-\lambda K_{u} $. In the sequel, we use both representations, depending on what is more convenient in the context.

Throughout this paper,  every union will be understood as multiset union with multiplicities of elements preserved.

As the main result of the paper, we are to prove the following theorems.

\begin{Theorem} \label{main result}
Let $v$, $w$ and $u$ be positive integers. Let $\lambda\equiv 1\pmod{2}$. There exists a $\triangle$-decomposition of $\lambda K_v-\lambda K_w-\lambda K_u$, i.e., a $ \Diamond $-ITS$_\lambda (v;w,u)$, if and only if $v\geq w+u+\max\{u,w\}$, $ \lambda \left({v\choose 2}-{u\choose 2}-{w\choose 2}\right) \equiv 0 \pmod 3$ and $\lambda (v-w) \equiv \lambda (v-u) \equiv \lambda (v-1) \equiv 0 \pmod 2$.
\end{Theorem}

\begin{Theorem} \label{main result even}
Let $v$, $w$ and $u$ be positive integers such that
\begin{displaymath}
v\geq \left\{ \begin{array}{ll}
	w+4u+2, & \textrm{if $w<2u+1$},\\
	2w+2u+1, & \textrm{if $w\geq 2u+1$}.\\
\end{array} \right.
\end{displaymath}
There exists a $\triangle$-decomposition of $\lambda K_v-\lambda K_w-\lambda K_u$, i.e., a $ \Diamond $-ITS$_\lambda (v;w,u)$, if and only if $ \lambda \left({v\choose 2}-{u\choose 2}-{w\choose 2}\right) \equiv 0 \pmod 3$ and $\lambda (v-w) \equiv \lambda (v-u) \equiv \lambda (v-1) \equiv 0 \pmod 2$.
\end{Theorem}

\section{Preliminary}

\subsection{Necessary conditions}\label{sec:nec}

\begin{Lemma}\label{necessary}
Let $v$, $w$ and $u$ be positive integers. If $\lambda K_v-\lambda K_w-\lambda K_u$  has  a $\triangle$-decomposition, then $v\geq w+u+\max\{u,w\}$, $\lambda\left({v\choose 2}-{u\choose 2}-{w\choose 2}\right) \equiv 0\pmod3$ and $\lambda (v-w) \equiv \lambda (v-u) \equiv \lambda (v-1) \equiv 0 \pmod 2$.
\end{Lemma}

\proof Let $(V;W,U;\mathcal{B})$ be a $ \Diamond $-ITS$_\lambda (v;w,u)$, where $W=\{x_1,x_2,\ldots,x_w \}$, $U=\{y_1,y_2,\ldots,y_u\}$ and $V\setminus(W \cup U)=\{z_1,z_2,\ldots,z_{v-w-u}\}$. Fix an element $y\in U$, and consider the triples of the form $\{x_i,y,z_j\}$ in $\mathcal{B}$ for $1\leq i\leq w$ and $1\leq j\leq v-w-u$. It follows that $V\setminus(W \cup U)$ contains at least $w$ elements. Similar argument holds when we fix an element $x\in W$. Hence $v-w-u\geq \max\{u,w\}$.

Counting pairs that are to appear in triples of $\mathcal{B}$, we have $\lambda\left({v\choose 2}-{u\choose 2}-{w\choose 2}\right) \equiv 0\pmod{3}$. Counting triples involving a specified element, we have $\lambda (v-w) \equiv \lambda (v-u) \equiv \lambda (v-1) \equiv 0 \pmod 2$. \qed

For convenience, in what follows we always assume that $w\geq u$.

\begin{Lemma} \label{w=u=2}
There is no $\triangle$-decomposition of $\lambda K_6-\lambda K_2-\lambda K_2$ for any $\lambda\equiv 0\pmod{6}$.
\end{Lemma}

\proof Let $(V;W,U;\mathcal{B})$ be a $ \Diamond $-ITS$_\lambda (6;2,2)$, where $W=\{x_1,x_2\}$, $U=\{y_1,y_2\}$ and $V\setminus(W \cup U)=\{z_1,z_2\}$. There are $4\lambda$ triples of the form $\{x_i,y_1,z_j\}$ or $\{x_i,y_2,z_j\}$ in $\mathcal{B}$ for $1\leq i,j\leq 2$. Since $|\mathcal{B}|=13\lambda/3$ and there are only two elements in $V\setminus( U \cup W)$, one cannot find the remaining $\lambda/3$ triples, a contradiction. \qed

Since the existences of a $ \Diamond $-ITS$_{\lambda_1} (v;w,u)$ and a $ \Diamond $-ITS$_{\lambda_2} (v;w,u)$ imply the existence of a $ \Diamond $-ITS$_{\lambda_1+\lambda_2} (v;w,u)$, by Lemma \ref{necessary} and Theorem \ref{lambda=1}, to complete the proof of Theorems \ref{main result} and \ref{main result even}, it suffices to establish the existence of a $ \Diamond $-ITS$_{\lambda} (v;w,u)$ for any $v\geq 2w+u$ and
\begin{itemize}
\item $ \lambda=2 $, $vwu\equiv0\pmod2$, ${v\choose 2}-{u\choose 2}-{w\choose 2}\equiv 0\pmod3$;
\item $ \lambda=3 $, $vwu\equiv1\pmod2$, ${v\choose 2}-{u\choose 2}-{w\choose 2}\equiv 1,2\pmod3$; and
\item $ \lambda=6 $, $vwu\equiv0\pmod2$, ${v\choose 2}-{u\choose 2}-{w\choose 2}\equiv 1,2\pmod3$.
\end{itemize}

\subsection{Notation and working lemmas}

Let $G$ be a graph (or multigraph). Let $ V(G) $ and $ E(G) $ denote its vertex set and edge set, respectively. The number of edges in $ G $ is written as $\epsilon(G)$. For a subgraph $H$ of $G$, denote by $G-H$ the subgraph of $G$ obtained by deleting the edges of $H$.

For any vertex $x\in V(G)$, the degree of $ x $ is denoted by $\deg_{_G}(x)$. $G$ is said to be \textit{near-regular} if $|\deg_{_G}(x)-\deg_{_G}(y)| \leq 2 $ for all $x,y\in V(G)$.
For any $x,y\in V(G)  $ and $ x\neq y $, denote by $ E(x, y) $ the set of all edges in $ G $ which are incident both at $ x $ and $ y $. Let $\mu(G)= \max\{|E(x,y)|:x,y\in V(G)\} $ be the maximum multiplicity of $ G $.

\begin{Lemma}\label{chromatic} {\rm\cite[Lemma 9.6 and Theorem 9.12]{triple}}
There exists a near-regular multigraph of order $ s $ with $q$ edges and maximum multiplicity $3$ that has a $\triangle$-decomposition if and only if $ q\equiv0\pmod3 $ and $0\leq q\leq 3\lfloor\frac{s}{3}\lfloor \frac{3(s-1)}{2}\rfloor\rfloor$.
\end{Lemma}

Set $\Delta (G)=\max\{\deg_{_G}(x): x\in V(G)\}$. A \textit{proper edge coloring} of $G$ is a coloring of its edges satisfying that adjacent edges receive distinct colors. For any graph $ G $, let $\chi'(G)$ be the minimum number of colors required for all possible proper edge colorings of $G$.

\begin{Lemma}\label{vizing}{\rm \cite[Theorem 17.5]{bm}}
$\Delta(G)\leq \chi'(G)\leq \Delta(G)+\mu(G)$.
\end{Lemma}


A set of edges of $ G $ assigned by the same color is called a \textit{color class}. A proper edge coloring is \textit{equitable} if the sizes of any two color classes differ by at most one. The following lemma concerns about simple graphs, but its proof is also appropriate for multigraphs.

\begin{Lemma}{\rm \cite[Lemma 2.3]{h17}} \label{equi}
If a $($multi-$)$ graph has a proper edge coloring with $t$ colors, then it has an equitable proper edge coloring with $t$ colors.
\end{Lemma}

A \textit{$k$-factor} of $G$ is a spanning $k$-regular subgraph of $G$, and a \textit{$k$-factorization} of $G$ partitions the edges of $G$ into disjoint $k$-factors. 
A $1$-factorization of a $k$-regular multigraph $ G $ is equivalent to an edge coloring of $ G $ with $k$ colors.

\begin{Lemma} \label{bipartite 1-factorization} {\rm\cite[Corollary 1.11]{triple}}
	A regular bipartite multigraph has a $1$-factorization.
\end{Lemma}

 A $2$-factor of a multigraph $ G $ is a collection  of cycles (allow cycles of length 2) that spans all vertices of $ G $. Clearly, for a $2$-factorable graph,  it must be $2k$-regular for some integer $k$.  And this necessary condition is also sufficient.
\begin{Lemma} \label{2-factorization} {\rm\cite[Theorem 1.22]{triple}}
Every regular  multigraph of even degree has a $2$-factorization.
\end{Lemma}

\section{Meeting the bound}\label{meeting}

By Lemma \ref{necessary}, a $ \Diamond$-ITS$_{\lambda}(v;w,u) $ exists only if $v\geq 2w+u$. In this section, we deal with the case of $v=2w+u$.

\begin{Lemma}\label{2w+u,2}
A $\Diamond$-ITS$_2(2w+u-2z;w,u;z)$ with $w\geq u >z \geq 0$ exists whenever $u\geq 2z+1$ and $(w-z)(u-2z-1) \equiv 0 \pmod{3}$.
\end{Lemma}

\proof By Theorem \ref{ITS-one hole}, there is a $\Diamond$-ITS$_{2}(2w-u;w-u+z,0)$, $(V;L,\emptyset;{\mathcal B})$, for any $w\geq u >z \geq 0$, $u\geq 2z+1$ and $(w-z)(u-2z-1) \equiv 0 \pmod{3}$.
Let $L=\{l_1,l_2,\ldots,l_{w-u+z}\}$ and for each $1\leq i\leq w-u+z$, form a multiset of pairs ${\cal P}_i=\{\{x,y\}:\{l_i,x,y\}\in{\mathcal B}\}$. Let ${\mathcal B}_1=\{B\in {\mathcal B}:B\cap L=\emptyset\}$. Let $V_1=V\setminus L$. Define sets $W=\{b_1,b_2,\ldots,b_{w-z}\}$, $U=\{c_1,c_2,\ldots,c_{u-z}\}$ and $Z=\{d_1,d_2,\ldots,d_{z}\}$. We shall construct the required $\Diamond$-ITS on $V_1\cup W\cup U\cup Z$ with two holes $W\cup Z$ and $U\cup Z$.

Define ${\mathcal B}_2=\{\{d_i,x,y\}:\{x,y\}\in {\cal P}_i,1\leq i\leq z\}$. Form a multigraph $G_1=(V_1,E_1)$ where $E_1=\bigcup_{i=z+1}^{w-u+z}{\cal P}_i$. Then $G_1$ is a $2(w-u)$-regular multigraph. Form a complete bipartite 2-multigraph $G_2=(V_1\cup U,E_2)$ with $V_1$ and $U$ as its two partite sets; the degrees of each point in $V_1$ and $U$ are $2(u-z)$ and $2(w-z)$, respectively. Now take the union of $G_1$ and $G_2$ to form a multigraph $G=(V_1\cup U,E_1\cup E_2)$. $G$ is a $2(w-z)$-regular multigraph. By Lemma \ref{2-factorization}, $G$ admits a $2$-factorization ${\cal F}_1,{\cal F}_2,\ldots,{\cal F}_{w-z}$. Define ${\mathcal B}_3=\{\{b_j,x,y\}:\{x,y\}\in {\cal F}_j,1\leq j\leq w-z\}$. Then $(V_1\cup W \cup U\cup Z;W\cup Z,U\cup Z;{\mathcal B}_1 \cup {\mathcal B}_2\cup {\mathcal B}_3)$ is a $\Diamond$-ITS$_2(2w+u-2z;w,u;z)$. \qed

A $\Diamond$-ITS$_\lambda(v;0,0)$ is often called a {\em triple system} of order $v$ with index $\lambda$, denoted by a TS$_{\lambda}(v)$. An {\em $ \alpha $-parallel class} of a TS$_{\lambda}(v)$ is a collection of its blocks containing every point of this design exactly $ \alpha $ times. A TS$_{\lambda}(v)$ is {\em $ \alpha $-resolvable} if the collection of its blocks can be partitioned into $ \alpha $-parallel classes.

\begin{Lemma}{\rm\cite{jmv}}\label{3resolvable}
A $3$-resolvable TS$_{\lambda}(v)$ exists if and only if $ \lambda(v-1)\equiv0\pmod6 $.
\end{Lemma}

\begin{Lemma}\label{2w+u,3}
	A $\Diamond$-ITS$_3(2w+u;w,u)$ with $w\geq u\geq 1$ exists whenever $w$ and $u$ are odd.
\end{Lemma}

\proof The case of $w=1$ is trivial. For any odd $w\geq 3$, take a $3$-resolvable TS$_3(w)$ from Lemma \ref{3resolvable}, $(V,{\cal B})$, with $3$-parallel classes ${\cal P}_1,{\cal P}_2,\ldots,{\cal P}_{(w-1)/2}$. Let  $V=\{a_1,a_2,\ldots,a_w\}$, $W=\{b_1,b_2,\ldots,b_w\}$, and $U=\{c_1,c_2,\ldots,c_u\}$.
Define
\begin{center}\tabcolsep 0.01in
	\begin{tabular}{lll}
		${\cal B}_1$ & = & $\bigcup\limits_{l=1}^{(u-1)/2}{\cal P}_l$;\\
		${\cal B}_2$ & = & $\{\{a_i,a_j,b_k\}, \{a_i,b_j,a_k\}, \{b_i,a_j,a_k\}:\{a_i,a_j,a_k\}\in{\cal P}_l,(u+1)/2\leq l \leq (w-1)/2\}$;\\
		${\cal B}_3$ & = & $\{\{c_u,a_i,b_i\},\{c_u,a_i,b_i\},\{c_u,a_i,b_i\}:1\leq i \leq w\}$;\\
		${\cal B}_4$ & = & $\{\{c_{2l},a_i,b_j\}, \{c_{2l},a_j,b_k\}, \{c_{2l},a_k,b_i\}, \{c_{2l-1},a_i,b_k\}, \{c_{2l-1},a_j,b_i\},\{c_{2l-1},a_k,b_j\}:$ \\
		&& $\{a_i,a_j,a_k\}\in{\cal P}_l,1\leq l \leq (u-1)/2\}$.\\
	\end{tabular}
\end{center}
Then $(V\cup W \cup U;W,U;{\cal B}_1 \cup {\cal B}_2 \cup {\cal B}_3 \cup {\cal B}_4)$ is a $\Diamond$-ITS$_3(2w+u;w,u)$. \qed

\begin{Lemma}\label{2w+u,6}
A $\Diamond$-ITS$_6(2w+u;w,u)$ exists for any $w\geq u\geq 1$ except for $w=u=2$.
\end{Lemma}

\proof The case of $w=1$ is trivial. When $w=2$ and $u=1$, the conclusion follows from Theorem \ref{ITS-one hole}. When $w=u=2$, there is no $\Diamond$-ITS$_6(6;2,2)$ by Lemma \ref{w=u=2}.

For any $w\geq 3$, take a $3$-resolvable TS$_6(w)$ from Lemma \ref{3resolvable}, $(V,{\cal B})$, with $3$-parallel classes ${\cal P}_1,{\cal P}_2,\ldots,{\cal P}_{w-1}$. Let  $V=\{a_1,a_2,\ldots,a_w\}$, $W=\{b_1,b_2,\ldots,b_w\}$, and $U=\{c_1,c_2,\ldots,c_u\}$. Define
\begin{center}\tabcolsep 0.01in
\begin{tabular}{lll}
${\cal B}_1$ & = & $\bigcup\limits_{l=1}^{u-1}{\cal P}_l$;\\
${\cal B}_2$ & = & $\{\{a_i,a_j,b_k\}, \{a_i,b_j,a_k\}, \{b_i,a_j,a_k\}:\{a_i,a_j,a_k\}\in{\cal P}_l,u\leq l \leq w-1\}$;\\
${\cal B}_3$ & = & $\{\{c_1,a_i,b_i\},\{c_1,a_i,b_i\},\{c_1,a_i,b_i\}, \{c_u,a_i,b_i\},\{c_u,a_i,b_i\},\{c_u,a_i,b_i\}:1\leq i \leq w\}$;\\
${\cal B}_4$ & = & $\{\{c_l,a_i,b_j\}, \{c_l,a_j,b_k\}, \{c_l,a_k,b_i\}, \{c_{l+1},a_i,b_k\}, \{c_{l+1},a_j,b_i\},\{c_{l+1},a_k,b_j\}:$ \\
&& $\{a_i,a_j,a_k\}\in{\cal P}_l,1\leq l \leq u-1\}$.\\
\end{tabular}
\end{center}
Then $(V\cup W \cup U;W,U;{\cal B}_1 \cup {\cal B}_2 \cup {\cal B}_3 \cup {\cal B}_4)$ is a $\Diamond$-ITS$_6(2w+u;w,u)$. \qed

\section{Near the bound}\label{near bound}

In this section, we deal with the cases of $ v=2w+u+x $ with $ x\in \{1,2,3,4\} $. We begin by providing the definitions of  {\rm IGDDs}.

An {\em incomplete group divisible design} with block size $ k $ and index $ \lambda $, briefly a $ (k, \lambda) $-IGDD, is a quintuple $ (V, \mathcal{G}, H, H', \mathcal{B})  $ where $V$ is a $v$-set, $\mathcal{G}$ is a partition of $V$ into {\em groups}, $H$ and $H'$ are two disjoint subsets of $V$ (called {\em holes}), and $ \mathcal{B} $ is a family of $ k $-subsets of $ V $ (called {\em blocks}) satisfying that each block meets each group and each hole in at most one element, and any two elements of $ V $ from different groups and different holes appear in exactly $\lambda$ blocks. An IGDD is of {\em type} $ (g_1; h_1,h'_1)^{a_1} (g_2; h_2,h'_2)^{a_2} \cdots (g_t; h_t,h'_t)^{a_t }$ if there are $ a_i $ groups of size $ g_i $, $1\leq i\leq t$, each of which intersects $H$ and $H'$ in exactly $h_i$ and $h'_i$ elements, respectively.

The following lemma was stated in \cite{cor} by using the terminology of Latin squares.

\begin{Lemma}\label{IGDD small}{\rm\cite[Theorem 1.4]{cor}}
Let $ n$, $m $ and $ l $ be positive integers with $m\geq  l $. Then there exists a $(3,1)$-IGDD of type $(n;m,l)^3$ for any $ n\geq 2m+l $.
\end{Lemma}

A $ (k, \lambda) $-IGDD $ (V, \mathcal{G}, H, H', \mathcal{B})  $ is called a \textit{group divisible design} if $ H= H'=\emptyset $, and is denoted by a $ (k, \lambda) $-GDD.

\begin{Lemma}{\rm\cite{chang,wy}}\label{GDD2}
There exists a $ (3, \lambda) $-GDD of type $ g^t w^1 $ if and only if all of the following conditions are satisfied:
	$ (1) $ if $ g > 0 $, then $ t \geq 3 $, or $ t = 2  $ and $ w = g $, or $ t = 1  $ and $ w = 0 $, or $ t = 0 $;
    $ (2) $ $ w \leq g(t-1) $ or $ gt = 0 $;
	$ (3) $ $ \lambda g(t-1) + \lambda w \equiv 0 \pmod 2  $ or $ gt = 0 $;
	$ (4) $ $ \lambda gt \equiv 0 \pmod 2 $ or $ w = 0 $;
	$ (5) $ $ \lambda g^2t(t-1)/2 + \lambda gtw \equiv 0 \pmod 3 $.
\end{Lemma}

 We can construct some $\Diamond$-ITS$_\lambda(v;u,u)$s by employing Lemma \ref{GDD2}.

\begin{Lemma}\label{w=u small}
If $\lambda$ is even, $ \lambda({v\choose 2}-2{u\choose 2})\equiv 0 \pmod 3$ and  $3u\leq v\leq 5u$, then there exists a $\Diamond$-ITS$_\lambda(v;u,u)$ for any $ u\geq3 $.
\end{Lemma}

\proof  By Lemma \ref{GDD2}, when $0<3u\leq v\leq 5u$, a $ (3,\lambda) $-GDD of type $ u^3(v-3u)^1$ exists for any even $\lambda$. Place a TS$_\lambda(u)  $ (from Theorem \ref{ITS-one hole}) on one group of size $ u $ and a TS$ _\lambda(v-3u) $ (from Theorem \ref{ITS-one hole}) on the group of size $ v-3u $.
Note that when $\lambda\not\equiv 0\pmod{6}$, $ \lambda({v\choose 2}-2{u\choose 2})\equiv 0 \pmod 3$ implies $ u,v\equiv0,1\pmod3  $. \qed

To construct more ITSs,  we need to introduce a new type of incomplete triple systems with three holes. A $ \Diamond $-ITS$_\lambda (v;w,u,y;z_1,z_2) $ is a quintuple $ (V;W,U,Y;\mathcal{B}) $, where $ V $ is a $v$-set,  $W, U, Y\subseteq V$ (called {\em holes}) with $|W|=w$, $|U|=u$ and $|Y|=y$ for which $W\cap U=\emptyset  $, $ |W\cap Y|=z_1 $, $ |U\cap Y|=z_2 $, and $\mathcal{B}$ is a collection of $3$-subsets (called {\em triples}) of $V$ in which every pair of distinct elements from $V$ not contained within $U$, $W$ or $ Y $ appears in exactly $ \lambda $ triples in $\mathcal{B}$, while pairs within $U$, $W$ or $ Y $ appear in no triple.

\begin{Remark}\label{rek:lattice}
A $ \Diamond $-ITS$_\lambda (3n+c;n+a,n+b,c;a,b) $ is also a $ \Diamond $-ITS$_\lambda (3n+c;n+b,n+a,c;b,a) $.
\end{Remark}

\begin{Lemma}\label{cons lattice}
If there exist a $(3,\lambda)$-GDD of type $n^3c^1$ and a TS$_\lambda (n) $, then there exists a $ \Diamond $-ITS$_\lambda (3n+c;n,n,c;0,0) $. 
\end{Lemma}
\proof Start with a $(3,\lambda)$-GDD of type $n^3c^1$. Then place a TS$_\lambda (n) $ on one group of size $ n $. \qed

\begin{Lemma}\label{lattice-small}
There exists a $ \Diamond $-ITS$_2 (v;w,u,y;z_1,z_2 ) $
for $ (v,w,u,y,z_1,z_2)\in\{(11,3,3,2,0,0)$, $(11,4,3,2,1,0)$, $(12,3,3,3,0,0)$, $(13,3,3,4,0,0)$, $(13,4,3,4,1,0)$, $(14,$ $5, 4,5,2,3)$, $(18,7,4,5,2,3)$, $(22,6,6,4,0,0)$, $(23,8,7,5,2,3)$, $(23,8,7,8,2,6)$, $(24,8,6,$ $6, 2,0)\}$.
\end{Lemma}

\proof By taking $ (n,c,\lambda)\in\{(3,2,2),(3,3,2),(3,4,2),(6,4,2)\}$, apply Lemma \ref{cons lattice} together with Theorem \ref{ITS-one hole} and Lemma \ref{GDD2} to get the cases of $ (v,w,u,y,z_1,z_2)\in\{(11,3,3,2,0,0)$, $(12,3,3,3,0,\\0),(13,3,3,4,0,0)$, $(22,6,6,4,0,0)\}$. For the other cases, see Appendix A. \qed

\begin{Construction}\label{con:cons1}
Let $a,b,c\geq0$ with $c\geq a+b$. Suppose that there exist
\begin{itemize}
\item[$(1)$] a $(3,\lambda)$-IGDD of type $(2m+l;m,l)^1(3h;h,h)^{t-1}(2h;h,0)^{s-t}$,
\item[$(2)$] a $\Diamond$-ITS$_\lambda(2m+l+c;m+a,l+b)$,
\item[$(3)$] a $\Diamond$-ITS$_\lambda(2h+c;h+a,c;a)$ when $s>t$, and
\item[$(4)$] a $ \Diamond $-ITS$_\lambda (3h+c;h+a,h+b,c;a,b)$ when $ t>1 $.
\end{itemize}
Then there exists a $\Diamond$-ITS$_\lambda(2sh+ht-3h+2m+l+c;hs-h+m+a,ht-h+l+b)$.
\end{Construction}

\proof Let $ (V,\mathcal{G}, H,  H', \mathcal{B})$ be a $(3,\lambda)$-IGDD of type $(2m+l;m,l)^1(3h;h,h)^{t-1}(2h;h,0)^{s-t}$, where $\mathcal{G}=\{G_1,G_2,\ldots,G_s\}$ satisfying that $|G_1|=2m+l$, $|G_i|=3h$ for $2\leq i\leq t$, and $|G_i|=2h$ for $t+1\leq i\leq s$. Let $\{\infty_1,\infty_2,\ldots,\infty_c\}\cap V=\emptyset$. We shall construct the required $\Diamond$-ITS on $ V\cup \{\infty_1,\ldots,\infty_c\} $ with two holes $ H\cup \{\infty_1,\ldots,\infty_a\} $ and $ H'\cup \{\infty_{a+1},\ldots,\infty_{a+b}\} $. For $ 1\leq i\leq s $, let $H_i=G_i\cap H  $ and $H'_i=G_i\cap H'  $. Note that $H'_i=\emptyset$ for $t+1\leq i\leq s$.

Firstly, place on $G_1\cup \{\infty_1,\ldots,\infty_c\}$ a $\Diamond$-ITS$_\lambda(2m+l+c;m+a,l+b)$, aligning the hole of size $m+a$ on $H_1\cup \{\infty_1,\ldots,\infty_a\}$ and the hole of size $l+b$ on $H'_1\cup \{\infty_{a+1},\ldots,\infty_{a+b}\}$. Then place a $ \Diamond $-ITS$_\lambda (3h+c;h+a,h+b,c;a,b)$ on $G_i\cup \{\infty_1,\ldots,\infty_c\}$ for each $2\leq i\leq t$, aligning the third hole on the infinite points, the first hole on $H_i\cup \{\infty_1,\ldots,\infty_a\}$ and  the second hole on $H'_i\cup\{\infty_{a+1},\ldots,\infty_{a+b}\}$. Finally, place a $\Diamond$-ITS$_\lambda(2h+c;h+a,c;a)$ on $G_i\cup \{\infty_1,\ldots,\infty_c\}$ for each $t<i\leq s$, aligning the hole of size $c$ on the infinite points, the hole of size $h+a$ on $H_i\cup \{\infty_1,\ldots,\infty_a\}$, and hence the common points on $\{\infty_1,\ldots,\infty_a\}$. \qed

\subsection{Constructions for IGDDs}

By Construction \ref{con:cons1}, various IGDDs with different types are extremely needed. This subsection is devoted to constructing IGDDs.

Let $ \alpha $ be a positive integer. Given a $ (k,  \lambda) $-GDD on the set $ V $ with the group set $ \mathcal{G}=\{G_1,G_2,\ldots,G_t\} $, an \textit{$ \alpha $-partial parallel class} associated with group $ G_i $, $ 1\leq i\leq t $, is a subcollection of blocks such that every element of $ V\setminus G_i $ occurs exactly $ \alpha $ times and the elements of $ G_i $ do not occur at all.

A $ (k,  \lambda) $-GDD is said to be a \textit{$ (k,  \lambda) $-frame} if its block set can be partitioned into $ 1 $-partial parallel classes. In a $(k,\lambda)$-frame of type $g^n$, each 1-partial parallel class contains $g(n-1)/k$ blocks, and every block belongs to a 1-partial parallel class associated with a group that is disjoint from the block. It is known that for any group $G$ of a $(k,\lambda)$-frame of type $g^n$, the number of 1-partial parallel classes associated with $G$ is $\lambda g/(k-1)$ (see \cite[Lemma 2.2.2]{fmy}).

\begin{Lemma}\label{3-frame} {\rm \cite{gm}}
There exists a $(3,\lambda)$-frame of type $g^n$ if and only if $n\geq 4$, $\lambda g \equiv 0 \pmod{2}$ and $g(n-1)\equiv 0 \pmod{3}$.
\end{Lemma}

Colbourn, Oravas and Rees \cite{cor} introduced a special class of frames, named eframes, to construct IGDDs. Let $m$ be even. A $(3,1)$-GDD of type $6^n m^1$ is said to be a \textit{$(3,1)$-eframe} if its block set can be partitioned into $n$ $3$-partial parallel classes, each associated with one of the groups of size $6$, and $\frac{m}{2}$ $1$-partial parallel classes associated with the group of size $m$.

\begin{Lemma}\label{eframe-cor}{\rm \cite[Corollary 4.10]{cor}}
Let $ n\geq 4 $ and $ n\notin\{5,6,8,10,11,14,17\} $. Then there exist $(3,1)$-eframes of types $6^n 2^1$ and $6^n 4^1$.
\end{Lemma}

Now we shall generalize the definition of eframes. Let $\lambda$ be even. A $(3,\lambda)$-GDD of type $h^nm^1$ is said to be a \textit{$(3,\lambda)$-eframe} if its block set can be partitioned into $\frac{\lambda n}{2}$ $h$-partial parallel classes, each associated with one of the groups of size $h$, and $\frac{\lambda m}{2}$ $1$-partial parallel classes associated with the group of size $m$. By the same argument as that in \cite[Lemma 2.2.2]{fmy}, one can see that for any group $G$ of size $h$ in a $(3,\lambda)$-eframe of type $h^nm^1$, the number of $h$-partial parallel classes associated with $G$ is $\frac{\lambda}{2}$.

\begin{Remark}\label{rek:eframe}
$\lambda$ copies of a $(3,1)$-eframe of type $6^n m^1$ produce a $(3,\lambda)$-eframes of type $6^n m^1$.
\end{Remark}

\begin{Lemma}\label{lem:eframe-direct}
There exists a $(3,2)$-eframe of type $h^n m^1$ for $(h,n,m)\in\{(3,5,1),(3,5,2),(3, 6,1)$, $(3,6,2),(6,5,5),(6,5,8)\}$.
\end{Lemma}

\proof We give an explicit construction for each eframe in Appendix B. \qed

\begin{Construction}\label{cons:IGDD}
Let $ \lambda\equiv0\pmod2 $.  If there exists a $(3,\lambda)$-eframe of type $h^{s-1}m^1$, then there exists a $(3,\lambda)$-IGDD of type $(2m+l;m,l)^1(3h;h,h)^{t-1}(2h;h,0)^{s-t}$ for any  $ 0\leq l\leq m $ and $1\leq t\leq s$.
\end{Construction}

\proof Take $2s+t$ pairwise disjoint sets $J_1,J_2,\ldots,J_s$, $H_1,H_2,\ldots,H_s$ and $H'_1,H'_2,\ldots,H'_t$ such that $|J_1|=|H_1|=m$, $|H'_1|=l$, $|J_i|=|H_i|=h$ for $2\leq i\leq s$, and $|H'_i|=h$ for $2\leq i\leq t$. We shall construct the required IGDD on $(\bigcup_{i=1}^{s} (J_i\cup H_i))\cup(\bigcup_{i=1}^t H'_i)$ with groups $J_i\cup H_i\cup H'_i$ for $1\leq i\leq t$, and $J_i\cup H_i$ for $t+1\leq i\leq s$, where holes are $\bigcup^s_{i=1}H_i$ and $\bigcup^t_{i=1}H'_i$.

First construct a $(3,\lambda)$-eframe of type $h^{s-1}m^1$ on $\bigcup_{i=1}^s J_i$ with groups $J_1,J_2,\ldots,J_s$. Let ${\cal P}_{11},{\cal P}_{12},\ldots,{\cal P}_{1\frac{\lambda m}{2}}$ be 1-partial parallel classes missing the group $J_1$. For $2\leq i\leq s$, let ${\cal P}_{i1},{\cal P}_{i2},\ldots,{\cal P}_{i\frac{\lambda}{2}}$ be $h$-partial parallel classes missing the group $J_i$.

Now we form blocks of the required IGDDs. Set
$${\cal B}_1=\left(\bigcup_{j=1}^{\frac{\lambda l}{2}}{\cal P}_{1j}\right)\cup\left(\bigcup_{i=2}^t\bigcup_{j=1}^{\frac{\lambda}{2}}{\cal P}_{ij}\right).$$
Let $\phi$ be a bijection from $\bigcup_{i=1}^s J_i$ to $\bigcup_{i=1}^s H_i$ which maps each element in $J_i$ to some element in $H_i$ for each $1\leq i\leq s$. Set
\begin{center}
	${\cal B}_2=\bigcup\limits_{\substack{(i,j)\in I\\ \{x,y,z\}\in{\cal P}_{ij}}}\{\{\phi(x),y,z\},\{x,\phi(y),z\},\{x,y,\phi(z)\}\},$
\end{center}
where $I=\left\{(1,j):\lambda l/2<j\leq \lambda m/2\right\} \cup\left\{(i,j):t+1\leq i\leq s,1\leq j\leq \lambda /2\right\}$. For $1\leq i\leq t$, define a multigraph $G_i$ with vertex set $\left(\bigcup_{r=1}^s(J_r\cup H_r)\right)\backslash(J_i\cup H_i) $ and edge set
\begin{center}
	$E(G_1)=\bigcup\limits_{\substack{1\leq j\leq \frac{\lambda l}{2}\\ \{x,y,z\}\in{\cal P}_{1j}}}\{\{x,\phi(y)\},\{x,\phi(z)\},\{y,\phi(x)\},\{y,\phi(z)\},\{z,\phi(x)\},\{z,\phi(y)\}\}$
\end{center}
and for $2\leq i\leq t$,
\begin{center}
	$E(G_i)=\bigcup\limits_{\substack{1\leq j\leq \frac{\lambda}{2}\\ \{x,y,z\}\in{\cal P}_{ij}}}
	\{\{x,\phi(y)\},\{x,\phi(z)\},\{y,\phi(x)\},\{y,\phi(z)\},\{z,\phi(x)\},\{z,\phi(y)\}\}. $
\end{center}
Then for $1\leq i\leq t$, $G_i$ is a bipartite multigraph with bipartition $\{(\bigcup_{r=1}^s J_r)\setminus J_i,(\bigcup_{r=1}^s H_r)\setminus H_i)\}$, and $G_i$ is regular of degree $\lambda|H'_i|$. By Lemma \ref{bipartite 1-factorization}, $G_i$ admits a 1-factorization into $\lambda|H'_i|$ 1-factors ${\cal F}_{\alpha\beta}^{(i)}$, $1\leq\alpha\leq |H'_i|$, $1\leq\beta\leq\lambda$. Let $H'_i=\{c_{i1},c_{i2},\ldots,c_{i|H'_i|}\}$. Set
\begin{center}
	${\cal B}_3=\bigcup\limits_{i=1}^{t}\bigcup\limits_{\alpha=1}^{|H'_i|}\bigcup\limits_{\beta=1}^{\lambda}\bigcup
	\limits_{\{x,y\}\in{\cal F}_{\alpha\beta}^{(i)}}\{\{c_{i\alpha},x,y\}\}.$	
\end{center}
It is readily checked that ${\cal B}_1 \cup {\cal B}_2 \cup {\cal B}_3$ forms the block set of a $(3,\lambda)$-IGDD of type $(2m+l;m,l)^1(3h;h,h)^{t-1}(2h;h,0)^{s-t}$. \qed

\begin{Construction}\label{special,GDD}
 If there exists a $(3,6)$-GDD of type $ 2^{s-1}1^1 $ that has $1$ $3$-partial parallel class associated with the group of size $1$, and $1$ $6$-partial parallel class for each group of size $2$, then there exists a $(3,6)$-IGDD of type $(2+l;1,l)^1(6;2,2)^{t-1}$ $(4;2,0)^{s-t}$ for $l\in\{0,1\}$ and $1\leq t\leq s$.

\end{Construction}

\proof (sketch only) The proof is very similar to that in Construction \ref{cons:IGDD}. First define $ 2s+t$ pairwise disjoint sets and a bijection $ \phi $ as in Construction \ref{cons:IGDD}. The key step is to combine the vertices in $ H_i' $ with the vertices in  $ (\bigcup_{r=1}^s J_r)\setminus J_i$ and $(\bigcup_{r=1}^s H_r)\setminus H_i) $. To realize it, a bipartite $ \lambda|H_i | $-regular multigraph is defined. In what follows, we just give the required bipartite multigraphs.

 Let ${\cal P}_{1}$ be the $3$-partial parallel class missing the group $J_1$ and ${\cal P}_{i}$ be the $6$-partial parallel class missing the group $J_i$ for $ 2\leq i\leq s $. If $ |H_i'| \neq 0$, we define bipartite $ 6 |H_i'| $-regular multigraphs $G_i$ with edge set
\begin{center}
	$E(G_i)=\bigcup\limits_{ \{x,y,z\}\in{\cal P}_{i}}
	\{\{x,\phi(y)\},\{x,\phi(z)\},\{y,\phi(x)\},\{y,\phi(z)\},\{z,\phi(x)\},\{z,\phi(y)\}\}. $
\end{center}
If $ |H_i'|  =0$, then $G_i$ is taken as an empty graph. \qed

\begin{Lemma}\label{GDDs,partial}
	There exists a $(3,6)$-GDD of type $ 2^{s-1}1^1 $ for $ s\in\{5,6\}$ that has $1$ $3$-partial parallel class associated with the group of size $1$, and $1$ $6$-partial parallel class for each group of size $2$.
\end{Lemma}

\proof We give an explicit construction for each GDD in Appendix C. \qed
\begin{Theorem}\label{special,IGDD}
Let $s\in\{5,6\}$. Then there exists a $(3,6)$-IGDD of type $(2+l;1,l)^1(6;2,2)^{t-1}$ $(4;2,0)^{s-t}$ for  $l\in\{0,1\}$ and $1\leq t\leq s$.
\end{Theorem}

\proof Take a  $(3,6)$-GDD of type $ 2^{s-1}1^1 $ for $ s\in\{5,6\} $ that has 1 3-partial parallel class associated with the group of size 1, and 1 6-partial parallel class for each group of size 2 (from Lemma \ref{GDDs,partial}). Then apply Construction \ref{special,GDD}. \qed

\begin{Theorem}\label{IGDD}
Let $s\geq4$, $\lambda \equiv 0 \pmod{2}$ and $h(s-1)\equiv 0 \pmod{3}$. Then there exists a $(3,\lambda)$-IGDD of type $(2h+l;h,l)^1(3h;h,h)^{t-1}(2h;h,0)^{s-t}$ for any $0\leq l\leq h$ and $1\leq t\leq s$.
\end{Theorem}

\proof Clearly, a $ (3,\lambda) $-frame of type $ h^s $ with even $\lambda$ yields a $ (3,\lambda) $-eframe of type $ h^{s-1}h^1 $. Then apply Construction \ref{cons:IGDD} together with the existence of $ (3,\lambda) $-frames of type $ h^s $ from Lemma \ref{3-frame}, we get the desired IGDDs. \qed

\begin{Theorem}\label{IGDD3}
Let $s\geq5$ and $ s\notin\{6,7,9,11,12,15,18\} $. Let $ m\in\{2,4\} $. Then there exists a $(3,2)$-IGDD of type $(2m+l;m,l)^1(18;6,6)^{t-1}(12;6,0)^{s-t}$ for any $ 0\leq l
	\leq m$  and $1\leq t\leq s$.
\end{Theorem}

\proof Take $(3,\lambda)$-eframes of types $6^{s-1} 2^1$ and $6^{s-1} 4^1$ (from Lemma \ref{eframe-cor} and Remark \ref{rek:eframe}). Then apply Construction \ref{cons:IGDD} to complete the proof. \qed

\begin{Theorem}\label{IGDD-add}
Let $s\in\{6,10\}$. Then there exists a $(3,2)$-IGDD of type $(16+l;8,l)^1(18;6,6)^{t-1}$ $(12;6,0)^{s-t}$ for any $0\leq l\leq 8$ and $1\leq t\leq s$.
\end{Theorem}

\proof By Lemma \ref{lem:eframe-direct}, there exists a $(3,2)$-eframe of type $6^5 8^1$. By \cite[Theorem 4.9]{cor}, there exists a $(3,1)$-eframe of type  $6^9 8^1$, which yields a $(3,2)$-eframe of type $6^9 8^1$ by Remark \ref{rek:eframe}.  Start from these eframes and then apply Construction \ref{cons:IGDD} to complete the proof. \qed

A $(k,1)$-GDD of type $g^k$ is said to be a {\em transversal design}, and denoted by a TD$(k,n)$, which is equivalent to $k-2$ mutually orthogonal Latin squares of order $n$.

\begin{Lemma}\label{eframe}
Let $n\equiv 0,1\pmod 3$, $n\geq 4$ and $n\not\in \{6,10\}$. Let $0\leq m< n $ when $n\equiv 1\pmod 3$ and $1\leq m\leq n  $ when $n\equiv 0\pmod 3$. Then there exists a $(3,2)$-eframe of type $3^n m^1$.
\end{Lemma}

\proof It is known that a TD$(5,n)$ exists for any $n\geq 4$ and $n\not\in \{6,10\}$ \cite{ACD}. Construct a TD$(5,n)$ on $V=\{0,1,\ldots,n-1\}\times\{1,2,3,4,5\} $ with groups $ \{0,1,\ldots,n-1\}\times\{i\} $ for $ i\in\{1,2,3,4,5\} $. Denote by $ \mathcal{B} $ its block set.

Let $ \{\infty_1,\infty_2,\ldots,\infty_m\}$ be a set that is disjoint from $V$. We shall construct a $(3,2)$-eframe of type $3^n m^1$ on $ (\{0,1,\ldots,n-1\}\times\{1,2,3\} ) \cup\{\infty_1,\infty_2,\ldots,\infty_m\},$ where $\{j\}\times\{1,2,3\} $ forms a group of size $ 3 $ for $0\leq j<n $, and $\{\infty_1,\infty_2,\ldots,\infty_m\}  $ forms the group of size $ m $.

For convenience, write the element $(j,i)$ in $V$ as $ j_i $. Without loss of generality, assume that $ \{j_1,j_2,j_3,j_4,0_5\}\in \mathcal{B}$ for $ 0\leq j<n $. Set
$$\mathcal{H}_j=\{\{a_1,b_2,c_3\}:\{a_1,b_2,c_3,j_4,l_5\}\in\mathcal{B},1\leq l<n\}  $$
for $ 0\leq j<n $. Then ${\cal H}_j$ is a 1-partial parallel class with respect to the ``group" $ \{j_1,j_2,j_3\} $. Further define
$$\mathcal{P}_l=\{\{a_1,b_2,c_3\}:\{a_1,b_2,c_3,j_4,l_5\}\in\mathcal{B},0\leq j<n\}  $$
for $ 1\leq l<n $. Then ${\cal P}_l$ covers each element in $ \{0,1,\ldots,n-1\}\times\{1,2,3\}$ exactly once. One can check that $ \mathcal{H}_j\cap \mathcal{P}_l$ contains exactly one triple for $ 0\leq j<n $ and $ 1\leq l<n $.

Write $ \beta=m $ if $ n\equiv1\pmod3 $ and $ \beta=m-1 $ if $ n\equiv0\pmod3 $.

If $ n\equiv1\pmod3 $, for any $ i\in\{1,2,3\} $, construct a $ (3,2) $-frame of type $ 1^n $ (from Lemma \ref{3-frame}) on $\{0,1,\ldots,n-1\}\times\{i\} $ with groups $\{j_i\}$, $0\leq j<n$. Denote by $ \mathcal{F}_{ji} $ its 1-partial parallel class with respect to the group $\{j_i\}$. Let
$$\mathcal{C}_{j1} = \mathcal{F}_{j1}\cup\mathcal{F}_{j2}\cup\mathcal{F}_{j3}. $$
For $ 0\leq j<n $, let
$$\mathcal{C}_{j2} =\bigcup\limits_{\substack{1\leq l\leq\beta\\ \{a_1,b_2,c_3\}\in{{\cal H}_j}\cap{\cal P}_l}}\{\{a_1,b_2,\infty_l\},\{a_1,c_3,\infty_l\},\{b_2,c_3,\infty_l\}\}$$
and
$$\mathcal{C}_{j3} =\bigcup\limits_{\substack{\beta<l<n \\\{a_1,b_2,c_3\}\in{{\cal H}_j}\cap{\cal P}_l}}\{\{a_1,b_2,c_3\},\{a_1,b_2,c_3\}\}.$$
Take $ \mathcal{C}_j=\mathcal{C}_{j1}\cup\mathcal{C}_{j2}\cup\mathcal{C}_{j3} $. Then $\mathcal{C}_j$ is a
3-partial parallel class with respect to the group $\{j_1,j_2,j_3\}$. For $ 1\leq l\leq \beta $, let
$$ 	\mathcal{D}_l =\bigcup\limits_{\substack{0\leq j<n \\ \{a_1,b_2,c_3\}\in{{\cal H}_j}\cap{\cal P}_l}}\{\{a_1,b_2,c_3\}\} .$$
Then $\mathcal{D}_l=\mathcal{P}_l$ is a 1-partial parallel class with respect to the group $\{\infty_1,\infty_2,\ldots,\infty_m\}  $. One can check that $(\bigcup_{j=0}^{n-1}{\cal C}_j)\cup(\bigcup_{l=1}^{\beta}{\cal D}_l)$ forms a $(3,2)$-eframe of type $3^n m^1$.

If $ n\equiv0\pmod3 $, for any $ i\in\{1,2,3\} $, construct a $ (3,2) $-frame of type $ 1^{n+1} $ on $(\{0,1,\ldots,n-1\}\times\{i\})\cup \{\infty_m\} $ with groups $\{j_i\}$, $0\leq j<n$, and $\{\infty_m\}$. Denote by $ \mathcal{G}_{ji} $ and $ \mathcal{G}_{\infty i} $ its 1-partial parallel class with respect to the group $\{j_i\}$ and the group $ \{\infty_m\}$, respectively. Let
$$\mathcal{C}'_{j1} = \mathcal{G}_{j1}\cup\mathcal{G}_{j2}\cup\mathcal{G}_{j3}.$$
Take $ \mathcal{C}'_j=\mathcal{C}'_{j1}\cup\mathcal{C}_{j2}\cup\mathcal{C}_{j3} $ (note that $\beta=m-1$ now). Then $\mathcal{C}'_j$ is a 3-partial parallel class with respect to the group $\{j_1,j_2,j_3\}$.
For $ 1\leq l\leq \beta $, let ${\cal D}'_l={\cal P}_l$. Let
$$ \mathcal{D}'_m =\mathcal{G}_{\infty 1}\cup \mathcal{G}_{\infty 2}\cup\mathcal{G}_{\infty 3}. $$
Then $\mathcal{D}'_l$, $1\leq l\leq m$, are 1-partial parallel classes with respect to the group $\{\infty_1,\infty_2,\ldots,\infty_m\}$. One can check that $(\bigcup_{j=0}^{n-1}{\cal C}'_j)\cup(\bigcup_{l=1}^{m}{\cal D}'_l)$ forms a $(3,2)$-eframe of type $3^n m^1$. \qed


\begin{Lemma}\label{eframe1}
Let $n\equiv 0,1\pmod 3$, $n\geq 4$ and $n\not\in \{6,10\}$. Let $0\leq m< n  $ when $n\equiv 1\pmod 3$ and $1\leq m\leq n  $ when $n\equiv 0\pmod 3$. Then there exists a $(3,2)$-eframe of type $6^n (2m-1)^1$.
\end{Lemma}

\proof Construct a TD$(5,n)$ on $V=\{0,1,\ldots,n-1\}\times\{1,2,3,4,5\} $ with groups $ \{0,1,\ldots,n-1\}\times\{i\} $ for $ i\in\{1,2,3,4,5\} $. Then by the same argument as that in the proof of Lemma \ref{eframe}, we can define $ \mathcal{H}_j$ for $0\leq j<n$ and $\mathcal{P}_l$ for $1\leq l<n$. We still write the element $(j,i)$ in $V$ as $ j_i $.

Let $ \{\infty_1,\infty_2,\ldots,\infty_m\}$ be a set that is disjoint from $V$. We shall construct a $(3,2)$-eframe of type $6^n (2m-1)^1$ on $(((\{0,1,\ldots,n-1\}\times\{1,2,3\} ) \cup\{\infty_1,\infty_2,\ldots,\infty_{m-1}\})\times \{0,1\})\cup \{\infty_m\}$, where $\{j\}\times\{1,2,3\}\times \{0,1\} $ (or written as $\{j_1,j_2,j_3\}\times\{0,1\}$) forms a group of size $ 6 $ for $0\leq j<n $, and $(\{\infty_1,\infty_2,\ldots,\infty_{m-1}\}\times \{0,1\})\cup \{\infty_m\}$ forms the group of size $ 2m-1 $.

For $ 0\leq j<n $ and $ 1\leq l<n $, consider the triple $ \{a_1,b_2,c_3\}\in \mathcal{H}_j\cap \mathcal{P}_l $. Form a $ (3,2) $-GDD of type $ 2^3$ (from Lemma \ref{GDD2}) on $\{a_1,b_2,c_3\}\times \{0,1\} $ with groups $ \{x\}\times \{0,1\} $ for $ x\in \{a_1,b_2,c_3\} $. Denote the collection of its blocks by $ \mathcal{B}_{jl} $. Note that each element in $\{a_1,b_2,c_3\}\times \{0,1\} $ occurs in exactly four blocks in $ \mathcal{B}_{jl} $.

Write $ \beta=m $ if $ n\equiv1\pmod3 $ and $ \beta=m-1 $ if $ n\equiv0\pmod3 $.

If $ n\equiv1\pmod3 $, for any $ i\in\{1,2,3\} $, construct a $ (3,2) $-frame of type $ 2^n $ (from Lemma \ref{3-frame}) on $\{0_i,1_i,\ldots,(n-1)_i\}\times \{0,1\} $ with groups $\{j_i\}\times \{0,1\}$, $0\leq j<n$. Its two 1-partial parallel classes with respect to the group $\{j_i\}\times\{0,1\}$ are denote by $ \mathcal{F}'_{ji} $ and $ \mathcal{F}''_{ji} $. Write $\mathcal{F}_{ji}=\mathcal{F}'_{ji} \cup \mathcal{F}''_{ji}  $. Let
$$ \mathcal{C}_{j1}=\mathcal{F}_{j1}\cup \mathcal{F}_{j2}\cup\mathcal{F}_{j3}. $$
For $ 0\leq j<n $, let
\begin{align*}
&\mathcal{C}_{j2} =\bigcup\limits_{\substack{\beta<l<n \\\{a_1,b_2,c_3\}\in{{\cal H}_j}\cap{\cal P}_l}}\mathcal{B}_{jl};\\
&\mathcal{C}_{j3} =\bigcup\limits_{\substack{1\leq l<\beta\\ \{a_1,b_2,c_3\}\in{{\cal H}_j}\cap{\cal P}_l}}\bigcup\limits_{r,s\in\{0,1\}}\{\{(x,r),(y,1-r),(\infty_l,s)\}:
(x,y)\in\{(a_1,b_2),(a_1,c_3),(b_2,c_3)\}\};\\
&\mathcal{C}_{j4} = \bigcup\limits_{\{a_1,b_2,c_3\}\in{{\cal H}_j}\cap{\cal P}_\beta}\bigcup\limits_{r\in\{0,1\}}\{\{(a_1,r),(b_2,1-r),(c_3,r)\},\{(a_1,r),(b_2,r),(c_3,1-r)\},\\
&~~~~~~~~~~~~~~~~~~~~~~\{(a_1,r),(b_2,1-r),\infty_m\},\{(a_1,r),(c_3,1-r),\infty_m\},
\{(b_2,r),(c_3,r),\infty_m\} \}.
\end{align*}
Take $ \mathcal{C}_j=\mathcal{C}_{j1}\cup\mathcal{C}_{j2}\cup\mathcal{C}_{j3}\cup\mathcal{C}_{j4}$. Then $\mathcal{C}_j$ is a 6-partial parallel class with respect to the group $\{j_1,j_2,j_3\}\times \{0,1\}$.
For $ 1\leq l<\beta $, let
$$\mathcal{D}_l =\mathcal{D}_{l+\beta-1}=\bigcup\limits_{\substack{0\leq j<n \\ \{a_1,b_2,c_3\}\in{{\cal H}_j}\cap{\cal P}_l}}\bigcup\limits_{r\in\{0,1\}}\{\{(a_1,r),(b_2,r),(c_3,r)\}\} $$
and
\begin{center}
$ 	\mathcal{D}_{2\beta-1} =\bigcup\limits_{ \{a_1,b_2,c_3\}\in{{\cal H}_j}\cap{\cal P}_\beta}\bigcup\limits_{r\in\{0,1\}}\{\{(a_1,r),(b_2,r),(c_3,r)\}\}  $.
\end{center}
Then $\mathcal{D}_1,\mathcal{D}_2,\ldots,\mathcal{D}_{2\beta-1}$ are 1-partial parallel classes with respect to the group $(\{\infty_1,\ldots,\infty_{m-1}\}\times \{0,1\})\cup \{\infty_m\}$. One can check that $(\bigcup_{j=0}^{n-1}{\cal C}_j)\cup(\bigcup_{l=1}^{2\beta-1}{\cal D}_l)$ forms a $(3,2)$-eframe of type $6^n (2m-1)^1$.

If $ n\equiv0\pmod3 $, for any $ i\in\{1,2,3\} $, construct a $ (3,2) $-frame of type $ 2^{n+1} $ on $\{0_i,1_i,\ldots,(n-1)_i,\infty_{m-1}\})\times \{0,1\} $ with groups $\{j_i\}\times \{0,1\}$ for $0\leq j<n$, and $\{\infty_{m-1}\}\times\{0,1\}$. Denote by $ \mathcal{G}'_{ji} $, $ \mathcal{G}''_{ji} $ and $ \mathcal{G}'_{\infty i} $, $ \mathcal{G}''_{\infty i} $ its 1-partial parallel classes with respect to the group $\{j_i\}\times \{0,1\}$  and the group $\{\infty_{m-1}\}\times\{0,1\}$, respectively. Write $ \mathcal{G}_{ji}=\mathcal{G}'_{ji}\cup \mathcal{G}''_{ji} $. Let
$$ \mathcal{C}'_{j1}=\mathcal{G}_{j1}\cup \mathcal{G}_{j2}\cup\mathcal{G}_{j3}. $$
Take
$ \mathcal{C}'_j=\mathcal{C}'_{j1}\cup\mathcal{C}_{j2}\cup\mathcal{C}_{j3}\cup\mathcal{C}_{j4} $  (note that $\beta=m-1$ now). Then $\mathcal{C}'_j$ is a 6-partial parallel class with respect to the group $\{j_1,j_2,j_3\}\times\{0,1\}$. For $ 1\leq l< \beta $, let
$$\mathcal{D}'_l =\mathcal{D}'_{l+\beta-1}=\bigcup\limits_{\substack{0\leq j<n \\ \{a_1,b_2,c_3\}\in{{\cal H}_j}\cap{\cal P}_l}}\bigcup\limits_{r\in\{0,1\}}\{\{(a_1,r),(b_2,r),(c_3,r)\}\}. $$
Let $\mathcal{D}'_{2m-1} =\mathcal{G}''_{\infty 1}\cup \mathcal{G}''_{\infty 2}\cup\mathcal{G}''_{\infty 3}$, $\mathcal{D}'_{2m-2} =\mathcal{G}'_{\infty 1}\cup \mathcal{G}'_{\infty 2}\cup\mathcal{G}'_{\infty 3}$ and
$$ \mathcal{D}'_{2m-3} =\bigcup\limits_{\substack{0\leq j<n \\ \{a_1,b_2,c_3\}\in{{\cal H}_j}\cap{\cal P}_\beta}}\bigcup\limits_{r\in\{0,1\}}\{\{(a_1,r),(b_2,r),(c_3,r)\}\}.$$
Then $\mathcal{D}'_1,\mathcal{D}'_2,\ldots,\mathcal{D}'_{2m-1}$ are 1-partial parallel classes with respect to the group $(\{\infty_1,\infty_2,\ldots,$ $\infty_{m-1}\}\times \{0,1\})\cup \{\infty_m\}$. One can check that $(\bigcup_{j=0}^{n-1}{\cal C}'_j)\cup(\bigcup_{l=1}^{2m-1}{\cal D}'_l)$ forms a $(3,2)$-eframe of type $6^n (2m-1)^1$.  \qed

\begin{Theorem}\label{IGDD4}
Let $ s\equiv 1,2\pmod{3} $, $s\geq 5$ and $s\notin \{7,11\}$. Let $0\leq m<s-1$ when $s\equiv 2\pmod{3}$ and $1\leq m\leq s-1$ when $s\equiv 1\pmod{3}$. Then there exist
\begin{itemize}
\item [$ (1) $] a $(3,2)$-IGDD of type $(2m+l;m,l)^1(9;3,3)^{t-1}(6;3,0)^{s-t}$ for any $ 0\leq l\leq m $ and $1\leq t\leq s$;
\item [$ (2) $] a $(3,2)$-IGDD of type $(2(2m-1)+l;2m-1,l)^1(18;6,6)^{t-1}(12;6,0)^{s-t}$ for any $ 0\leq l\leq 2m-1 $ and $1\leq t\leq s$.
\end{itemize}
\end{Theorem}

\proof Take $(3,2)$-eframes of types $3^{s-1} m^1$ (from Lemma \ref{eframe}) and $6^{s-1} (2m-1)^1$ (from Lemma \ref{eframe1}). Then apply Construction \ref{cons:IGDD} to complete the proof. \qed

\begin{Lemma}\label{eframe2}
\begin{itemize}
\item [$ (1) $] Let $ n\geq 4 $ and $ n\notin\{8,10,11,14,17\} $. Then there exist $(3,2)$-eframes of types $3^n 1^1$ and $3^n 2^1$.
\item [$ (2) $] Let $ n\geq 4 $ and $ n\notin\{6,8,10,11,14,17\} $. Then there exists a $(3,2)$-eframe of type $6^n 5^1$.

	\end{itemize}
\end{Lemma}

\proof By Lemma \ref{lem:eframe-direct}, there exists a $(3,2)$-eframe of type $h^n m^1$ for $(h,n,m)\in\{(3,5,1),(3,5,2),\\ (3,6,1),(3,6,2),(6,5,5)\}$. When $n\equiv0,1\pmod3 $, $ n\geq 4 $ and $n\notin\{6,10\}$, the required conclusion follows from Lemmas \ref{eframe} and \ref{eframe1}. Assume that $n\equiv 2\pmod3 $ and $n\geq 20$.

Let $m\in\{1,2\}$. Apply Lemma \ref{eframe} to form a $(3,2)$-eframe of type $3^{n-4} (12+m)^1$, and then fill in the group of size $ 12+m $ with a $(3,2)$-eframe of type $3^4 m^1$ (from Lemma \ref{eframe}) to obtain a $(3,2)$-eframe of types $3^n m^1$.

Apply Lemma \ref{eframe1} again to form a $(3,2)$-eframe of type $6^{n-4} 29^1$, and then fill in the group of size $ 29 $ with a $(3,2)$-eframe of type $6^4 5^1$ (from Lemma \ref{eframe1}) to obtain a $(3,2)$-eframe of type $6^n 5^1$.
\qed

\begin{Theorem}\label{IGDD2}
\begin{itemize}
\item [$ (1) $] Let $s\geq5$ and $ s\notin\{9,11,12,15,18\} $. Let $ m\in\{1,2\} $. Then there exists a $(3,2)$-IGDD of type $(2m+l;m,l)^1(9;3,3)^{t-1}(6;3,0)^{s-t}$ for any $ 0\leq l\leq m $ and $1\leq t\leq s$.
\item [$ (2) $] Let $s\geq5$ and $ s\notin\{7,9,11,12,15,18\} $. Then there exists a $(3,2)$-IGDD of type $(10+l;5,l)^1(18;6,6)^{t-1}(12;6,0)^{s-t}$ for any $ 0\leq l\leq 5 $ and $1\leq t\leq s$.

\end{itemize}
\end{Theorem}

\proof Take $(3,2)$-eframes of types $3^{s-1} 1^1$, $3^{s-1} 2^1$ (from Lemma \ref{eframe2}(1)), and $6^{s-1} 5^1$ (from Lemma \ref{eframe2}(2)). Then apply Construction \ref{cons:IGDD} to complete the proof. \qed

\subsection{$ \lambda=2 $}

Now, we are ready to construct some $\Diamond$-ITS$_2(2w+u+x,w,u) $ for $ x\in \{1,2,3,4\} $.

\begin{Lemma}\label{2w+u+x-small-lambda=2}
Let $(w,u,x)\in\{(6,3,1), (9,3,1), (9,6,1),(8,6,1),(11,9,1),(9,4,2),(10,9,2),(3,2$, $3),(6,4,3),(7,4,3),(8,4,3),(8,7,3), (9,8,3),(10,4,3),(10,7,3),(11,4,3),(11,10,3)\}$. Then there exists a $\Diamond$-ITS$_2(2w+u+x;w,u)$.
\end{Lemma}

\proof We give an explicit construction for each ITS in Appendix D.1. \qed

\begin{Lemma}\label{2w+u+x}
Let $u\equiv 0,1\pmod{3}$ and $u\geq 3$. Let $x$ be a nonnegative positive integer. If $w+x\geq 2u$, $w+x\equiv 0\pmod{u}$ and $x(w+u+x)\equiv0\pmod{3}$, then there exists a $\Diamond$-ITS$_2(2w+u+x;w,u)$.
\end{Lemma}

\proof Take a $ (3,2) $-GDD of type $ u^{\frac{w+x}{u}+1} w^1$ from Lemma \ref{GDD2}. Then place a TS$_2(u)$ (from Theorem \ref{ITS-one hole}) on each group of size $ u $ except one. \qed

\begin{Lemma}\label{w:small}
Let $ w\in\{26,32,35,44,53\} $ and $(x,u\pmod{3})\in\{(1,0),(3,1)\} $. If $w>u>1$, then there exists a $\Diamond$-ITS$_2(2w+u+x;w,u)$.
\end{Lemma}

\proof Apply Construction \ref{con:cons1} with $(a,b,c)=(0,0,x) $, $ h=3 $,
$$ m=\left\{
\begin{array}{lll}
5, & \hbox{if $w=26,32,44,53$,} \\
8, & \hbox{if $w=35$,} \\
\end{array}
\right.$$
$ s=\frac{w-m}{3}+1 $,
$$ t=\left\{
\begin{array}{lll}
\lfloor \frac{u}{3}\rfloor+1, & \hbox{if $3\lfloor \frac{u}{3}\rfloor\leq w-m$,} \\
s, & \hbox{otherwise,} \\
\end{array}
\right.$$
and $ l=u-3(t-1)$, where the needed $(3,2)$-IGDDs of type $(2m+l;m,l)^1(9;3,3)^{t-1}(6;3,0)^{s-t}$ are from Theorem \ref{IGDD4}(1), the needed $\Diamond$-ITS$_2(2m+l+x;m,l)$s can be found in the following table
\begin{center}
\begin{tabular}{|lllc|lllc|}\hline
$x$ & $m$ & $l$ & Source & $x$ & $m$ & $l$ & Source \\\hline
$1$ & $5$ & $0$ & Theorem \ref{ITS-one hole} & $3$ & $5$ & $1$ & Theorem \ref{ITS-one hole} \\
$1$ & $5$ & $3$ & Lemma \ref{2w+u+x} & $3$ & $5$ & $4$ & Lemma \ref{2w+u+x}  \\\hline
$1$ & $8$ & $0$ & Theorem \ref{ITS-one hole} & $3$ & $8$ & $1$ & Theorem \ref{ITS-one hole} \\
$1$ & $8$ & $3$ & Lemma \ref{2w+u+x} & $3$ & $8$ & $4$ & Lemma \ref{2w+u+x-small-lambda=2} \\
$1$ & $8$ & $6$ & Lemma \ref{2w+u+x-small-lambda=2} & $3$ & $8$ & $7$ & Lemma \ref{2w+u+x-small-lambda=2} \\
\hline
\end{tabular} ,
\end{center}
the needed $\Diamond$-ITS$_2(6+x;3,x)$s come from Theorem \ref{ITS-one hole} if $ x=1$ and from Lemma \ref{w=u small} if $x=3$, and the needed $\Diamond$-ITS$_2(9+x;3,3,x;0,0)$s are from Lemma \ref{w=u small} if $x=1$ and from Lemma \ref{lattice-small} if $x=3$. \qed

\begin{Lemma}\label{2w+u+1,2}
Let $w\geq u>1$ and $u(w-1) \equiv 0 \pmod{3}$. Then there exists a $\Diamond$-ITS$_2(2w+u+1;w,u)$.
\end{Lemma}

\proof The case of $w=u$ follows from Lemma \ref{w=u small}. Assume that $w> u>1$.

For $w \equiv 1 \pmod{3}$, apply Construction \ref{con:cons1} with $(m,l,h)=(1,1,1)$, $(s,t)=(w,u)$ and $(a,b,c)=(0,0,1)$, where the needed $(3,2)$-IGDD of type $(3;1,1)^{u}(2;1,0)^{w-u}$ is from Theorem \ref{IGDD}, and the needed $\Diamond$-ITS$_2(4;1,1)$, $\Diamond$-ITS$_2(3;1,1;0)$ and $\Diamond$-ITS$_2(4;1,1,1;0,0)$ are from Theorem \ref{ITS-one hole}.

For $w \equiv u\equiv 0 \pmod{3}$, if $(w,u)\in\{(6,3), (9,3), (9,6) \}$, apply Lemma \ref{2w+u+x-small-lambda=2}. If $w\geq 12$, apply Construction \ref{con:cons1} with $(m,l,h)=(3,3,3)$, $(s,t)=(w/3,u/3)$ and $(a,b,c)=(0,0,1)$, where the needed $(3,2)$-IGDD of type $(9;3,3)^{t}(6;3,0)^{s-t}$ is from Theorem \ref{IGDD}, the needed $\Diamond$-ITS$_2(10;3,3)$ and $\Diamond$-ITS$_2(10;3,3,1;0,0)$ are from Lemma \ref{w=u small}, and the needed $\Diamond$-ITS$_2(7;3,1;0)$ is from Theorem \ref{ITS-one hole}.

For $w \equiv 2 \pmod{3}$ and $u \equiv 0 \pmod{3}$, if $(w,u)\in\{(8,6),(11,9)\} $, apply Lemma \ref{2w+u+x-small-lambda=2}.
If $(w,u)\in\{ (5,3),(8,3),(11,3),(11,6)\} $, apply Lemma \ref{2w+u+x}. If  $ w\geq14 $ and $ w\notin\{26,32,35,44,53\} $, apply Construction \ref{con:cons1} with $(m,l,h)=(2,2,3)$, $(s,t)=((w+1)/3,u/3)$ and $(a,b,c)=(0,1,2)$, where the needed $(3,2)$-IGDD of type $(6;2,2)^1(9;3,3)^{t-1}(6;3,0)^{s-t}$ is from Theorem \ref{IGDD2}(1), the needed $\Diamond$-ITS$_2(8;2,3)$ is from Lemma \ref{2w+u,2}, and the needed $\Diamond$-ITS$_2(11;3,4,2;0,1)$ is from Remark \ref{rek:lattice} and Lemma \ref{lattice-small}. If $ w\in\{26,32,35,44,53\} $, apply Lemma \ref{w:small}.  \qed

\begin{Lemma}\label{2w+u+2,2}
Let $w\geq u>1$ and $(u+1)(w+1) \equiv 2 \pmod{3}$. Then there exists a $\Diamond$-ITS$_2(2w+u+2;w,u)$.
\end{Lemma}

\proof The case of $w=u$ follows from Lemma \ref{w=u small}. Assume that $w> u>1$.

For $(w,u)\in\{(4,3),(6,4),(10,3),(10,6)\} $, apply Lemma \ref{2w+u+x}. For $(w,u)\in\{(7,3),(9,7)\}$, apply Theorem \ref{lambda=1}. For $(w,u)=(7,6) $, take a $(3,2)$-IGDD of type $(7;2,2)^3$ from Lemma \ref{IGDD small}, adjoin one infinite point, and fill in three groups using a $\Diamond$-ITS$_2(8;3,2)$ (from Lemma \ref{2w+u,2}). For $(w,u)\in\{(9,4),(10,9)\}$, apply Lemma \ref{2w+u+x-small-lambda=2}.

For $w\geq 12$, apply Construction \ref{con:cons1} with $s=\lfloor w/3\rfloor$, $t=\lceil u/3\rceil$, $m=h=3$ and
$$ (l,a,b,c)=\left\{
\begin{array}{ll}
(1,0,0,2), & \hbox{$u \equiv 1\pmod{3}$ and $w \equiv 0\pmod{3}$,} \\
(3,1,0,4), & \hbox{$u \equiv 0\pmod{3}$ and $w \equiv 1\pmod{3}$,} \\
\end{array}
\right.$$
where the needed $(3,2)$-IGDD of type $(6+l;3,l)^1(9;3,3)^{t-1}(6;3,0)^{s-t}$ is from Theorem \ref{IGDD}, the needed $\Diamond$-ITS$_2(9;3,1)$ and $\Diamond$-ITS$_2(13;4,3)$ are from Theorem \ref{ITS-one hole} and Lemma \ref{2w+u+x}, respectively, the needed $\Diamond$-ITS$_2(8;3,2;0)$ and $\Diamond$-ITS$_2(10;4,4;1)$ are from Lemma \ref{2w+u,2}, and the needed $\Diamond$-ITS$_2(11;3,3,2;0,0)$ and $\Diamond$-ITS$_2(13;4,3,4;1,0)$ are from Lemma \ref{lattice-small}. \qed

\begin{Lemma}\label{2w+u+3,2}
Let $w\geq u>1$ and $w(u-1)\equiv 0 \pmod{3}$. Then there exists a $\Diamond$-ITS$_2(2w+u+3;w,u)$.
\end{Lemma}

\proof The case of $w=u$ follows from Lemma \ref{w=u small}. Assume that $w> u>1$.

For $(w,u)\in\{(3,2),(6,4),(7,4),(8,4),(8,7)$, $(9,8),(10,4),(10,7),(11,4),(11,10)\}$, apply Lemma \ref{2w+u+x-small-lambda=2}.
For $(w,u)\in\{(5,4),(6,3),(9,3),(9,4),(9$, $6),(11,7)\}$, apply Lemma \ref{2w+u+x}. For $(w,u)\in\{(6,5)$, $(9,5),(9,7)\}$, take a $(3,2)$-IGDD of type $(n;m,l)^3$ with $(n,m,l)\in \{(6,2,1),(8,3,1),(9,3$, $2)\}$ (from Lemma \ref{IGDD small}), adjoin $ c\in\{2,2,1\} $ infinite points, and fill in three groups using a $\Diamond$-ITS$_2(n+c;m,l+c)$, which exists by Lemmas \ref{2w+u,2} and \ref{2w+u+1,2}. For $(w,u)\in\{(6,2),(9,2)\}$, start with a $\Diamond$-ITS$_2(2w+5;w,5)$ (from Lemma \ref{2w+u,2}), and then place a $\Diamond$-ITS$_2(5;2,0)$ (from Theorem \ref{ITS-one hole}) on the hole of size 5.

For $w \equiv 0 \pmod{3}$ and $w\geq 12$, apply Construction \ref{con:cons1} with $s=w/3$, $t=\lceil u/3\rceil$, $m=h=3$, $ l=u+3-3t $ and $ (a,b,c)=(0,0,3) $, where the needed $(3,2)$-IGDD of type $(6+l;3,l)^1(9;3,3)^{t-1}(6;3,0)^{s-t}$ is from Theorem \ref{IGDD}, the needed $\Diamond$-ITS$_2(9+l;3,l)$ is from Theorem \ref{ITS-one hole} or has been constructed in the above two paragraphs, the needed $\Diamond$-ITS$_2(9;3,3;0)$ is from Lemma \ref{2w+u,2}, and the needed $\Diamond$-ITS$_2(12;3,3,3;0,0)$ is from Lemma \ref{lattice-small}.

For $w \equiv 1,2 \pmod{3}$, $w\geq 13$ and $u \equiv 1 \pmod{3}$, if $\lceil w/3\rceil \notin\{9,11,12,15,18\}$, then apply Construction \ref{con:cons1} with $(s,t)=(\lceil w/3\rceil,(u+2)/3)$, $(a,b,c)=(0,0,3)$, $(l,h)=(1,3)$, and $m\in\{1,2\}$ such that $m\equiv w\pmod{3}$, where the needed $(3,2)$-IGDD of type $(2m+1;m,1)^1(9;3,3)^{t-1}(6;3,0)^{s-t}$ is from Theorem \ref{IGDD2}(1), and the needed $\Diamond$-ITS$_2(6;1,1)$ and $\Diamond$-ITS$_2(8;$ $2,1)$ are from Theorem \ref{ITS-one hole}.

For $ w\in\{26,32,35,44,53\} $ and $u \equiv 1 \pmod{3}$, apply Lemma \ref{w:small}. For $ w\in\{25,31,34,43,52\} $ and $u \equiv 1 \pmod{3}$, apply Construction \ref{con:cons1} with $ h=3 $, $(a,b,c)=(0,0,3)$,
$$ m=\left\{
\begin{array}{lll}
4, & \hbox{if $w=25,31,43,52$,} \\
7, & \hbox{if $w=34$,} \\
\end{array}
\right.$$
$ s=\frac{w-m}{3}+1 $,
$$ t=\left\{
\begin{array}{lll}
\frac{u+2}{3}, & \hbox{if $u\leq w-m+1$,} \\
s, & \hbox{if $w-m+1< u\leq w-3$,} \\
\end{array}
\right.$$
and $ l=u-3(t-1)$, where the needed $(3,2)$-IGDD of type $(2m+l;m,l)^1(9;3,3)^{t-1}(6;3,0)^{s-t}$ is from Theorem \ref{IGDD4}(1), the needed $\Diamond$-ITS$_2(12;4,1)$ and $\Diamond$-ITS$_2(18;7,1)$ are from Theorem \ref{ITS-one hole}, and the needed $\Diamond$-ITS$_2(21;7,4)$ is from Lemma \ref{2w+u+x-small-lambda=2}. \qed

\begin{Lemma}\label{2w+u+4,2}
Let $w> u>1$, $w \equiv 2 \pmod{3}$ and $u \equiv 0 \pmod{3}$. Then there exists a $\Diamond$-ITS$_2(2w+u+4;w,u)$.
\end{Lemma}

\proof By the analysis in the last paragraph of Section 2.1, it suffices to examine the existence of a $\Diamond$-ITS$_2(2w+u+4;w,u)$ for $wu\equiv0\pmod2$ and $w> u>1$.

For $ (w,u)=\{(8,3),(8,6),(20,3),(20,6),(20,12)\} $, apply Lemma \ref{2w+u+x}. For $(w,u)\in\{(11,6)$, $(17,6),(20,9),(20,15)\}$, take a $(3,2)$-IGDD of type $(n;m,l)^3$ with $(n,m,l)\in \{(9,3,1),(13,5,1)$, $(15,6,1),(18,6,4)\}$ (from Lemma \ref{IGDD small}), adjoin $ c\in\{5,5,8,5\} $ infinite points, place a $\Diamond$-ITS$_2(n+c;m+2,l+c-2)$ (from Lemma \ref{2w+u,2}) on one group, and fill in the other two groups using a $\Diamond$-ITS$_2(n+c;m+2,l+c-2,c;2,c-2)$ (from Lemma \ref{lattice-small}).

For $ w\in\{14,23,41,50,68,86,104\} $ and $u \equiv 0 \pmod{3}$, or $ (w,u)=(17,12)$, set $ m=\frac{w-2}{3} $, $ l=\frac{u}{3} $ and $ n=2m+l+2 $. Start from a $(3,2)$-IGDD of type $(n;m,l)^3$ (from Lemma \ref{IGDD small}), adjoin $ 2 $ infinite points, and fill in three groups using a $\Diamond$-ITS$_2(n+2;m+2,l)$ (from Lemma \ref{2w+u,2}).

For $(w,u)=(20,18)$, apply Construction \ref{con:cons1} with $s=t=7$, $(h,m,l)=(3,2,0)$ and $(a,b,c)=(0,0,4)$, where the needed $(3,2)$-IGDD of type $(4;2,0)(9;3,3)^6$ is from Theorem \ref{IGDD2}(1), the needed $\Diamond$-ITS$_2(8;2,0)$ is from Theorem \ref{ITS-one hole}, and the needed  $\Diamond$-ITS$_2(13;3,3,4;0,0)$ is from Lemma \ref{lattice-small}.

For  $w\geq 26$ and $ w\notin\{32,38,41,50,53,62,65,68,71,86,89,104,107\} $, apply Construction \ref{con:cons1} with $s=\lceil w/6\rceil$, $t=\lceil (u+3)/6\rceil$, $h=6$ and
$$ (m,l,a,b,c)=\left\{
\begin{array}{llll}
(2,0,0,0,4), & \hbox{$u \equiv 0\pmod{6}$ and $w \equiv 2\pmod{6}$,} \\
(2,1,0,2,6), & \hbox{$u \equiv 3\pmod{6}$ and $w \equiv 2\pmod{6}$,} \\
(5,0,0,0,4), & \hbox{$u \equiv 0\pmod{6}$ and $w \equiv 5\pmod{6}$,}
\end{array}
\right.$$
where the needed $(3,2)$-IGDD of type $(2m+l;m,l)^1(18;6,6)^{t-1}(12;6,0)^{s-t}$ is from Theorem \ref{IGDD2}(2) if $ m=5 $ and from Theorem \ref{IGDD3} if $ m=2 $, the needed $\Diamond$-ITS$_2(8;2,0)$ and $\Diamond$-ITS$_2(14;5,0)$ are from Theorem \ref{ITS-one hole}, the needed $\Diamond$-ITS$_2(11;3,2)$ is from Lemma \ref{2w+u+3,2}, the needed $\Diamond$-ITS$_2(16;6,4)$ and $\Diamond$-ITS$_2(18;6,6)$ are from Lemma \ref{2w+u,2}, and the needed $\Diamond$-ITS$_2(22;6,6,4;0,0)$ and $\Diamond$-ITS$_2(24;6$, $8,6;0,2)$ are from Lemma \ref{lattice-small} and Remark \ref{rek:lattice}.

For $ w=32$, $ u \equiv 0 \pmod{3} $ and $ u\leq 30 $, let $ l=\frac{u}{3} $, take a $(3,2)$-IGDD of type $(22+l;10,l)^3$ from Lemma \ref{IGDD small}, adjoin two infinite point $ \{\infty_1,\infty_2\} $, and place on each group together with $ \{\infty_1,\infty_2\} $ the triples of a $\Diamond$-ITS$_2(24+l;12,l)$ (from Lemma \ref{2w+u,2}).

For $ w\in\{38,53,62,65,71,89,107\} $, $u \equiv 0 \pmod{3}$ and $w>u>1$, apply Construction \ref{con:cons1} with $ h=6 $,  $(a,b,c)=(0,0,4)$,
$$ m=\left\{
\begin{array}{lll}
8, & \hbox{if $w=38,62$,} \\
11, & \hbox{if $w=53,65,89,107$,} \\
17, & \hbox{if $w=71$,} \\
\end{array}
\right.$$
$ s=\frac{w-m}{6}+1 $,
$$ t=\left\{
\begin{array}{lll}
 \lfloor u/6\rfloor+1, & \hbox{if $6\lfloor u/6\rfloor\leq w-m$}, \\
s, & \hbox{otherwise,} \\
\end{array}
\right.$$
and $ l=u-6(t-1)$, where the needed $(3,2)$-IGDDs of type $(2m+l;m,l)^1(18;6,6)^{t-1}(12;6,0)^{s-t}$ are from Theorem \ref{IGDD-add} if $m=8$ and from Theorem \ref{IGDD4}(2) if $m\in\{11,17\}$, the needed $\Diamond$-ITS$_2(2m+l+4;m,l)$s
for $m\in\{8,11,17\}$, $l\equiv 0\pmod{3}$ and $0\leq l\leq m-2$
are from Theorem \ref{ITS-one hole} if $l=0$, and otherwise from the the former part of the current proof, the needed $\Diamond$-ITS$_2(16;6,4;0)$ is from Lemma \ref{2w+u,2}, and the needed $\Diamond$-ITS$_2(22;6,6,4;0,0)$ is from Lemma \ref{lattice-small}. \qed

\subsection{ $\lambda=3$}

The following lemma can be obtained straightforwardly from Lemma \ref{main result1 }. We write it here because we will use it in Lemma \ref{th:hole1}.

\begin{Lemma}\label{2w+u+2,3}
Let $u,w$ and $v$ be odd positive integers and $3\leq u\leq w$. Then there exists a $\Diamond$-ITS$_3(2w+u+2;w,u)$.
\end{Lemma}

\subsection{ $ \lambda=6 $}



\begin{Lemma}\label{2w+u+x-small-lambda=6}
Let $(w,u)\in\{(2,2),(3,2),(5,2),(5,4),(6,2),(6,4),(6,5),(8,2),(8,4),(8,5),(8,7)\}$. Then there exists a $\Diamond$-ITS$_6(2w+u+1;w,u)$.
\end{Lemma}
\proof For each case except for $ (w,u)=(8,2)$, we give an explicit construction in Appendix D.2. For $ (w,u)=(8,2)$, take a $\Diamond$-ITS$_6(19;8,7;2)$ from Lemma \ref{2w+u,2} by repeating blocks, and fill in the hole of size $7$ with a $\Diamond$-ITS$_6(7;2,2)$.
 \qed

\begin{Lemma}\label{w:small:lambda=6}
	Let $ w\in\{26,32,35,44,53\} $ and $u\equiv1,2\pmod{3}$. If $w>u>1$, then there exists a $\Diamond$-ITS$_6(2w+u+1;w,u)$.
\end{Lemma}

\proof Apply Construction \ref{con:cons1} with $(a,b,c)=(0,0,1) $, $ h=3 $,
$$ m=\left\{
\begin{array}{lll}
5, & \hbox{if $w=26,32,44,53$,} \\
8, & \hbox{if $w=35$,} \\
\end{array}
\right.$$
$ s=\frac{w-m}{3}+1 $,
$$ t=\left\{
\begin{array}{lll}
\lfloor \frac{u}{3}\rfloor+1, & \hbox{if $3\lfloor \frac{u}{3}\rfloor\leq w-m$,} \\
s, & \hbox{otherwise,} \\
\end{array}
\right.$$
and $ l=u-3(t-1)$, where the needed $(3,6)$-IGDDs of type $(2m+l;m,l)^1(9;3,3)^{t-1}(6;3,0)^{s-t}$ are from Theorem \ref{IGDD4}(1) (by repeating blocks of a $(3,2)$-IGDD with the same type), the needed $\Diamond$-ITS$_6(2m+l+1;m,l)$s can be found in the following table
\begin{center}
\begin{tabular}{|cccc|}\hline
$x$ & $m$ & $l\in$ & Source  \\\hline
$1$ & $5$ & $\{1\}$ & Theorem \ref{ITS-one hole}  \\
$1$ & $5$ & $\{2,4\}$ & Lemma \ref{2w+u+x-small-lambda=6} \\
$1$ & $8$ & $\{1\}$ & Theorem \ref{ITS-one hole}  \\
$1$ & $8$ & $\{2,4,5,7\}$ & Lemma \ref{2w+u+x-small-lambda=6}  \\
		\hline
	\end{tabular} ,
\end{center}
the needed $\Diamond$-ITS$_6(7;3,1)$ comes from Theorem \ref{ITS-one hole}, and the needed $\Diamond$-ITS$_6(10;3,3,1;0,0)$ is from Lemma \ref{w=u small}. \qed

\begin{Lemma}\label{2w+u+1,6}
There exists a $\Diamond$-ITS$_6(2w+u+1;w,u)$ for any $ w\geq u>1 $.
\end{Lemma}

\proof The case of $w=u$ follows from Lemma \ref{w=u small}. Assume that $w> u>1$. By the analysis in the last paragraph of Section 2.1, it suffices to examine the existence of a $\Diamond$-ITS$_6(2w+u+1;w,u)$ for $u(w-1) \equiv 1,2 \pmod{3}$.

For $w \equiv 0 \pmod{3}$, $w\geq 12$ and $ u\equiv1,2\pmod3 $, apply Construction \ref{con:cons1} with $(s,t)=(w/3,\lceil u/3\rceil)$, $(m,h)=(3,3)$, $ (a,b,c)=(0,0,1) $ and $ l\in\{1,2\} $ such that $ l\equiv u\pmod3 $, where the needed $(3,6)$-IGDD of type $(6+l;3,l)^1(9;3,3)^{t-1}(6;3,0)^{s-t}$ is from Theorem \ref{IGDD}, the needed  $\Diamond$-ITS$_6(7;3,1)$ and $\Diamond$-ITS$_6(8;3,1)$ are  from Theorem \ref{ITS-one hole}, the needed $\Diamond$-ITS$_6(9;3,2)$ is from Lemma \ref{2w+u+x-small-lambda=6}, and the needed $\Diamond$-ITS$_6(10;3,3,1;0,0)$ is just a $\Diamond$-ITS$_6(10;3,3)$  which is from Lemma \ref{w=u small}.

For $ w\in\{26,32,35,44,53\} $ and $ u\equiv1,2\pmod3 $, apply Lemma \ref{w:small:lambda=6}. For $w \equiv 2 \pmod{3}$, $ w\geq14 $, $ w\not\in\{26,32,35,44,53\} $ and $ u\equiv1,2\pmod3 $, apply Construction \ref{con:cons1} with $(s,t)=((w+1)/3,\lceil u/3\rceil)$, $ (m,h)=(2,3)$, $ (a,b,c)=(0,0,1) $ and $ l\in\{1,2\} $ such that $ l\equiv u\pmod3 $, where the needed $(3,6)$-IGDD of type $(4+l;2,l)^1(9;3,3)^{t-1}(6;3,0)^{s-t}$ comes from Theorem \ref{IGDD2}(1), the needed  $\Diamond$-ITS$_6(6;2,1)$ is from Theorem \ref{ITS-one hole}, and the needed  $\Diamond$-ITS$_6(7;2,2)$ is from Lemma \ref{2w+u+x-small-lambda=6}.

For $(w,u)\in\{(2,2),(3,2),(5,2),(5,4),(6,2),(6,4), (6,5), (8,2),(8,4),(8,5),(8,7)\}$, apply Lemma \ref{2w+u+x-small-lambda=6}. For $ w\in\{9,11\} $, apply Construction \ref{con:cons1} with $s=(w+1)/2$, $t=\lfloor u/2\rfloor+1$, $ m=1 $, $h=2$, $ l=u-2(t-1)$ and $ (a,b,c)=(0,0,1)$, where the needed $(3,6)$-IGDD of type $(2+l;1,l)^1(6;2,2)^{t-1}(4;2,0)^{s-t}$ is from Theorem \ref{special,IGDD}, the needed  $\Diamond$-ITS$_6(5;2,1)$, $\Diamond$-ITS$_6(3;1,0)$ and $\Diamond$-ITS$_6(4;1,1)$ are  from Theorem \ref{ITS-one hole}, and the needed $\Diamond$-ITS$_6(7;2,2,1;0,0)$ is just a $\Diamond$-ITS$_6(7;2,2)$ which is from Lemma \ref{2w+u+x-small-lambda=6}. \qed

\section{The cases of large $ v $}

\begin{Lemma}\label{hole}
If a $\Diamond$-ITS$_\lambda(v;w,b)$ and a $\Diamond$-ITS$_\lambda(b;u,0)$ exist, then a $\Diamond$-ITS$_\lambda(v;w,u)$ exists.
\end{Lemma}

\proof Place a $\Diamond$-ITS$_\lambda(b;u,0)$ on the hole of size $b$ of a $\Diamond$-ITS$_\lambda(v;w,b)$. \qed

\begin{Lemma} \label{th:hole1}
Let $\lambda\in\{2,3,6\}$. Let $w\geq u>1$ and $v\geq 2\max\{w,2u+1\}+w$. If $\lambda({v\choose 2}-{u\choose 2}-{w\choose 2}) \equiv 0\pmod3$ and $\lambda (v-w) \equiv \lambda (v-u) \equiv \lambda (v-1) \equiv 0 \pmod 2$, then there exists a $\Diamond$-ITS$_\lambda(v;w,u)$.
\end{Lemma}

\proof Since $v\geq 2\max\{w,2u+1\}+w$, we shall pick up proper $a\geq \max\{w-1,2u+1\}$ and $x\in \{0,1,2,3,4\}$ such that (1) $v-w=2a+x$ and $\lambda({a\choose 2}-{u\choose 2})\equiv0\pmod{3}$; (2) there exist a $\Diamond$-ITS$_\lambda(2a+w+x;a,w)$ and a $\Diamond$-ITS$_\lambda(a;u,0)$. Then one can use Lemma \ref{hole} to obtain a $\Diamond$-ITS$_\lambda(v;w,u)$. Note that by Theorem \ref{ITS-one hole}, the condition $\lambda({a\choose 2}-{u\choose 2})\equiv0\pmod{3}$ is necessary for the existence of a $\Diamond$-ITS$_\lambda(a;u,0)$.

When $\lambda=2$, $\lambda({v\choose 2}-{u\choose 2}-{w\choose 2}) \equiv 0\pmod3$ yields $(v-w)(v+w-1)\equiv u(u-1)\pmod{3}$, which implies the following congruence conditions
\begin{center}
\begin{tabular}{|lll|}\hline
$v-w\equiv 0\pmod{3}$ & $u\equiv 0,1\pmod{3}$ & \\\hline
$v-w\equiv 1\pmod{3}$ & $u\equiv 0,1\pmod{3}$ & $w\equiv 0\pmod{3}$ \\
& $u\equiv 2\pmod{3}$ & $w\equiv 1\pmod{3}$ \\\hline
$v-w\equiv 2\pmod{3}$ & $u\equiv 0,1\pmod{3}$ & $w\equiv 1\pmod{3}$ \\
& $u\equiv 2\pmod{3}$ & $w\equiv 0\pmod{3}$ \\
\hline
\end{tabular} .
\end{center}
When $ u\equiv0,1\pmod 3 $, the choices for $ (x,a \pmod 3) $ are (note that $a=(v-w-x)/2$)
\begin{center}
\begin{tabular}{|c|l|l|l|l|l|l|}\hline
$v-w \pmod{6}$ & 0 & 1 & 2 & 3 & 4 & 5 \\
$x$ & 0 & 1 & 0 & 1 & 2 & 3 \\
$a \pmod{3}$ & 0 & 0 & 1 & 1 & 1 & 1 \\
\hline
\end{tabular} .
\end{center}
When $ u\equiv2\pmod 3 $, the choices for $ (x,a \pmod 3) $ are
\begin{center}
\begin{tabular}{|c|l|l|l|l|l|l|}\hline
$v-w \pmod{6}$ & 1 & 2 & 4 & 5 \\
$x$ & 3 & 4 & 0 & 1 \\
$a \pmod{3}$ & 2 & 2 & 2 & 2 \\
\hline
\end{tabular} .
\end{center}
Since $v\geq 2\max\{w,2u+1\}+w$ and $a=(v-w-x)/2$,  we have $ a\geq  \max\{w,2u+1\}-\frac{x}{2}$. Clearly $ a\geq  \max\{w,2u+1\}$ when $ x\in\{0,1\} $, and $ a\geq  \max\{w,2u+1\}-1$ when $ x\in\{2,3\} $. If $x=2$, then $ u\equiv0,1\pmod3 $, $v-w \equiv4\pmod{6}$ and $a \equiv1\pmod{3}  $. In this case, $ w\equiv0\pmod3 $, so $ \max\{w,2u+1\}-1\equiv 0,2 \pmod3$, which implies $ a\geq  \max\{w,2u+1\}$. Similar arguments show that $ a\geq  \max\{w,2u+1\}$ for the two cases of $x=3$. If $x=4$, then $ u\equiv2\pmod3 $, $v-w \equiv2\pmod{6}$ and $a \equiv2\pmod{3} $. In this case, $ w\equiv0\pmod3 $ and $ a\geq  \max\{w,2u+1\}-2$. If $ 2u+1\geq w$, then  $ \max\{w,2u+1\}-2\equiv 0\pmod 3 $ yields $ a\geq  \max\{w,2u+1\}$. If $ 2u+1< w$, then $ \max\{w,2u+1\}-2\equiv 1\pmod 3 $ yields $ a\geq  \max\{w-1,2u+1\}$. The needed $\Diamond$-ITS$_2(2a+w+x;a,w)$s are from Lemmas \ref{2w+u,2}, \ref{2w+u+1,2}, \ref{2w+u+2,2}, \ref{2w+u+3,2} and \ref{2w+u+4,2} (note that when $ a=w-1 $, the needed $\Diamond$-ITS$_2(3w+2;w-1,w)$ is from Lemma \ref{2w+u+3,2}). The needed $\Diamond$-ITS$_2(a;u,0)$s are from Theorem \ref{ITS-one hole}.

When $ \lambda=3 $, $\lambda (v-w) \equiv \lambda (v-u) \equiv \lambda (v-1) \equiv 0 \pmod 2$ yields $vuw\equiv 1\pmod{2}$. The choices for $ (x,a \pmod 6) $ are
\begin{center}
		\begin{tabular}{|c|l|l|l|l|l|l|}\hline
			$v-w \pmod{12}$ & 0 & 2 & 4 & 6 & 8 & 10 \\
			$x$ & 2 & 0 & 2 & 0 & 2 & 0 \\
			$a \pmod{6}$ & 5 & 1 & 1 & 3 & 3 & 5 \\
			\hline
		\end{tabular} .
\end{center}
Since $v\geq 2\max\{w,2u+1\}+w$, $a=(v-w-x)/2$ and $x\in\{0,2\}$, we have $ a\geq  \max\{w,2u+1\}-1$. Note that $ a\equiv w \equiv 2u+1\equiv1\pmod2$, so $ a\geq  \max\{w,2u+1\}$. The needed $\Diamond$-ITS$_3(2a+w+x;a,w)$s are from Lemmas \ref{2w+u,3} and \ref{2w+u+2,3}. The needed $\Diamond$-ITS$_3(a;u,0)$s are from Theorem \ref{ITS-one hole}.

When $ \lambda=6 $, the choices for $ (x,a \pmod 3) $ are
\begin{center}
\begin{tabular}{|c|l|l|l|l|l|l|}\hline
$v-w \pmod{6}$ & 0 & 1 & 2 & 3 & 4 & 5 \\
$x$ & 0 & 1 & 0 & 1 & 0 & 1 \\
$a \pmod{3}$ & 0 & 0 & 1 & 1 & 2 & 2 \\
\hline
\end{tabular} .
\end{center}
Since $x\in\{0,1\}$,  we have $ a\geq  \max\{w,2u+1\}$. The needed $\Diamond$-ITS$_6(2a+w+x;a,w)$s are from Lemmas \ref{2w+u,6} and \ref{2w+u+1,6}. The needed $\Diamond$-ITS$_6(a;u,0)$s are from Theorem \ref{ITS-one hole}. \qed

\begin{Lemma} \label{th:hole2}
Let $\lambda\in\{2,3,6\}$. Let $w\geq 2u+1>3$ and $v\geq 2w+2u+1$. If $\lambda({v\choose 2}-{u\choose 2}-{w\choose 2}) \equiv 0\pmod3$ and $\lambda (v-w) \equiv \lambda (v-u) \equiv \lambda (v-1) \equiv 0 \pmod 2$, then there exists a $\Diamond$-ITS$_\lambda(v;w,u)$.
\end{Lemma}

\proof When $ v\geq 3w $, the conclusion follows from Lemma \ref{th:hole1}. Assume that $2w+2u+1\leq v\leq 3w$. Since $v\geq 2w+2u+1$, we shall pick up proper $b\geq 2u+1$ and $x\in \{0,1,2\}$ such that (1) $v-2w=b+x$ and $\lambda{b\choose 2}\equiv\lambda{u\choose 2}\pmod 3$; (2) there exist a $\Diamond$-ITS$_\lambda(2w+b+x;w,b)$ and a $\Diamond$-ITS$_\lambda(b;u,0)$. Then one can use Lemma \ref{hole} to obtain a $\Diamond$-ITS$_\lambda(v;w,u)$. Note that $b=v-2w-x$ and $2w+2u+1\leq v\leq 3w$, so $2u+1-x\leq b\leq w-x$.

When $\lambda=2$, $\lambda({v\choose 2}-{u\choose 2}-{w\choose 2}) \equiv 0\pmod3$ yields $(v-w)(v-2w-1)\equiv u(u-1)\pmod{3}$, which implies the following congruence conditions
\begin{center}
\begin{tabular}{|lll|}\hline
$v-2w\equiv 0\pmod{3}$ & $u\equiv 0,1\pmod{3}$ & $w\equiv 0\pmod{3}$ \\
& $u\equiv 2\pmod{3}$ & $w\equiv 1\pmod{3}$  \\\hline
$v-2w\equiv 1\pmod{3}$ & $u\equiv 0,1\pmod{3}$ & \\\hline
$v-2w\equiv 2\pmod{3}$ & $u\equiv 0,1\pmod{3}$ & $w\equiv 1\pmod{3}$ \\
& $u\equiv 2\pmod{3}$ & $w\equiv 0\pmod{3}$ \\
\hline
\end{tabular} .
\end{center}
When $ u\equiv 0,1\pmod 3 $, the choices for $ (x,b \pmod 3) $ are (note that $b=v-2w-x$)
\begin{center}
\begin{tabular}{|c|l|l|l|l|l|l|}\hline
$v-2w \pmod{3}$ & 0 & 1 & 2 \\
$x$ & 0 & 0& 1 \\
$b \pmod{3}$ & 0 & 1 & 1 \\
\hline
\end{tabular} .
\end{center}
When $ u\equiv 2\pmod 3 $, the choices for $ (x,b \pmod 3) $ are
\begin{center}
\begin{tabular}{|c|l|l|l|l|l|l|}\hline
$v-2w \pmod{3}$ & 0 & 2 \\
$x$ & 1 & 0 \\
$b \pmod{3}$ & 2 & 2 \\
\hline
\end{tabular} .
\end{center}
Clearly $2u+1\leq b\leq w$ if $ x=0 $. If $x=1$, when $ v-2w \equiv2\pmod{3} $ and $ u\equiv 0,1\pmod 3 $, we have $ w\equiv 1\pmod 3 $ and $b \equiv1\pmod{3}$. Then $b\geq 2u+1-x=2u$ implies that $ b\geq 2u+1 $. When $ v-2w \equiv0\pmod{3} $ and $ u\equiv 2\pmod 3 $, we have $ w\equiv 1\pmod 3 $ and $b \equiv2\pmod{3}$. Then $b\geq 2u+1-x=2u$ implies that $ b\geq 2u+1 $. The needed $\Diamond$-ITS$_2(2w+b+x;w,b)$s are from Lemmas \ref{2w+u,2} and \ref{2w+u+1,2}. The needed $\Diamond$-ITS$_2(b;u,0)$s are from Theorem \ref{ITS-one hole}.

When $ \lambda=3 $, $\lambda (v-w) \equiv \lambda (v-u) \equiv \lambda (v-1) \equiv 0 \pmod 2$ yields $vuw\equiv 1\pmod{2}$. The choices for $ (x,b \pmod 3) $ are
\begin{center}
	\begin{tabular}{|c|l|l|l|l|l|l|}\hline
			$v-2w \pmod{6}$ & 1 & 3& 5 \\
			$x$ & 0 & 0 & 0 \\
			$b \pmod{6}$ & 1 & 3 & 5 \\
			\hline
	\end{tabular} .
\end{center}
The needed $\Diamond$-ITS$_3(2w+b+x;w,b)$s are from Lemmas \ref{2w+u,3}. The needed $\Diamond$-ITS$_3(b;u,0)$s are from Theorem \ref{ITS-one hole}.

When $ \lambda=6 $, the choices for $ (x,b \pmod 3) $ are
\begin{center}
\begin{tabular}{|c|l|l|l|}\hline
$v-2w \pmod{3}$ & 0 & 1& 2 \\
$x$ & 0 & 0& 0 \\
$b \pmod{3}$ & 0 & 1 & 2 \\
\hline
\end{tabular} .
\end{center}
The needed $\Diamond$-ITS$_6(2w+b+x;w,b)$ are from Lemma \ref{2w+u,6}. The needed $\Diamond$-ITS$_6(b;u,0)$s are from Theorem \ref{ITS-one hole}. \qed

\begin{Corollary}\label{cor:hole}
Let $\lambda\in\{2,3,6\}$. Let $w\geq u>1$. If $\lambda({v\choose 2}-{u\choose 2}-{w\choose 2}) \equiv 0\pmod3$, $\lambda (v-w) \equiv \lambda (v-u) \equiv \lambda (v-1) \equiv 0 \pmod 2$ and
\begin{displaymath}
v\geq \left\{ \begin{array}{ll}
w+4u+2, & \textrm{if $w<2u+1$},\\
2w+2u+1, & \textrm{if $w\geq 2u+1$},\\
\end{array} \right.
\end{displaymath}
then there exists a $\Diamond$-ITS$_\lambda(v;w,u)$.
\end{Corollary}

\proof If $w<2u+1$, apply Lemma \ref{th:hole1}. Otherwise, apply Lemma \ref{th:hole2}. \qed

\begin{Corollary}\label{w=u}
Let $ \lambda\in \{2,6\}$. If $ \lambda({v\choose 2}-2{u\choose 2})\equiv 0 \pmod 3$ and  $v\geq 3u$, then there exists a $\Diamond$-ITS$_\lambda(v;u,u)$ for any $ u\geq3 $.
\end{Corollary}

\proof The cases of $ v\geq 5u+2 $ and $ v\leq 5u $ are from Corollary \ref{cor:hole} and Lemma \ref{w=u small}, respectively. For $ v=5u+1 $, take a $ (3,\lambda) $-GDD of type $ u^5 1^1 $ (from Lemma \ref{GDD2}), and then place a TS$_\lambda(u)  $ (from Theorem \ref{ITS-one hole}) on three groups of size $ u $. \qed

\section{The remaining cases of $ \lambda=3 $}

When $\lambda=3$, by Lemma \ref{necessary}, $v$, $w$ and $u$ must be odd integers. Lemma \ref{2w+u,3} deals with the case of $ v=2w+u $. Combining Corollary \ref{cor:hole},  it suffices to consider
\begin{displaymath}
2w+u+2\leq v< \left\{ \begin{array}{ll}
w+4u+2, & \textrm{if $w<2u+1$},\\
2w+2u+1, & \textrm{if $w\geq 2u+1$},\\
\end{array} \right.
\end{displaymath}
and $vwu\equiv 1\pmod{2}$.

\begin{Lemma}\label{lem:pro1}
Let $s$ be a positive integer. Let $u$ and $w$ be odd positive integers. If $3\leq u\leq w$ and
\begin{displaymath}
w-\frac{1}{2}(u-1)\leq s< \left\{ \begin{array}{ll}
	\frac{w+1}{2}+u, & \textrm{if $w<2u+1$},\\
	w, & \textrm{if $w\geq 2u+1$},\\
\end{array} \right.
\end{displaymath}
then there exists a $\triangle$-decomposition of $3K_s-N$ where $N$ is a subgraph of $3 K_s$ satisfying that
\begin{itemize}
\item [$ (1) $] $N$ contains a subgraph $N'$ with $\epsilon(N')=\lceil \frac{3}{2}  u(s+\frac{1}{2}(u+1)-w)\rceil$ and $\chi'(N')\leq 3 u$;
\item[$ (2) $] $\epsilon(N)=\frac{3}{2} w(w-1)-3us+3\epsilon(N')$
\end{itemize}
except for $ (w,u,s)\in\{(3,3,2),(3,3,3),(3,3,4),(5,5,7)\}$.
\end{Lemma}

\proof Let $n'=\epsilon(N')= \lceil \frac{3}{2}  u\phi\rceil$ with
$$\phi=s+\frac{1}{2}(u+1)-w.$$ Let  $n=\epsilon(N)=\frac{3}{2} w(w-1)-3 us+3n'$.  Clearly, $\phi\geq 1$ since $s\geq w-\frac{1}{2}(u-1)$.

Firstly, we shall show that there exists a near-regular spanning subgraph $N$ of $3 K_s$ with $\epsilon(N)=n$ edges such that there exists a $\triangle$-decomposition of $3 K_s-N$. Write $f(s,n)=3{s\choose 2}-n$, which is the number of edges in $3 K_s-N$. By Lemma \ref{chromatic}, it suffices to show that $f(s,n)\equiv 0 \pmod 3$ and $0\leq f(s,n)\leq 3\lfloor\frac{s}{3}\lfloor \frac{3(s-1)}{2}\rfloor\rfloor$.
\begin{itemize}
\item By the definition of $f(s,n)$, we can see that $f(s,n) \equiv 0 \pmod 3$.
\item We now show that $f(s,n)\geq 0$. We have
	\begin{align}\label{L1-1}
	f(s,n) & =3{s\choose 2}-\frac{3}{2}w(w-1)+3 us-3\lceil \frac{3}{2}  u\phi\rceil\\\nonumber
	&\geq3{s\choose 2}-\frac{3}{2}w(w-1)+3 us-\frac{9}{2}  u\phi-\frac{3}{2}.
	\end{align}
	Note that $s=\phi-\frac{1}{2}(u+1)+w$. Then	
	$$f(s,n)\geq \frac{3}{8}\left[4\phi(\phi-2)+(8\phi+3u-4)(w-u)+u(w-4)-1\right],
	$$
	which is positive when $w\geq 5$ and $\phi\geq 2$ due to $ w\geq u\geq 3 $. For $ \phi=1 $, we have
	$$f(s,n)\geq \frac{3}{8}\left[(3u+4)(w-u)+u(w-4)-5\right],
	$$
	which is nonnegtive when $w\geq 5$. For $ (w,u)=(3,3) $, we have $ s\in\{2,3,4\} $ and $ \phi=s-1$. By substituting the values of $w$, $u $, $ s$ and $\phi$ in \eqref{L1-1}, we have $ f(s,n)\geq 0 $ except for $ (w,u,s)=(3,3,2) $.

	\item Next we show that $f(s,n)\leq 3\lfloor\frac{s}{3}\lfloor \frac{3(s-1)}{2}\rfloor\rfloor$. To do so, we evaluate
	\begin{align*}
3\lfloor\frac{s}{3}\lfloor \frac{3(s-1)}{2}\rfloor\rfloor-f(s,n)&\geq \frac{3 s(s-2)}{2}-3{s\choose 2}+\frac{3}{2} w(w-1)-3 us+\frac{9}{2}  u\phi\\
	& =-\frac{3s}{2}+\frac{3}{2}  w(w-1)-3 us+\frac{9}{2}  u\phi.
	\end{align*}
Recall that $s=\phi-\frac{1}{2}(u+1)+w$. Then
	\begin{align*}
	3s\lfloor \frac{s-1}{2}\rfloor-f(s,n)&\geq \frac{3}{2}\left((w-u-1)^2+\frac{u}{2}(\phi-1)+\phi(\frac{u}{2}-1)-\frac{1}{2}\right),
	\end{align*}
	which is positive since $ w\geq u\geq 3 $ and  $\phi\geq 1$.
\end{itemize}

Secondly, we need to show that $N$ contains a subgraph $N'$ satisfying $\epsilon(N')=n'$ and $\chi'(N')\leq 3 u$. We distinguish two cases.

CASE 1. $\Delta(N)\leq 3(u-1)$. We have
	\begin{align*}
	\epsilon(N)-n'&=n-n'\\
	& =  \frac{3}{2} w(w-1)-3 us+ 2 \lceil\frac{3}{2}u(s+\frac{1}{2}(u+1)-w)\rceil \\
	& \geq \frac{3}{2}  w(w-1)+\frac{3}{2} u(u+1)-3 uw\\
	&=\frac{3}{2} (w-u)(w-u-1)\geq 0.
	\end{align*}
Take any subgraph $ N' $ of $ N $ with $\epsilon(N')=n' $. By Lemma \ref{vizing}, $\chi'(N')\leq 3(u-1)+3=3 u $.

CASE 2. $\Delta(N)\geq 3(u-1)+1$. By Lemma \ref{vizing}, there exists a proper edge coloring of $N$ with $\Delta(N)+3$ colors. Let $N^*$ be a subgraph of $N$ with edge set given by the $3 u$ largest color classes. If we can show that  $\epsilon(N^*)\geq n'$, then there exists a subgraph $N'$  of $ N $ satisfying  $\chi'(N')\leq 3 u$.
	
From the definition of $N^*$, we have
\begin{align*}
	\epsilon(N^*)-n'& \geq\frac{3 u}{\Delta(N)+3}\epsilon(N)-n' = \frac{3 u}{\Delta(N)+3}\epsilon(N)-\lceil\frac{3}{2}   u\phi\rceil.
\end{align*}
Since $N$ is a near-regular spanning subgraph, for each vertex $x$ of $N$, we have $\deg_N(x)\geq\Delta(N)-2$. It follows that $\epsilon(N)\geq\frac{1}{2}s(\Delta(N)-2)$. Then
\begin{align*}
\epsilon(N^*)-n' &\geq \frac{3}{2} us\left(\frac{\Delta(N)-2}{\Delta(N)+3}\right)-\lceil\frac{3}{2}   u\phi\rceil.
\end{align*}
Since $\Delta(N)\geq3(u-1)+1$, we have $\frac{\Delta(N)-2}{\Delta(N)+3}\geq \frac{3u-4}{3 u+1}$. Thus
\begin{align}
	\epsilon(N^*)-n' &\geq \frac{3}{2} u\left(\frac{(3u-4)s}{3 u+1}-\phi\right) -\frac{1}{2}\nonumber\\
	&=\frac{3}{2} u\left(\frac{(3u-4)s}{3 u+1}-s-\frac{1}{2}(u+1)+w\right)-\frac{1}{2}\nonumber \\
	&=\frac{3}{2} u\left(-\frac{5s}{3 u+1}-\frac{1}{2}(u+1)+w\right)-\frac{1}{2}.\label{L2}
\end{align}

If $w\geq 2u+1$, then by assumption, $s<w$ and so
\begin{align*}
		\epsilon(N^*)-n'
		&>\frac{3}{2} u\left(-\frac{5w}{3 u+1}-\frac{1}{2}(u+1)+w\right)-\frac{1}{2} \\
		& = \frac{3}{2} u\left((1-\frac{5}{3 u+1})w-\frac{1}{2}(u+1)\right)-\frac{1}{2}\\
		&\geq \frac{3}{2} u\left((1-\frac{5}{3 u+1})(2u+1)-\frac{1}{2}(u+1)\right)-\frac{1}{2}\\
		&= \frac{3 u}{4(3 u+1)}\left(9u^2-14u-9\right)-\frac{1}{2},
\end{align*}	
which is positive since $ u\geq 3 $.

If $w<2u+1$, then by assumption, $s<\frac{1}{2}(w+1)+u$, i.e., $ s\leq\frac{1}{2}(w-1)+u$, and so
\begin{align*}
		\epsilon(N^*)-n'
		&\geq\frac{3}{2} u\left(-(\frac{1}{2}(w-1)+u)\frac{5}{3 u+1}-\frac{1}{2}(u+1)+w\right)-\frac{1}{2} \\
		& = \frac{3 u}{4(3 u+ 1)}\left((6u-3)w-3u^2-14u+4\right)-\frac{1}{2}\\
        & \geq \frac{3 u}{4(3 u+ 1)}\left((6u-3)u-3u^2-14u+4\right)-\frac{1}{2}\\
		& = \frac{3 u}{4(3 u+ 1)}\left(3u^2-17u+4\right)-\frac{1}{2},
\end{align*}
which is positive when $ u\geq7 $. Finally it suffices to do more careful calculation for $u\in\{3,5\}$. Recall that at this point, $u\leq w\leq 2u-1$ and  $w-\frac{u-1}{2}\leq s\leq\frac{1}{2}(w-1)+u$. If $u=3$, then $(w,s)\in\{(3,2),(3,3),(3,4),(5,4),(5,5)\}$. If $ u=5 $, then $(w,s)\in\{(5,3),(5,4),(5,5),(5,6),(5,7),(7,5),(7,$ $6),(7,7),(7,8),(9,7),(9,8),(9,9)\}$.
Using the inequality \eqref{L2}, one can check that $\epsilon(N^*)-n'\geq 0$ for all required $(w,u,s)$ except for $(w,u,s)=\{(3,3,2),(3,3,3),(3,3,4),(5,5,7)\}$.
\qed

Let $ \beta $ be an edge coloring of a multigraph $G$ with color set $C$. Let $x\in V (G)$ and $c \in C$. We use the notation $\eta_c(x)$ to denote the number of edges of $G$ that are incident with $x$ and colored $c$. If $\eta_c(x)\leq \lambda$ for each $x \in V (G)$ and each $c \in C$, then $ \beta $ is said to be {\em $\lambda$-proper} (cf. \cite{repacking}).

\begin{Proposition}\label{recolor}
If $ \beta $ is an equitable proper edge coloring of $ G $ with $ 3 k $ colors, then there  exists a $ 3$-proper edge coloring of $ G $ with $ k $ colors such that the sizes of any two color classes differ by at most $2$.
\end{Proposition}

\proof Let $\beta$ be an equitable proper edge coloring of $ G $ with color set $C=\{C_1,C_2,\ldots,C_x\}\cup\{C'_1,C'_2,\ldots,C'_y\}$ with $ x+y=3 k $, where $|C_i|=a$ for $1\leq i\leq x$, $|C'_j|=b$ for $1\leq j\leq y$, and $ |a-b|\leq 1 $. W.l.o.g., assume that $0\leq x\leq y \leq 3k$. For each $ 1\leq i\leq \lfloor\frac{x}{3}\rfloor $,  define two new color classes $ D_i= C_{1+3(i-1)} \cup C'_{1+3(i-1)}\cup C'_{2+3(i-1)} $ whose size is $ a+2b $, and $ D'_{i}=C_{2+3(i-1) } \cup C_{3+3(i-1)}\cup C'_{3+3(i-1)}$ whose size is $ 2a+b $.  We deal with the remaining color classes by merging three of them into one. Write the size of the new color classes as $ s_i $ for $ 1\leq i\leq k $. Then $ s_i\in\{a+2b,2a+b,3b\} $. So we get a $ 3$-proper edge coloring of $ G $ with $ k $ colors such that the sizes of any two color classes differ by at most $2$. \qed

\begin{Lemma}{\rm\cite[Definition 3.1 and Lemma 4.1]{repacking}}.\label{repacking} Let $\{U,S\}$ be a partition of a set $R$. Let $H_0$ be a subgraph of $\lambda K_R$. If there exists a $\triangle$-decomposition $\mathcal{B}_0$ of $H_0$, then there exists a subgraph $H$ of $\lambda K_R$ admitting a $\triangle$-decomposition $\mathcal{B}$ such that
\begin{itemize}
	\item $\epsilon(H)=\epsilon(H_0)$ and hence $|\mathcal{B}|=|\mathcal{B}_0|$;
    \item  for all $x,y\in U$, $xy\in E(H)$ if and only if $xy\in E(H_0)$;
	\item  for all $x,y,z\in U$, $\{x,y,z\}\in \mathcal{B}$ if and only if $\{x,y,z\}\in \mathcal{B}_0$;
	\item  for all $x\in U$, $\deg_H(x)=\deg_{H_0}(x)$;
	\item  for all $x,y\in S$, $|\deg_H(x)-\deg_H(y)|\leq 2$.
\end{itemize}	
\end{Lemma}

\begin{Lemma}\label{lem:pro2}
Let $s$ be a positive integer. Let $u$ and $w$ be odd positive integers. Let $U$ and $S$ be two disjoint sets with $|U|=u$ and $|S|=s$. If $3\leq u\leq w$ and
\begin{displaymath}
	w-\frac{1}{2}(u-1)\leq s< \left\{ \begin{array}{ll}
	\frac{w+1}{2}+u, & \textrm{if $w<2u+1$},\\
	w, & \textrm{if $w\geq 2u+1$},\\
	\end{array} \right.
\end{displaymath}
then there exists a $\triangle$-decomposition of $3 K_{U\cup S}-3 K_U-M$ where $M$ is a subgraph of $3 K_{U\cup S}-3 K_U$ such that
\begin{itemize}
    \item[$ (1) $] $\epsilon(M)=\frac{3}{2} w(w-1)$ and
    \item[$ (2) $] $3(2w-(u+s)-1)\leq\deg_M(x)\leq 3(w-1)$ for all $x\in U\cup S$
\end{itemize}
except for $ (w,u,s)\in\mathcal \{(3,3,2),(3,3,3),(3,3,4),(5,5,7)\}$.
\end{Lemma}

\proof Since $u,w$ and $s$ satisfy the conditions of Lemma \ref{lem:pro1}, there exists a $\triangle$-decomposition $\mathcal{A}$ of $3 K_S-N$ where $N$ is a subgraph of $3 K_S$ such that
\begin{itemize}
	\item [$ (1) $]$N$ contains a subgraph $N'$ with $\epsilon(N')=\lceil\frac{3}{2}  u(s+\frac{1}{2}(u+1)-w)\rceil$ and $\chi'(N')\leq 3 u$;
	\item [$ (2) $]$\epsilon(N)=\frac{3}{2} w(w-1)-3us+3\epsilon(N')$.
\end{itemize}
Since $\chi'(N')\leq 3 u$, by Lemma \ref{equi},  there exists an equitable proper edge coloring  of $N'$ with $ 3 u $ colors. It follows that by Proposition \ref{recolor}, there exists a $ 3 $-proper edge coloring $ \beta $ of $N'$ with color set $U$ such that the sizes of any two color classes differ by at most $2$.

Let $\mathcal{A}'=\{\{x,y,\beta(xy)\}:xy\in E(N')\}.$ Then $\mathcal{A}\cup \mathcal{A}'$ is a $\triangle$-decomposition of $3 K_{U\cup S}-3 K_U-M_0$ for some subgraph $M_0$ of $3 K_{U\cup S}-3 K_U$. Clearly $\epsilon(M_0)=\epsilon(N)-\epsilon(N')+3 us-2\epsilon(N')=\frac{3}{2} w(w-1)$. Note that for any $x,y\in U$, the numbers of edges in $N'$ colored by $ x $ and $ y $ differ by at most $2$, so $|\deg_{M_0}(x)-\deg_{M_0}(y)|\leq 4$.

Applying Lemma \ref{repacking} with $H_0=3 K_{U\cup S}-3 K_U-M_0$ and $\mathcal{B}_0=\mathcal{A}\cup \mathcal{A}'$, we get a new $\triangle$-decomposition $\cal B$ of $3 K_{U\cup S}-3 K_U-M$ where $M$ is a subgraph of $3 K_{U\cup S}-3 K_U$ such that
\begin{itemize}
\item [$ (1) $]$ \epsilon(M)=\epsilon(M_0) $;
\item [$ (2) $]$ \deg_M(x)=\deg_{M_0}(x) $ for all $x\in U$; and
\item [$ (3) $]$ |\deg_M(x)-\deg_M(y)|\leq 2 $ for all $x,y\in S$.
\end{itemize}
Thus $\epsilon(M)=\epsilon(M_0)=\frac{3}{2}  w(w-1)$. Note that since $\deg_M(x)=\deg_{M_0}(x)$ for all $x\in U$, we have $|\deg_M(x)-\deg_M(y)|\leq 4$ for all $x,y\in U$.

It remains to show that $3(2w-(u+s)-1)\leq\deg_M(x)\leq 3(w-1)$ for all $x\in U\cup S$. Firstly we can see that
\begin{align*}
\sum_{x\in U}\deg_M(x)&=\sum_{x\in U}\deg_{M_0}(x)=3 us-2\epsilon(N')\\
&=3 us-2\lceil\frac{3}{2} u(s+\frac{1}{2}(u+1)-w)\rceil,
\end{align*}
and
\begin{align*}
\sum_{x\in S}\deg_M(x)&=2\epsilon(M)-\sum_{x\in U}\deg_M(x) \\
&= 3 w(w-1)-3 us+2\lceil\frac{3}{2} u(s+\frac{1}{2}(u+1)-w)\rceil.
\end{align*}

We now show that $\deg_M(x)\leq 3(w-1)$ for all $x\in U\cup S$.
\begin{itemize}
\item If there were a vertex $ x\in U $ with $\deg_M(x)\geq 3w-2$, then $\deg_M(x)\geq 3w-6$ for any $ x\in U $ since $|\deg_M(x)-\deg_M(y)|\leq 4$ for all $x,y\in U$. Thus $\sum_{x\in U}\deg_M(x)\geq (u-1)(3 w-6)+3w-2=u(3w-6)+4$. On the other hand,
	\begin{align*}
	u(3 w-6)+4- \sum_{x\in U}\deg_M(x)&\geq u(3 w-6)+4 + 3 u(\frac{1}{2}(u+1)-w)\\
	&\geq \frac{3}{2} u(u-3)+4,
	\end{align*}
	which is positive since $u\geq 3$, a contradiction.
\item If there were a vertex $ x\in S $ with $\deg_M(x)\geq 3w-2$, then $\deg_M(x)\geq 3w-4$ for any  $ x\in S $ since $|\deg_M(x)-\deg_M(y)|\leq $2 for all $x,y\in S$. Thus $\sum_{x\in S}\deg_M(x)\geq (s-1)(3w-4)+3w-2=s(3w-4)+2$. On the other hand,
	\begin{align*}
	s(3 w-4)+2-\sum_{x\in S}\deg_M(x)&\geq s(3 w-4)-3 w(w-1)+3 uw-\frac{3}{2} u(u+1)+1\\
	&= 3(s-w+\frac{1}{2}(u-1))( w-3)+\frac{3}{2}u(w-u)\\
	&~~~-\frac{9}{2}w+3u+5s-\frac{7}{2}\\
		&\geq5(w-\frac{1}{2}(u-1))-\frac{9}{2}w+3u-\frac{7}{2}\\
		&=\frac{1}{2}(w+u)-1>0,
	\end{align*}
    where the second inequality follows from the fact that $s\geq w-\frac{1}{2}(u-1)$ and $ w\geq u\geq3 $. A contradiction occurs.
\end{itemize}

Finally we show that $\deg_M(x)\geq3(2w-(u+s)-1)$ for all $x\in U\cup S$.
\begin{itemize}
	\item If there were a vertex $ x\in U $ with $\deg_M(x)\leq 3(2w-(u+s)-1)-1$, then $\deg_M(x)\leq3(2w-(u+s)-1)+3$ for any  $ x\in U $ since $|\deg_M(x)-\deg_M(y)|\leq 4$ for all $x,y\in U$. Thus $\sum_{x\in U}\deg_M(x)\leq 3(u-1)(2w-(u+s))+3(2w-(u+s)-1)-1=3u(2w-(u+s))-4$. On the other hand,
	\begin{align*}
		&\sum_{x\in U}\deg_M(x)-3u(2w-(u+s))+4\\
		&=3 us-2\lceil\frac{3}{2} u(s+\frac{1}{2}(u+1)-w)\rceil-3u(2w-(u+s))+4\\
		&\geq \frac{3}{2} u(2s-2w+u-1)+3,
	\end{align*}
	which is positive since $s\geq w-\frac{1}{2}(u-1)$, a contradiction.

	\item Note that for $x\in S$, the existence of a $\triangle$-decomposition of $3 K_{U\cup S}-3 K_U-M$ implies $\deg_M(x)\equiv 3(u+s-1)\pmod 2$. If there were a vertex $x\in S$ with $\deg_M(x)\leq  3(2w-(u+s)-1)-2$, then $\deg_M(x)\leq3(2w-(u+s)-1)$ for any $ x\in S $ since $|\deg_M(x)-\deg_M(y)|\leq 2$ for all $x,y\in S$. Thus $\sum_{x\in S}\deg_M(x)<3 s(2w-(u+s)-1)$. On the other hand,		
	\begin{align*}
		&\sum_{x\in S}\deg_M(x)- 3 s(2w-(u+s)-1)\\
		&= 3 w(w-1)-3 us+2\lceil\frac{3}{2} u(s+\frac{1}{2}(u+1)-w)\rceil- 3s(2w-(u+s)-1)\\
		&\geq 3(s+\frac{1}{2}u-w)^2+\frac{3}{4}(u^2+2u-4w+4s),
	\end{align*}
	which is positive since $s\geq w-\frac{1}{2}(u-1)$, a contradiction. \qed
	\end{itemize}

For two vertex-disjoint multigraphs $G$ and $\lambda K_t$, denote by $G\vee\lambda K_t$ a new graph having vertex set $V(G)\cup V(K_t)$ and edge (multi-) set $E(G)\cup E(\lambda K_t)\cup \{\underline{\lambda} xy\colon x\in V(G),y\in V(K_t)\}$, where $\underline{\lambda} xy$ means that $xy$ occurs $\lambda$ times in the edge set.

\begin{Lemma}\label{laststep} \rm{\cite[Lemma 6.6]{repacking}}
Let $r$, $ t $ and $w$ be positive integers satisfying $w\leq r$ and $t=r-w+1$. Let $\lambda$ be a positive integer. Let $R$, $W$ and $T$ be disjoint sets with $|R|=r$, $|W|=w$, and $|T|=t$. Let $L$ be a bipartite multigraph with maximum multiplicity $ \lambda$ having bipartition $\{R,W\}$ satisfying
	\begin{itemize}
		\item[$ (a) $] $\deg_{L}(x)\equiv \lambda t \pmod 2$  for all $x\in R$;
		\item[$ (b) $] $\deg_{L}(x)=\lambda t$ for all $x\in W$; and
		\item[$ (c) $] $\deg_{L}(x)\leq \lambda t$ for all $x\in R$.
	\end{itemize}
	Then there exists a $\triangle$-decomposition of $L\vee \lambda K_t$.
\end{Lemma}

\begin{Lemma} \label{main result1 }
Let $u,w$ and $v$ be odd positive integers and $3\leq u\leq w$. If
\begin{displaymath}
2w+u+2\leq v< \left\{ \begin{array}{ll}
	w+4u+2,& \textrm{if $w<2u+1$,}\\
	2w+2u+1, & \textrm{if $w\geq 2u+1$,}\\
\end{array} \right.
\end{displaymath}
then there exists a $\triangle$-decomposition of $3 K_v-3 K_w-3 K_u$.
\end{Lemma}

\proof Let $s=\frac{1}{2}(v-1)-u$ and $t=\frac{1}{2}(v+1)-w$. So $s+t=v-u-w$ and $u+s+1=w+t$.  Let $ U $, $ W $, $ S $ and $ T $ be pairwise disjoint sets with $|U|=u$, $|W|=w$, $|S|=s$ and $|T|=t$. Let $ V=U\cup W\cup S\cup T $. We shall give a $\triangle$-decomposition of $3 K_V-3 K_W-3 K_U$.

Since $s=\frac{1}{2}(v-1)-u$, we have
\begin{displaymath}
w-\frac{1}{2}(u-1)\leq s< \left\{ \begin{array}{ll}
\frac{w+1}{2}+u, & \textrm{if $w<2u+1$},\\
w, & \textrm{if $w\geq 2u+1$}.\\
\end{array} \right.
\end{displaymath}
Then for $ (w,u,s)\notin \mathcal \{(3,3,2),(3,3,3),(3,3,4),(5,5,7)\}$, $w$, $u$ and $s$ satisfy the conditions in Lemma \ref{lem:pro2}. Thus there exists a $\triangle$-decomposition $\mathcal{B}_1$ of $3 K_{U\cup S}-3 K_U-M$ where $M$ is a subgraph of $3 K_{U\cup S}-3 K_U$ such that
$\epsilon(M)=\frac{3}{2}  w(w-1)$ and $3(2w-(u+s)-1)\leq\deg_M(x)\leq 3(w-1)$ for all $x\in U\cup S$.

Note that $\Delta(M)\leq 3(w-1)$, and so there exists an equitable proper edge coloring of $M$ with $ 3 w $ colors by lemmas \ref{vizing} and \ref{equi}. Also, due to $\epsilon(M)=\frac{3}{2}  w(w-1)$, each color class  has $\frac{1}{2} (w-1)$ edges. Merge every $ 3 $ color classes into a new one. Then we get a $ 3 $-proper edge cororing $\gamma  $ with color set $W$.  Let
$$\mathcal{B}_2=\{\{x,y,\gamma(xy)\}:xy\in E(M)\}.$$
Then $\mathcal{B}_1\cup\mathcal{B}_2$ is a $\triangle$-decomposition of $3 K_{U\cup W\cup S}-3 K_U-3 K_W-L$ where $L$ is some subgraph of $3 K_{U\cup W\cup S}-3 K_U-3 K_W$. Note that $L$ is a bipartite graph with bipartition $\{U\cup S,W\}$.

If there is a $\triangle$-decomposition of $L\vee 3 K_T$, then there exists a $\triangle$-decomposition of $3 K_V-3 K_U-3 K_W$. To get a $\triangle$-decomposition of $L\vee 3 K_T$, we use Lemma \ref{laststep} with $R=U\cup S$. Let $ |R|=r $. Then $ r=u+s $, $t=r-w+1$ (recall that $ u+s+1=w+t $) and $ w\leq r $ (since $ r=u+s=\frac{v-1}{2} $ and $ v\geq 2w+u+2 $). It remains to check that $L$ satisfies the conditions (a), (b) and (c) of Lemma \ref{laststep}.

\textbf{Condition (a).} Let $x\in R$. If $x\in U$, then the existence of a $\triangle$-decomposition of $3 K_{U\cup W\cup S}-3 K_U-3 K_W-L$ implies $\deg_L(x)\equiv 3(w+s) \pmod 2$. Since $ u\equiv1\pmod2 $, $\deg_L(x)\equiv 3(w+s+u+1) \equiv 3(u+s-w+1)\equiv3t \pmod 2$.  If $x\in S$, then $\deg_L(x)\equiv 3(u+s-1+w)\equiv3t \pmod 2$.

\textbf{Condition (b).} Let $x\in W$. Then $\deg_L(x)=3(u+s)-2\gamma_{\mathcal{B}_2}(x)$, where $\gamma_{\mathcal{B}_2}(x)$ is the number of triangles in $\mathcal{B}_2$ that contain $x$, i.e., the number of edges that colored by $x$ in $\gamma$. So $\gamma_{\mathcal{B}_2}(x)=\frac{3}{2}  (w-1)$. Thus $\deg_L(x)=3(u+s)-3(w-1)=3(u+s-w+1)=3 t$.

\textbf{Condition (c).} Let $x\in R$. Then $\deg_L(x)=\deg_M(x)+3 w-2\gamma_{\mathcal{B}_2}(x)$. Note that  $\gamma_{\mathcal{B}_2}(x)=\deg_M(x)$. Then $\deg_L(x)=3 w-\deg_M(x)$. Since $\deg_M(x)\geq 3(2w-(u+s)-1)$, it follows that $\deg_L(x)\leq 3(u+s-w+1)=3 t$.

Hence, by Lemma \ref{laststep}, there is a $\triangle$-decomposition  of $L\vee 3 K_T$.

Finally, it suffices to consider the cases of $ (w,u,s)\in \mathcal \{(3,3,2),(3,3,3),(3,3,4),(5,5,7)\}$. Since $ v=2(u+s)+1 $, we need to treat the cases of $ (w,u,v)\in\{(3,3,11),(3,3,13),(3,3,15),(5,5,25)\}$. For $ (w,u,v)=(3,3,11) $, we give an explicit construction in Appendix D.3. For $ (w,u,v)\in\{(3,3,13),(3,3,15)\}$, there exists a $\Diamond$-ITS$_3(v;w,u)$ by Theorem \ref{lambda=1}. For $ (w,u,v)=(5,5,25)$, start from a $ (3,3) $-GDD of type  $ 5^5 $ from Lemma \ref{GDD2}, and then fill in any three groups by using a $\Diamond$-ITS$_3(5)$ (from Theorem \ref{ITS-one hole}). \qed

\section{Conclusion}

Using the results of the previous sections, we are now ready to prove Theorems \ref{main result} and \ref{main result even}.

The necessity follows from Lemma \ref{necessary}. It suffices to consider the sufficiency. Theorem \ref{ITS-one hole} covers the case of $u=1$. The other cases are from Corollary \ref{cor:hole}, Lemmas \ref{2w+u,3} and \ref{main result1 }.



Finally, we remark that a $\Diamond $-ITS$_\lambda (v;w,u)$ is equivalent to a $ (3,\lambda )$-GDD of type $ u^1w^11^{v-w-u} $. Therefore, Theorems \ref{main result} and \ref{main result even} provide new existence results for GDDs of such type.

\appendix

\section{Appendix: $ \Diamond $-ITS$_\lambda (v;w,u,y;z_1,z_2 ) $s in Lemma $4.6$}\label{app lattic}

Denote by $I_{n}=\{0,1,\ldots,n-1\}$. $ A\cup A\cup \cdots \cup A$ ($h$ times) is denoted by $ \underline h A $.

\vspace{0.25cm}

\noindent A $ \Diamond $-ITS$_2 (11;4,3,2;1,0) $ on $I_{11}$ with holes $ \{0,1,2,3\} $, $\{5,6,7\}$ and $ \{3,4\} $:
\begin{center}\tabcolsep 0.065in
\begin{longtable}{llllllll}
\underline{2}\{3, 5, 10\}& \underline{2}\{3, 6, 9\}& \{1, 6, 10\}& \{2, 7, 9\}&\{0, 4, 7\}&\{0, 7, 10\}& \{0, 8, 10\}& \{1, 4, 6\}\\
\underline{2}\{2, 4, 5\}&\underline{2}\{3, 7, 8\}&\{2, 6, 8\}&\{2, 6, 10\}&
\{0, 4, 6\}&\{0, 6, 8\}&\{2, 7, 10\}&\{2, 8, 9\}\\
\underline{2}\{1, 5, 8\}&\underline{2}\{0, 5, 9\}&
\{1, 4, 7\}&\{1, 7, 9\}&\{1, 9, 10\}&\{4, 8, 9\}&\{4, 8, 10\}&\{4, 9, 10\}
\end{longtable}
\end{center}
A $ \Diamond $-ITS$_2 (13;4,3,4;1,0) $ on $I_{13}$ with holes $ \{0,1,2,3\} $, $ \{7,8,9\} $ and $ \{3,4,5,6\} $:
\begin{center}\tabcolsep 0.055in
\begin{longtable}{llllllll}
\underline{2}\{3, 7, 10\}&\underline{2}\{3, 8, 11\}&\underline{2}\{3, 9, 12\}&\underline{2}\{0, 5, 8\}&\{0, 11, 12\}&\{2, 9, 10\}&\{1, 6, 10\}&\{0, 10, 12\}\\
\underline{2}\{6, 7, 12\}&\underline{2}\{1, 5, 7\}&\underline{2}\{0, 4, 7\}&\{1, 9, 10\}&\{6, 10, 11\}&\{1, 4, 9\}&\{2, 5, 10\}&\{0, 10, 11\}\\
\underline{2}\{4, 8, 10\}&\underline{2}\{2, 7, 11\}&\underline{2}\{1, 8, 12\}&\{1, 4, 11\}&\{1, 6, 11\}&\{2, 4, 9\}&\{2, 4, 12\}&\{4, 11, 12\}\\
\underline{2}\{5, 9, 11\}&\underline{2}\{0, 6, 9\}&\underline{2}\{2, 6, 8\}&\{2, 5, 12\}&\{5, 10, 12\}
\end{longtable}
\end{center}
A $\Diamond$-ITS$_2(14;5,4,5;2,3)$  on $I_{14}$ with holes $ \{0,1,2,3,4\} $, $ \{5,6,7,8\} $ and $ \{3,4,5,6,7\} $:
\begin{center}\tabcolsep 0.04in
\begin{longtable}{llllllll}
		\underline{2}\{0, 5, 9\}& \underline{2}\{0, 6, 10\}&\underline{2}\{0, 7, 11\}&\{0, 8, 12\}&\{0, 8, 13\}&\{0, 12, 13\}&\{1, 5, 10\}&\{1, 5, 11\}\\\{1, 6, 9\}&\{1, 6, 11\}&\{1, 7, 12\}&\{1, 7, 13\}&\{1, 8, 9\}&\{1, 8, 10\}&\{1, 12, 13\}&\underline{2}\{2, 5, 12\}\\\{5, 10, 13\}&\{5, 11, 13\}&\{6, 9, 12\}&\{6, 11, 12\}&\underline{2}\{2, 6, 13\}&\{2, 7, 9\}&\{2, 7, 10\}&\{7, 10, 12\}\\\{2, 8, 10\}&\{7, 9, 13\}&\{2, 8, 11\}&\{2, 9, 11\}&\{3, 8, 9\}&\{3, 8, 11\}&\{4, 8, 13\}&\{4, 8, 12\}\\\{3, 9, 12\}&\{3, 10, 12\}&\{4, 11, 12\}&\{3, 10, 13\}&\{4, 9, 13\}&\{3, 11, 13\}&\{4, 9, 10\}&\{4, 10, 11\}\\\{9, 10, 11\}
\end{longtable}
\end{center}
A $\Diamond$-ITS$_2(18;7,4,5;2,3)$  on $I_{18}$ with holes $ \{0,1,\ldots,6\} $, $ \{7,8,9,10\} $ and $ \{5,6,7,8,9\} $:
\begin{center}\tabcolsep 0.02in
\begin{longtable}{llllllll}
		\underline{2}\{0, 7, 11\}&\underline{2}\{0, 8, 12\}&\underline{2}\{0, 9, 13\}&\underline{2}\{0, 10, 14\}&\{0, 15, 16\}&\{0, 15, 17\} &
		\{0, 16, 17\}&\underline{2}\{1, 7, 12\}\\\underline{2}\{1, 9, 11\}&\underline{2}\{1, 10, 13\}&\underline{2}\{1, 8, 14\}&\{1, 15, 16\}&
		\{1, 15, 17\}&\{1, 16, 17\}&\underline{2}\{2, 7, 13\}&\underline{2}\{3, 7, 15\}\\\underline{2}\{4, 7, 16\}&\underline{2}\{7, 14, 17\}&\underline{2}\{2, 8, 11\}&\underline{2}\{4, 8, 15\}&\{3, 8, 16\}&\{8, 13, 16\}&\{3, 8, 17\}&\{8, 13, 17\}\\\{2, 9, 17\}&\{2, 12, 17\} &\{2, 10, 15\}&\{2, 10, 16\}&\{2, 9, 15\}&\{2, 12, 14\}&\{2, 14, 16\}& \{9, 12, 15\}\\\{4, 9, 14\}&\{4, 9, 17\}&\{3, 9, 14\}&\{3, 9, 16\} &\{9, 12, 16\}&\{3, 10, 11\}&\{3, 10, 12\}&\{3, 13, 14\}\\\{3, 11, 12\} &\{3, 13, 17\}&\{4, 10, 11\}&\{4, 10, 17\}&\{4, 11, 13\}&\{4, 12, 13\}&
		\{4, 12, 14\}&\underline{2}\{5, 11, 17\}\\\{6, 10, 17\}&\{6, 12, 17\} &\{5, 10, 12\}&\{5, 10, 16\}&\{6, 10, 15\}&\{5, 12, 13\}&\underline{2}\{5, 14, 15\}&\{5, 13, 16\}\\\{11, 12, 15\}&\{6, 12, 16\}&\{6, 13, 15\}&\{11, 13, 15\}&\{6, 13, 14\}&\{6, 11, 14\}&\{6, 11, 16\}&\{11, 14, 16\}
	\end{longtable}
\end{center}
A $\Diamond$-ITS$_2(23;8,7,5;2,3)$  on $I_{23}$ with holes  $ \{0,1,\ldots,7\} $, $ \{8,9,\ldots,14\} $ and $ \{6,7,8,9,10\} $:
\begin{center}\tabcolsep 0.05in
\begin{longtable}{lllllll}
\underline{2}\{0, 8, 15\} & \underline{2}\{0, 9, 16\}&\underline{2}\{0, 10, 17\}&\underline{2}\{0, 11, 18\}&\underline{2}\{0, 12, 19\}&\underline{2}\{0, 13, 20\}&\underline{2}\{0, 14, 21\}\\
\{0, 21, 22\}& \underline{2}\{1, 8, 16\}&\underline{2}\{1, 10, 15\}&\underline{2}\{1, 11, 17\}&\underline{2}\{1, 12, 18\}&\underline{2}\{1, 13, 19\}&\underline{2}\{1, 9, 20\}\\
\{1, 14, 21\}&\{1, 14, 22\}&\{1, 21, 22\}&\underline{2}\{2, 8, 17\}&\underline{2}\{2, 9, 15\}&\underline{2}\{2, 12, 16\}&\underline{2}\{2, 13, 18\}\\
\underline{2}\{2, 14, 19\}&\{2, 10, 20\}&\{2, 10, 21\}&\{2, 11, 21\}&\{2, 11, 22\}&\{2, 20, 22\}&\underline{2}\{3, 8, 18\}\\
\underline{2}\{3, 9, 17\}&\underline{2}\{3, 14, 15\}&\underline{2}\{3, 12, 20\}&\underline{2}\{3, 13, 16\}&\{3, 10, 19\}&\{3, 10, 21\}&\{3, 11, 21\}\\
\{3, 11, 22\}&\{3, 19, 22\}&\underline{2}\{4, 8, 19\}&\{5, 8, 21\}&\{8, 20, 21\}&\{5, 8, 22\}&\{8, 20, 22\}\\
\underline{2}\{4, 9, 18\}&\{5, 9, 21\}&\{9, 19, 21\}&\{5, 9, 22\}&\{9, 19, 22\}&\underline{2}\{5, 12, 15\}&\underline{2}\{5, 13, 17\}\\
\underline{2}\{4, 10, 16\}&\underline{2}\{10, 18, 22\}&\underline{2}\{4, 11, 15\}&\underline{2}\{4, 13, 21\}&\{4, 12, 17\}&\{4, 12, 22\}&\{4, 17, 22\}\\
\underline{2}\{4, 14, 20\}&\{5, 10, 19\}&\{5, 10, 20\}&\{5, 14, 16\}&\{5, 14, 18\}&\{6, 13, 15\}&\{6, 13, 22\}\\
\{7, 13, 15\}&\{7, 13, 22\}&\{5, 11, 16\}&\{5, 11, 19\}&\{5, 18, 20\}&\{6, 11, 16\}&\{6, 11, 19\}\\
\underline{2}\{7, 11, 20\}&\{6, 12, 17\}&\{6, 12, 21\}&\{7, 12, 21\}&\{7, 12, 22\}&\{6, 14, 17\}&\{6, 14, 18\}\\
\{6, 16, 22\}&\{7, 14, 16\}&\{7, 14, 17\}&\{15, 16, 22\}&\{15, 17, 22\}&\{6, 15, 20\}&\{6, 18, 21\}\\
\{6, 19, 20\}&\{7, 15, 18\}&\{7, 16, 19\}&\{7, 17, 21\}&\{7, 18, 19\}&\{15, 16, 21\}&\{15, 17, 19\}\\
\{15, 18, 21\}&\{15, 19, 20\}&\{17, 18, 19\}&\{16, 19, 21\}&\{17, 20, 21\}&\{16, 17, 18\}&\{16, 17, 20\}\\
\{16, 18, 20\}
\end{longtable}
\end{center}
A $\Diamond$-ITS$_2(23;8,7,8;2,6)$  on $I_{23}$ with holes  $ \{0,1,\ldots,7\} $, $\{8,9,\ldots,14\} $ and  $ \{6,7,\ldots,13\} $:
\begin{center}\tabcolsep 0.045in
\begin{longtable}{lllllll}
\underline{2}\{0, 8, 15\}&\underline{2}\{0, 9, 16\}&\underline{2}\{0, 10, 17\}&\underline{2}\{0, 11, 18\}&\underline{2}\{0, 12, 19\}&\underline{2}\{0, 13, 20\}&\{0, 14, 21\}\\
\{0, 14, 22\}&\{0, 21, 22\}&\underline{2}\{1, 8, 16\}&\underline{2}\{1, 10, 15\}&\underline{2}\{1, 11, 17\}&\underline{2}\{1, 12, 18\}&\underline{2}\{1, 13, 19\}\\
\underline{2}\{1, 9, 20\}&\underline{2}\{1, 14, 21\}&\{1, 21, 22\}&\underline{2}\{2, 8, 17\}&\underline{2}\{2, 9, 15\}&\underline{2}\{2, 12, 16\}&
		\underline{2}\{2, 13, 18\}\\
\underline{2}\{2, 14, 19\}&\{2, 10, 20\}&\{2, 10, 21\}&\{2, 11, 21\}&\{2, 11, 22\}&\{2, 20, 22\} &
		\underline{2}\{3, 8, 18\}\\
\underline{2}\{3, 9, 17\}&\underline{2}\{3, 14, 15\}&\underline{2}\{3, 12, 20\}&\underline{2}\{3, 13, 16\}&	\{3, 10, 19\}&\{3, 10, 21\}&\{3, 11, 21\}\\
\{3, 11, 22\}&\{3, 19, 22\} &
		\underline{2}\{4, 8, 19\}&\{5, 8, 21\}&\{8, 20, 21\}&\{5, 8, 22\}&
		\{8, 20, 22\}\\
\underline{2}\{4, 9, 18\}&\{5, 9, 21\}&\{9, 19, 21\} &
		\{5, 9, 22\}&\{9, 19, 22\}&\underline{2}\{5, 12, 15\}&\underline{2}\{5, 13, 17\} \\
		\underline{2}\{4, 12, 21\}&\underline{2}\{12, 17, 22\}&\underline{2}\{13, 15, 21\}&\underline{2}\{4, 13, 22\} &
		\underline{2}\{4, 10, 16\}&\underline{2}\{10, 18, 22\}&\{4, 11, 15\}\\
\{4, 11, 20\}&\underline{2}\{4, 14, 17\}&
		\{4, 15, 20\}&\{5, 10, 19\}&\{5, 10, 20\}&\{6, 15, 22\}&\{6, 16, 22\} \\
		\{7, 15, 22\}&\{7, 16, 22\}&\{5, 11, 16\}&\{5, 11, 20\}&\{5, 14, 16\}&
		\{5, 14, 18\}&\{5, 18, 19\}\\
\{11, 15, 19\}&\{11, 16, 19\} &
		\underline{2}\{6, 14, 20\}&\{7, 14, 16\}&\{7, 14, 18\}&\{6, 15, 16\} &
		\{6, 17, 19\}\\
\{6, 17, 21\}&\{6, 18, 19\}&\{6, 18, 21\}&\{7, 15, 19\} &
		\{7, 18, 20\}&\{7, 17, 20\}&\{7, 17, 21\}\\
\{7, 19, 21\}&\{16, 19, 20\}&
		\{17, 19, 20\}&\{15, 16, 17\}&\{15, 17, 18\}&\{15, 18, 20\}&
		\{16, 17, 18\}\\\{16, 18, 21\}&\{16, 20, 21\}
	\end{longtable}
\end{center}
A $ \Diamond $-ITS$_2 (24;8,6,6;2,0) $ on $I_{24}$ with holes $ \{0,1,\ldots,7\} $, $ \{12,13,\ldots,17\}$ and $ \{6,7,\ldots,11\} $:
\begin{center}\tabcolsep 0.048in
\begin{longtable}{lllllll}
\underline{2}\{6, 12, 18\}& \underline{2}\{6, 13, 19\} &\underline{2}\{6, 14, 20\} &\underline{2}\{6, 15, 21\} &\underline{2}\{6, 16, 22\}&\underline{2}\{6, 17, 23\} &\underline{2}\{7, 12, 19\} \\
\underline{2}\{7, 14, 18\} &\underline{2}\{7, 15, 20\} &\underline{2}\{7, 16, 21\} &\underline{2}\{7, 17, 22\}  &\underline{2}\{7, 13, 23\}  &\underline{2}\{0, 8, 12\} &\underline{2}\{0, 9, 13\}\\
\underline{2}\{0, 10, 14\} &
		\underline{2}\{0, 11, 15\}&\underline{2}\{0, 16, 18\}&\underline{2}\{0, 17, 19\} &
		\underline{2}\{0, 20, 21\}&\underline{2}\{0, 22, 23\}&\underline{2}\{1, 9, 12\}\\
\underline{2}\{1, 8, 14\}&\underline{2}\{1, 10, 15\}&
		\underline{2}\{1, 11, 16\}&\underline{2}\{1, 17, 18\}&\underline{2}\{1, 13, 20\}&\underline{2}\{1, 19, 22\}&\underline{2}\{1, 21, 23\}\\
		\underline{2}\{2, 10, 12\}&\underline{2}\{2, 8, 13\}&\underline{2}\{2, 9, 16\}&\underline{2}\{2, 11, 17\}&\underline{2}\{2, 14, 19\} &
		\underline{2}\{2, 18, 21\}&\underline{2}\{2, 15, 22\}\\
\underline{2}\{2, 20, 23\}&\underline{2}\{3, 12, 23\}&\underline{2}\{3, 9, 14\} &
		\underline{2}\{3, 8, 15\}&\underline{2}\{3, 10, 16\}&\underline{2}\{3, 17, 20\}&\{3, 11, 19\}\\
\{3, 18, 19\}&\{3, 11, 18\} &
		\{3, 13, 21\}&\{3, 13, 22\}&\{3, 21, 22\}&\underline{2}\{4, 8, 17\}&
		\{5, 9, 17\}\\
\{9, 17, 21\}&\{5, 10, 17\}&\{10, 17, 21\}&\{4, 16, 19\} &
		\{4, 16, 23\}&\{5, 8, 21\}&\{8, 19, 21\}\\
\{8, 18, 22\}&\{8, 18, 22\} &
		\{8, 16, 20\}&\{8, 19, 20\}&\{5, 8, 23\}&\{8, 16, 23\}&\{5, 16, 19\} \\
		\{5, 16, 20\}&\{5, 14, 23\}&\{11, 14, 23\}&\{4, 14, 21\}&\{4, 14, 22\} &
		\{5, 14, 21\}&\{11, 14, 22\}\\
\{5, 12, 20\}&\{5, 12, 22\}&\{5, 9, 15\} &
		\{5, 15, 19\}&\{5, 10, 18\}&\{5, 13, 18\}&\{5, 11, 13\}\\
\{5, 11, 22\} &
		\{9, 21, 22\}&\{4, 9, 22\}&\{12, 20, 22\}&\{10, 13, 22\} &
		\{10, 20, 22\}&\{4, 11, 12\}\\
\{4, 12, 21\}&\{11, 12, 21\} &
		\{10, 19, 21\}&\{11, 13, 21\}&\{4, 10, 13\}&\{4, 13, 18\}&\{4, 9, 20\} \\
		\{9, 18, 20\}&\{4, 15, 18\}&\{4, 15, 19\}&\{4, 10, 23\}&\{4, 11, 20\} &
		\{9, 15, 23\}&\{15, 18, 23\}\\
\{9, 18, 19\}&\{9, 19, 23\} &
		\{10, 18, 23\}&\{10, 19, 20\}&\{11, 18, 20\}&\{11, 19, 23\}
	\end{longtable}
\end{center}

\section{Appendix: $(3,2)$-eframes of type $h^n m^1$ in Lemma 4.11}\label{app,eframe}

A $(3,2)$-eframe of type $3^5 1^1$ is constructed on  $I_{16}$ with the group set $\{\{i,5+i,10+i\}:0\leq i\leq4\}\cup\{\{15\}\}$. Only $14$ base blocks are listed:
	\begin{center}
		\begin{tabular}{llllll}
			$\{4,11,3\}$,
			&$\{13,4,6\}$,
			&$\{15,14,2\}$,
			&$\{8,6,7\}$,
			&$\{4,15,1\}$,
			&$\{7,14,3\}$,
			\\$\{9,11,12\}$,
			&$\{11,2,8\}$,
			&$\{1,12,3\}$,
			&$\{9,2,3\}$,&
			$\{1,13,14\}$,
			&$\{4,8,12\}$,
			\\$\{6,9,15\}$,
			&$\{7,11,13\}$.
		\end{tabular}
	\end{center}
	Let
	$\alpha=(0~\,3\cdots~12)(1~\,4\cdots~13)(2~\,5\cdots~14)$
	be a permutation on $I_{16}$. Let $G$ be the group generated by
	$\alpha$. All other blocks are obtained by developing these base
	blocks under the action of $G$. The first base block generates a 1-partial parallel class missing $\{15\}$.
	The other $13$ base blocks form a $3$-partial parallel class missing $\{0,5,10\}$, which generates $5$ $3$-partial parallel classes under the action of $G$.
	
\vskip16pt

\noindent A $(3,2)$-eframe of type $3^5 2^1$ is constructed on $(Z_5\times I_3)\cup \{\infty_1,\infty_2\}$ with the group set $\{\{i_0,i_1,i_2\}:0\leq i\leq4\}\cup\{\{\infty_1,\infty_2\}\}$. Only $16$ base blocks are listed:
\begin{center}
	\begin{tabular}{llllll}
		$\{1_2,0_1,3_0\}$,
		&$\{2_1,3_0,1_2\}$,
		&$\{1_1,2_1,3_0\}$,
		&$\{2_2,\infty_2,3_1\}$,
		&$\{\infty_1,2_0,1_0\}$,&
		$\{3_2,4_0,1_0\}$,\\
		$\{1_1,4_0,3_2\}$,
		&$\{2_1,4_1,\infty_1\}$,
		&$\{2_0,3_0,4_2\}$,
		&$\{1_2,3_2,2_2\}$,&
		$\{2_1,1_0,4_0\}$,
		&$\{4_2,1_2,\infty_1\}$,
		\\$\{4_2,3_0,\infty_2\}$,
		&$\{3_1,2_0,\infty_2\}$,
		&$\{3_1,4_1,1_2\}$,&
		$\{4_1,2_2,1_1\}$.
	\end{tabular}
\end{center}
All other blocks are obtained by developing these base
blocks by $(+1,-)$ modulo $(5,-)$. Each of the first $2$ base blocks generates a 1-partial parallel class missing $\{\infty_1,\infty_2\}$. The other $14$ base blocks form a $3$-partial parallel class missing $\{0_0,0_1,0_2\}$, which generates $5$ $3$-partial parallel classes by $(+1,-)$ modulo $(5,-)$.

\vskip16pt

\noindent A $(3,2)$-eframe of type $3^6 1^1$ is constructed on $I_{19}$ with the group set $\{\{i,6+i,12+i\}:0\leq i\leq5\}\cup\{\{18\}\}$. Only $17$ base blocks are listed:
	\begin{center}
		\begin{tabular}{llllll}
			$\{13,4,6\}$,
			&$\{16,15,5\}$,
			&$\{1,15,4\}$,
			&$\{2,16,18\}$,
			&$\{17,18,9\}$,
			&$\{9,10,14\}$,\\
			$\{15,11,13\}$,
			&$\{17,4,8\}$,
			&$\{3,11,2\}$,
			&$\{1,3,8\}$,
			&$\{1,5,14\}$,
			&$\{2,7,18\}$,\\
			$\{3,10,13\}$,
			&$\{4,5,8\}$,
			&$\{7,9,14\}$,
			&$\{7,10,11\}$,
			&$\{13,16,17\}$.
		\end{tabular}
	\end{center}
	Let
	$\alpha=(0~\,1\cdots~5)(6~\,7\cdots~11)(12~\,13\cdots~17)$
	be a permutation on $I_{19}$. Let $G$ be the group generated by
	$\alpha$. All other blocks are obtained by developing these base
	blocks under the action of $G$. The first base block generates a 1-partial parallel class missing $\{18\}$.
	The other $16$ base blocks form a $3$-partial parallel class missing $\{0,6,12\}$, which generates $6$ $3$-partial parallel classes under the action of $G$.

\vskip16pt

\noindent A $(3,2)$-eframe of type $3^6 2^1$ is constructed on $(Z_6\times I_3)\cup \{\infty_1,\infty_2\}$ with the group set $\{\{i_0,i_1,i_2\}:0\leq i\leq4\}\cup\{\{\infty_1,\infty_2\}\}$. Only $19$ base blocks are listed:
\begin{center}
	\begin{tabular}{llllll}
		$\{3_0,2_2,5_1\}$,
		&$\{3_2,5_0,0_1\}$,
		&$\{4_0,2_0,5_0\}$,
		&$\{1_1,2_2,4_0\}$,
		&$\{1_1,2_1,4_0\}$,&
		$\{5_1,3_1,1_2\}$,\\
		$\{2_2,3_2,4_2\}$,
		&$\{2_1,4_1,3_0\}$,
		&$\{5_1,4_2,\infty_1\}$,
		&$\{1_1,\infty_1,5_0\}$,
		&$\{1_0,2_0,3_2\}$,
		&$\{1_0,3_0,\infty_2\}$,\\
		$\{1_0,4_2,\infty_1\}$,
		&$\{2_0,3_2,5_2\}$,
		&$\{3_0,2_2,5_2\}$,&
		$\{5_0,3_1,4_1\}$,
		&$\{2_1,5_1,1_2\}$,
		&$\{3_1,1_2,\infty_2\}$,\\
		$\{4_1,5_2,\infty_2\}$.
	\end{tabular}
\end{center}
All other blocks are obtained by developing these base
blocks by $(+1,-)$ modulo $(6,-)$. Each of the first $2$ base blocks generates a 1-partial parallel class missing $\{\infty_1,\infty_2\}$. The other $17$ base blocks form a $3$-partial parallel class missing $\{0_0,0_1,0_2\}$, which generates $6$ $3$-partial parallel classes by $(+1,-)$ modulo $(6,-)$.

\vskip16pt

\noindent A $(3,2)$-eframe of type $6^5 5^1$ is constructed on $(Z_{10}\times I_3)\cup \{\infty_1,\infty_2,\ldots,\infty_5\}$ with the group set $\{\{i_0,i_1,i_2,(i+5)_0,(i+5)_1.(i+5)_2\}:0\leq i\leq4\}\cup\{\{\infty_1,\infty_2,\ldots,\infty_5\}\}$. Only $34$ base blocks are listed:
\begin{center}
	\begin{tabular}{llllll}
		$\{7_2,9_1,5_0\}$,
		&$\{6_2,7_1,3_0\}$,
		&$\{4_1,7_0,3_2\}$,
		&$\{8_2,0_1,4_0\}$,
		&$\{0_2,8_1,1_0\}$,
		&$\{8_0,7_2,\infty_1\}$,\\
		$\{3_0,9_2,\infty_4\}$,
		&$\{7_2,\infty_1,9_0\}$,
		&$\{6_2,\infty_3,3_2\}$,
		&$\{1_2,\infty_4,9_1\}$,
		&$\{\infty_1,6_1,7_1\}$,
		&$\{9_1,1_0,2_2\}$,\\
		$\{9_0,3_0,6_0\}$,
		&$\{\infty_4,7_0,8_1\}$,
		&$\{4_1,2_1,6_1\}$,
		&$\{1_2,8_1,9_2\}$,
		&$\{7_1,3_1,6_1\}$,
		&$\{6_0,9_2,2_1\}$,\\
		$\{2_1,8_2,6_2\}$,
		&$\{3_2,7_1,4_2\}$,
		&$\{4_2,3_1,7_2\}$,
		&$\{2_2,\infty_2,6_2\}$,
		&$\{3_2,\infty_5,2_2\}$,
		&$\{1_1,\infty_5,4_1\}$,\\
		$\{1_2,4_0,3_0\}$,
		&$\{9_1,7_0,1_0\}$,
		&$\{1_0,2_0,\infty_5\}$,
		&$\{8_2,4_0,6_0\}$,
		&$\{8_2,7_0,4_2\}$,
		&$\{2_0,4_0,3_1\}$,\\
		$\{2_0,4_1,\infty_2\}$,
		&$\{1_1,8_0,\infty_2\}$,
		&$\{1_1,8_0,\infty_3\}$,
		&$\{9_0,8_1,\infty_3\}$.
	\end{tabular}
\end{center}
All other blocks are obtained by developing these base
blocks by $(+1,-)$ modulo $(10,-)$. Each of the first $5$ base blocks generates a 1-partial parallel class missing $\{\infty_1,\infty_2,\ldots,\infty_5\}$.
The other $29$ base blocks generate a $6$-partial parallel class missing $\{0_0,0_1,0_2,5_0,5_1,5_2\}$ by $(+5,-)$ modulo $(10,-)$, which generates $5$ $6$-partial parallel classes by $(+1,-)$ modulo $(10,-)$.

\vskip16pt

\noindent A $(3,2)$-eframe of type $6^5 8^1$ is constructed on $(Z_{10}\times I_3)\cup \{\infty_1,\infty_2,\ldots,\infty_8\}$ with the group set $\{\{i_0,i_1,i_2,(i+5)_0,(i+5)_1.(i+5)_2\}:0\leq i\leq4\}\cup\{\{\infty_1,\infty_2,\ldots,\infty_8\}\}$. Only $40$ base blocks are listed:
\begin{center}
	\begin{tabular}{llllll}
		$\{2_0,0_1,1_2\}$,
		&$\{2_2,1_0,4_1\}$,
		&$\{6_2,8_0,4_1\}$,
		&$\{6_2,3_0,4_1\}$,
		&$\{1_0,3_2,2_1\}$,&
		$\{5_2,6_0,8_1\}$,
		\\$\{5_2,1_0,9_1\}$,
		&$\{5_1,9_2,6_0\}$,
		&$\{\infty_4,4_1,1_1\}$,
		&$\{4_1,2_0,\infty_1\}$,&
		$\{\infty_7,9_2,2_1\}$,
		&$\{4_0,8_0,\infty_4\}$,
		\\$\{3_2,\infty_1,7_2\}$,
		&$\{4_1,\infty_1,7_0\}$,
		&$\{8_0,\infty_6,6_0\}$,&
		$\{3_1,6_2,\infty_2\}$,
		&$\{7_0,1_0,3_2\}$,
		&$\{8_1,9_0,1_0\}$,
		\\$\{4_0,\infty_5,1_2\}$,
		&$\{2_2,3_1,9_2\}$,&
		$\{9_0,3_2,\infty_5\}$,
		&$\{6_1,4_2,\infty_8\}$,
		&$\{6_1,2_0,3_0\}$,
		&$\{\infty_5,1_1,9_1\}$,
		\\$\{\infty_6,2_1,1_1\}$,&
		$\{7_1,1_0,\infty_2\}$,
		&$\{\infty_4,4_2,2_2\}$,
		&$\{9_1,7_1,3_1\}$,
		&$\{6_2,\infty_6,8_2\}$,
		&$\{9_0,6_2,\infty_7\}$,\\
		$\{9_1,6_1,\infty_3\}$,
		&$\{8_1,7_1,1_2\}$,
		&$\{8_1,\infty_8,7_2\}$,
		&$\{4_0,7_0,\infty_8\}$,
		&$\{7_2,8_2,1_2\}$,&
		$\{8_2,9_2,\infty_3\}$,
		\\$\{2_0,3_0,\infty_3\}$,
		&$\{3_0,6_0,4_2\}$,
		&$\{2_1,8_0,\infty_7\}$,
		&$\{2_2,6_0,\infty_2\}$.
	\end{tabular}
\end{center}
All other blocks are obtained by developing these base
blocks by $(+1,-)$ modulo $(10,-)$. Each of the first $8$ base blocks generates a 1-partial parallel class missing $\{\infty_1,\infty_2,\ldots,\infty_8\}$.
The other $32$ base blocks generate a $6$-partial parallel class missing $\{0_0,0_1,0_2,5_0,5_1,5_2\}$ by $(+5,-)$ modulo $(10,-)$, which generates $5$ $6$-partial parallel classes by $(+1,-)$ modulo $(10,-)$.

\section{Appendix: $(3,6)$-GDDs of type $ 2^{s-1}1^1 $ in Lemma 4.14}

\noindent A $ (3,6) $-GDD of type $ 2^4 1^1 $ is constructed on $(Z_4\times I_2)\cup \{\infty\}$ with the group set $\{\{i_0,i_1\}:0\leq i\leq4\}\cup\{\{\infty\}\}$, which has 1 3-partial parallel class associated with the group of size 1, and 1 6-partial parallel class for each group of size 2. Only $16$ base blocks are listed:
\begin{center}
	\begin{tabular}{llllll}
		$\{0_1,1_0,3_1\}$,
		&$\{0_0,1_1,3_0\}$,
		&$\{1_0,3_0,\infty\}$,
		&$\{3_0,2_1,1_0\}$,
		&$\{\infty,3_1,2_1\}$,&
		$\{3_0,1_1,2_1\}$,
		\\$\{1_1,2_1,3_1\}$,
		&$\{3_0,1_1,2_0\}$,
		&$\{2_1,1_0,\infty\}$,
		&$\{1_1,2_1,3_0\}$,&
		$\{1_0,2_0,3_0\}$,
		&$\{1_0,2_0,3_1\}$,
		\\$\{1_0,2_0,\infty\}$,
		&$\{2_0,1_1,3_1\}$,
		&$\{2_0,3_1,\infty\}$,&
		$\{1_1,3_1,\infty\}$.
	\end{tabular}
\end{center}
All other blocks are obtained by developing these base
blocks by $(+1,-)$ modulo $(4,-)$. The first $2$ base blocks generate a $3$-partial parallel class missing $\{\infty\}$ by $(+1,-)$ modulo $(4,-)$.
The other $14$ base blocks form a $6$-partial parallel class missing $\{0_0,0_1\}$, which generate $4$ holey $6$-partial parallel classes.

\vskip16pt

\noindent A $ (3,6) $-GDD of type $ 2^5 1^1 $ is constructed on $(Z_5\times I_2)\cup \{\infty\}$ with the group set $\{\{i_0,i_1\}:0\leq i\leq4\}\cup\{\{\infty\}\}$, which has 1 3-partial parallel class associated with the group of size 1, and 1 6-partial parallel class for each group of size 2.  Only $20$ base blocks are listed:
\begin{center}
	\begin{tabular}{llllll}
		$\{4_0,0_1,3_1\}$,
		&$\{4_1,0_0,2_0\}$,
		&$\{4_0,\infty,2_0\}$,
		&$\{4_1,\infty,2_1\}$,
		&$\{1_1,4_0,3_0\}$,&
		$\{2_0,1_1,3_1\}$,
		\\$\{1_0,3_1,4_1\}$,
		&$\{1_0,3_1,4_0\}$,
		&$\{4_1,3_0,1_1\}$,
		&$\{4_1,3_1,2_1\}$,&
		$\{1_0,4_1,3_1\}$,
		&$\{1_0,\infty,4_1\}$,
		\\$\{2_0,1_1,3_1\}$,
		&$\{3_0,2_0,4_0\}$,
		&$\{2_1,3_0,1_0\}$,&
		$\{1_0,2_1,\infty\}$,
		&$\{2_0,3_0,4_0\}$,
		&$\{2_0,3_0,\infty\}$,
		\\$\{1_1,4_0,2_1\}$,
		&$\{1_1,2_1,\infty\}$.
	\end{tabular}
\end{center}
All other blocks are obtained by developing these base
blocks by $(+1,-)$ modulo $(5,-)$. The first $2$ base blocks generate a 3-partial parallel class missing $\{\infty\}$ by $(+1,-)$ modulo $(5,-)$.
The other $18$ base blocks form a $6$-partial parallel class missing $\{0_0,0_1\}$, which generates $5$ $6$-partial parallel classes.

\section{Appendix: ITSs of small orders}

\subsection{$\Diamond$-ITS$_2(2w+u+x,w,u) $s in  Lemma 4.24}

A $ \Diamond $-ITS$_2 (16;6,3) $ on $I_{16}$ with holes $ \{0,1,\ldots,5\} $ and $ \{6,7,8\} $:
\begin{center}\tabcolsep 0.05in
\begin{longtable}{lllllll}
		\underline{2}\{0, 6, 9\}& \underline{2}\{0, 7, 10\}& \underline{2}\{0, 8, 11\}&
		\underline{2}\{0, 12, 13\}& \underline{2}\{0, 14, 15\}&
		\underline{2}\{1, 6, 10\}&\underline{2}\{1, 8, 9\}\\\underline{2}\{1, 7, 12\}&\underline{2}\{1, 11, 14\}&\underline{2}\{1, 13, 15\}&
		\underline{2}\{2, 6, 11\}&\underline{2}\{2, 7, 9\}& \underline{2}\{2, 8, 13\}&
		\{2, 10, 14\}\\ \{2, 12, 14\}& \{2, 10, 15\}& \{2, 12, 15\}&
		\underline{2}\{3, 6, 12\}& \{9, 11, 12\} &\underline{2}\{5, 6, 15\} &\{5, 11, 12\} \\ \{4, 6, 13\}&
		\{4, 6, 14\} &\{6, 13, 14\}& \{3, 8, 10\} &\{3, 8, 15\} &\underline{2}\{5, 7, 13\} &\{5, 9, 12\}
		\\ \{3, 7, 14\}& \{4, 7, 14\}& \{3, 7, 11\} &\{4, 7, 15\}
		&\{7, 11, 15\}& \{3, 9, 14\}& \{9, 13, 14\} \\\{3, 10, 13\}& \{3, 11, 13\} &
		\{3, 9, 15\}& \{4, 8, 10\} &\{4, 8, 12\} &\{8, 12, 15\} &\underline{2}\{5, 8, 14\}\\ \{5, 10, 11\}
		& \{10, 12, 14\}& \{4, 10, 12\}& \{4, 11, 13\}& \{9, 10, 13\} &
		\{4, 9, 11\}& \{4, 9, 15\}\\\{10, 11, 15\} &\{5, 9, 10\}
	\end{longtable}
\end{center}
A $ \Diamond $-ITS$_2 (22;9,3) $ on $I_{22}$ with holes  $ \{0,1,\ldots,8\} $ and $ \{9,10,11\} $:
\begin{center}\tabcolsep 0.05in
	\begin{longtable}{lllllll}
		\underline{2}\{0, 9, 12\}&\underline{2}\{0, 10, 13\}&\underline{2}\{0, 11, 14\}& \underline{2}\{0, 15, 16\}&\underline{2}\{0, 17, 18\}&
		\{0, 19, 20\}&\{0, 19, 21\}\\\{0, 20, 21\}&\underline{2}\{1, 9, 13\}&
		\underline{2}\{1, 11, 12\}&\underline{2}\{1, 10, 15\}&\underline{2}\{1, 14, 17\} &\underline{2}\{1, 16, 18\}
		&\{1, 19, 20\}\\\{1, 19, 21\}&
		\{1, 20, 21\}&\underline{2}\{2, 9, 14\}&\underline{2}\{2, 10, 12\}&
		\underline{2}\{2, 11, 16\}&\underline{2}\{2, 13, 19\}&\{2, 15, 17\}\\
		\{2, 17, 20\}&\{2, 15, 20\}&\underline{2}\{2, 18, 21\}&\underline{2}\{3, 9, 15\}
		&\underline{2}\{3, 11, 13\}&\underline{2}\{3, 10, 18\}&\underline{2}\{3, 12, 19\}\\\{3, 14, 16\}&\{3, 14, 20\}&\{3, 16, 20\}
		&\underline{2}\{3, 17, 21\}&\underline{2}\{4, 9, 16\}&\underline{2}\{5, 9, 21\} &\underline{2}\{6, 9, 17\}\\\{7, 9, 18\}&\{7, 9, 19\}&
		\{8, 9, 19\}&\{8, 9, 20\}&\{9, 18, 20\}&\underline{2}\{4, 10, 14\}
		&\underline{2}\{4, 11, 19\}\\\underline{2}\{4, 12, 21\}&\underline{2}\{4, 13, 17\} &\{4, 15, 18\}&\{4, 15, 20\}&\{4, 18, 20\}&\underline{2}\{5, 10, 16\}&\underline{2}\{6, 10, 19\}\\\underline{2}\{7, 10, 21\}&\{8, 10, 17\}&\{8, 10, 20\}&\{10, 17, 20\}&\underline{2}\{5, 11, 20\}&\underline{2}\{5, 12, 17\}&\{5, 13, 14\}\\\{5, 13, 18\}&\{5, 14, 19\}&
		\{5, 15, 18\}&\{5, 15, 19\}&\{6, 11, 18\}&\{6, 11, 21\}&\{7, 11, 18\}\\
		\{7, 11, 17\}&\{8, 11, 17\}&\{8, 11, 15\}&\{11, 15, 21\}&
		\{15, 17, 19\}&\{7, 16, 17\}&\{16, 17, 19\}\\\{7, 16, 19\}&
		\{8, 18, 19\}&\{14, 18, 19\}&\{6, 12, 18\}&\underline{2}\{7, 12, 15\}&\{7, 13, 14\}&\{7, 13, 20\}\\\{7, 14, 20\}&\{8, 13, 18\}&\{12, 14, 18\}&\{6, 13, 15\}&\{6, 14, 15\}&\{8, 13, 15\}&\{14, 15, 21\}\\\{8, 12, 14\}& \{8, 12, 16\}&\{6, 14, 16\}&\{8, 14, 21\}&\{8, 16, 21\}&\{6, 13, 21\}&
		\{13, 16, 21\}\\\{6, 12, 20\}&\{6, 16, 20\}&\{12, 13, 16\}&
		\{12, 13, 20\}
	\end{longtable}
\end{center}
A $ \Diamond $-ITS$_2 (25;9,6) $ on $I_{25}$ with holes  $ \{0,1,\ldots,8\} $ and $ \{9,10,\ldots,14\} $:
\begin{center}\tabcolsep 0.08in
	\begin{longtable}{lllllll}
		$\{0,9,17\}$
		&$\{0,9,18\}$
		&$\{0,10,17\}$
		&$\{0,10,15\}$
		&$\{0,11,19\}$&$\{0,11,24\}$
		&$\{0,12,19\}$
		\\$\{0,12,21\}$
		&$\{0,13,20\}$
		&$\{0,13,24\}$
		&$\{0,14,20\}$
		&$\{0,14,18\}$
		&$\{0,15,23\}$
		&$\underline{2}\{0,16,22\}$\\$\{0,21,23\}$
		&$\{1,9,19\}$
		&$\{1,9,16\}$
		&$\{1,10,20\}$
		&$\{1,10,23\}$&$\{1,11,16\}$
		&$\{1,11,18\}$
		\\$\{1,12,19\}$
		&$\{1,12,22\}$
		&$\{1,13,20\}$&$\{1,13,23\}$&$\{1,14,24\}$
		&$\{1,14,15\}$
		&$\{1,15,21\}$
		\\$\{1,17,24\}$&$\{1,17,22\}$
		&$\{1,18,21\}$
		&$\{2,9,19\}$
		&$\{2,9,24\}$
		&$\{2,10,15\}$&$\{2,10,16\}$
		\\$\{2,11,23\}$
		&$\{2,11,21\}$
		&$\{2,12,15\}$
		&$\{2,12,22\}$&$\{2,13,23\}$&$\{2,13,16\}$
		&$\{2,14,21\}$
		\\$\{2,14,20\}$
		&$\{2,17,20\}$&$\{2,17,24\}$
		&$\{2,18,22\}$
		&$\{2,18,19\}$
		&$\{3,9,17\}$
		
		&$\{3,9,22\}$\\$\{3,10,18\}$
		&$\{3,10,20\}$
		&$\{3,11,17\}$
		&$\{3,11,24\}$
		&$\{3,12,16\}$&
		$\{3,12,23\}$
		&$\{3,13,22\}$
		\\
\end{longtable}
\end{center}
A $ \Diamond $-ITS$_2 (23;8,6) $ on $I_{23}$ with holes  $ \{0,1,\ldots,7\} $ and $ \{8,9,\ldots,13\} $:
\begin{center}\tabcolsep 0.05in
	\begin{longtable}{lllllll}
		$\{0,8,16\}$
		&$\{0,8,17\}$
		&$\{0,9,14\}$
		&$\{0,9,20\}$
		&$\{0,10,14\}$&$\{0,10,20\}$
		
		&$\{0,11,18\}$
		\\$\{0,11,21\}$
		&$\{0,12,18\}$
		&$\{0,12,16\}$&
		$\underline{2}\{0,13,22\}$&$\{0,15,19\}$
		&$\{0,15,21\}$
		&$\{0,17,19\}$\\
		$\{1,8,21\}$
		&$\{1,8,20\}$
		
		&$\{1,9,14\}$
		&$\{1,9,15\}$
		&$\{1,10,17\}$&
		$\{1,10,14\}$
		&$\{1,11,21\}$
		\\$\{1,11,19\}$
		&$\{1,12,20\}$
		&$\{1,12,22\}$&
		$\{1,13,17\}$
		&$\{1,13,18\}$
		&$\{1,15,22\}$
		&$\{1,16,19\}$
		
		\\$\{1,16,18\}$&
		$\{2,8,16\}$
		&$\{2,8,15\}$
		&$\{2,9,19\}$
		&$\{2,9,21\}$
		&$\{2,10,21\}$&
		
		$\{2,10,18\}$
		\\$\{2,11,20\}$
		&$\{2,11,22\}$
		&$\{2,12,20\}$
		&$\{2,12,22\}$&
		$\{2,13,16\}$
		&$\{2,13,14\}$
		&$\{2,14,17\}$
		\\$\{2,15,19\}$
		&$\{2,17,18\}$&
		$\{3,8,14\}$
		
		&$\{3,8,14\}$
		&$\{3,9,20\}$
		&$\{3,9,16\}$
		&$\{3,10,22\}$\\
		$\{3,10,15\}$
		&$\{3,11,17\}$
		
		&$\{3,11,15\}$
		&$\{3,12,18\}$
		&$\{3,12,21\}$&
		$\{3,13,19\}$
		&$\{3,13,17\}$
		\\
\end{longtable}
\end{center}
A $ \Diamond $-ITS$_2 (32;11,9) $ on $I_{32}$ with holes  $ \{0,1,\ldots,10\} $ and $ \{11,12,\ldots,19\} $:

\begin{center}\tabcolsep 0.05in
	\begin{longtable}{lllllll}
		$\underline{2}\{0,11,20\}$
		&$\{0,12,22\}$
		&$\{0,12,21\}$
		&$\{0,13,21\}$&
		$\{0,13,26\}$ &$\{0,14,29\}$
		&$\{0,14,30\}$
		\\$\{0,15,27\}$
		&$\{0,15,24\}$&
		$\{0,16,28\}$
		
		&$\{0,16,24\}$
		&$\{0,17,25\}$
		&$\{0,17,31\}$
		&$\{0,18,27\}$\\
		$\{0,18,30\}$
		
		&$\{0,19,23\}$
		&$\{0,19,25\}$
		&$\{0,22,23\}$
		&$\{0,26,29\}$&
		$\{0,28,31\}$
		
		&$\{1,11,29\}$
		\\$\{1,11,21\}$
		&$\{1,12,25\}$
		&$\{1,12,28\}$&
		$\{1,13,24\}$
		
		&$\{1,13,23\}$
		&$\{1,14,23\}$
		&$\{1,14,26\}$
		\\$\{1,15,21\}$&
		$\{1,15,24\}$
		
		&$\{1,16,25\}$
		&$\{1,16,22\}$
		&$\{1,17,27\}$
		&$\{1,17,30\}$&
		$\underline{2}\{1,18,31\}$
		
		\\$\{1,19,28\}$
		&$\{1,19,29\}$
		&$\{1,20,27\}$&
		$\{1,20,26\}$
		&$\{1,22,30\}$
		
		&$\{2,11,26\}$
		&$\{2,11,27\}$
		\\$\{2,12,21\}$&
		$\{2,12,20\}$
		&$\{2,13,28\}$
		
		&$\{2,13,21\}$
		&$\{2,14,24\}$
		&$\{2,14,25\}$&
		$\{2,15,22\}$
		\\$\{2,15,30\}$
		
		&$\{2,16,23\}$
		&$\{2,16,27\}$
		&$\{2,17,23\}$&
		$\{2,17,30\}$
		&$\{2,18,25\}$
		
		&$\{2,18,24\}$
		\\
\end{longtable}
\end{center}
A $ \Diamond $-ITS$_2 (24;9,4) $ on $I_{24}$ with holes  $ \{0,1,\ldots,8\} $ and $ \{9,10,11,12\} $:
\begin{center}\tabcolsep 0.05in
	\begin{longtable}{lllllll}
		$\{0,9,14\}$
		&$\{0,9,21\}$
		&$\{0,10,15\}$
		&$\{0,10,14\}$
		&$\{0,11,17\}$&$\{0,11,18\}$
		
		&$\{0,12,23\}$
		\\$\{0,12,21\}$
		&$\{0,13,19\}$
		&$\{0,13,22\}$&$\{0,15,19\}$
		&$\{0,16,23\}$
		&$\{0,16,20\}$
		&$\{0,17,20\}$
		\\$\{0,18,22\}$&$\{1,9,15\}$
		
		&$\{1,9,13\}$
		&$\{1,10,23\}$
		&$\{1,10,15\}$
		&$\{1,11,17\}$&$\{1,11,13\}$
		\\$\{1,12,19\}$
		
		&$\{1,12,16\}$
		&$\{1,14,21\}$
		&$\{1,14,20\}$&$\{1,16,18\}$
		&$\{1,17,22\}$
		&$\{1,18,21\}$
		
		\\$\{1,19,23\}$
		&$\{1,20,22\}$&$\{2,9,13\}$
		&$\{2,9,17\}$
		&$\{2,10,21\}$
		&$\{2,10,14\}$
		
		&$\{2,11,21\}$\\$\{2,11,15\}$
		&$\{2,12,15\}$
		&$\{2,12,22\}$
		&$\{2,13,17\}$
		&$\{2,14,20\}$&$\{2,16,19\}$
		&$\{2,16,23\}$
		\\$\{2,18,20\}$
		&$\{2,18,23\}$
		&$\{2,19,22\}$&$\{3,9,22\}$
		&$\{3,9,19\}$
		&$\{3,10,13\}$
		&$\{3,10,23\}$
		\\$\{3,11,20\}$&$\{3,11,13\}$
		&$\{3,12,22\}$
		&$\{3,12,18\}$
		&$\{3,14,23\}$
		&$\{3,14,21\}$&$\{3,15,17\}$
		\\
\end{longtable}
\end{center}
A $ \Diamond $-ITS$_2 (31;10,9) $ on $I_{31}$ with holes  $ \{0,1,\ldots,9\} $ and $ \{10,11,\ldots,18\} $:
\begin{center}\tabcolsep 0.04in
	\begin{longtable}{lllllll}
		$\{0,10,19\}$
		&$\{0,10,21\}$
		
		&$\underline{2}\{0,11,22\}$
		&$\{0,12,19\}$&$\{0,12,20\}$
		&$\{0,13,21\}$
		&$\{0,13,20\}$
		\\$\underline{2}\{0,14,24\}$
		&$\{0,15,27\}$
		&$\{0,15,25\}$
		&$\underline{2}\{0,16,28\}$
		&$\{0,17,29\}$&$\{0,17,26\}$
		&$\{0,18,29\}$
		\\$\{0,18,23\}$
		&$\{0,23,27\}$
		&$\{0,25,30\}$&$\{0,26,30\}$
		&$\{1,10,20\}$
		&$\{1,10,24\}$
		&$\{1,11,19\}$
		\\$\{1,11,21\}$&$\{1,12,20\}$
		&$\{1,12,21\}$
		&$\{1,13,27\}$
		&$\{1,13,26\}$
		&$\{1,14,27\}$&$\{1,14,22\}$
		\\$\{1,15,25\}$
		&$\{1,15,26\}$
		&$\{1,16,30\}$
		&$\{1,16,25\}$&$\{1,17,19\}$
		&$\{1,17,28\}$
		&$\{1,18,30\}$
		\\$\{1,18,29\}$
		&$\{1,22,28\}$&$\{1,23,24\}$
		&$\{1,23,29\}$
		&$\{2,10,22\}$
		&$\{2,10,19\}$
		&$\{2,11,19\}$\\$\{2,11,21\}$
		&$\{2,12,29\}$
		&$\{2,12,24\}$
		&$\{2,13,22\}$
		&$\{2,13,20\}$&$\{2,14,21\}$
		&$\{2,14,23\}$
		\\$\{2,15,23\}$
		&$\{2,15,26\}$
		&$\{2,16,27\}$&$\{2,16,25\}$
		&$\{2,17,28\}$
		&$\{2,17,26\}$
		&$\{2,18,27\}$
		\\
\end{longtable}
\end{center}
A $ \Diamond $-ITS$_2 (11;3,2) $ on $I_{11}$ with holes  $ \{0,1,2\} $ and $ \{3,4\} $:

\begin{center}\tabcolsep 0.05in
	\begin{longtable}{llllllll}
		\underline{2}\{0, 3, 5\}&\underline{2}\{0, 4, 6\}  &\underline{2}\{0, 7, 8\}&\underline{2}\{0, 9, 10\}&\{1, 3, 9\}&\{1, 3, 10\}&
		\{3, 6, 7\}&\{3, 7, 10\}\\\{2, 3, 6\}&\{2, 3, 8\}&\{3, 8, 9\} &
		\{1, 4, 7\}&\{1, 4, 10\}&\{1, 5, 7\}&\{1, 5, 6\}&\{1, 6, 8\}\\
		\{1, 8, 9\}&\{2, 4, 8\}&\{4, 5, 8\}&\{5, 8, 10\}&\{6, 8, 10\} &
		\{2, 4, 9\}&\{4, 7, 9\}&\{4, 5, 10\}\\\{2, 6, 10\}&\{2, 7, 10\}&
		\{2, 5, 7\}&\{2, 5, 9\}&\{5, 6, 9\}&\{6, 7, 9\}
	\end{longtable}
\end{center}
A $ \Diamond $-ITS$_2 (19;6,4) $ on $I_{19}$ with holes  $ \{0,1,\ldots,5\} $ and $ \{6,7,8,9\} $:
\begin{center}\tabcolsep 0.05in
	\begin{longtable}{lllllll}
		$\{0,6,10\}$
		&$\{0,6,11\}$
		&$\{0,7,16\}$
		&$\{0,7,18\}$
		&$\underline{2}\{0,8,17\}$&$\{0,9,13\}$
		&$\{0,9,11\}$
		\\$\{0,10,13\}$
		&$\{0,12,16\}$&$\{0,12,15\}$
		&$\{0,14,18\}$
		&$\{0,14,15\}$
		
		&$\{1,6,10\}$
		&$\{1,6,12\}$\\$\{1,7,12\}$
		&$\{1,7,14\}$
		&$\{1,8,14\}$
		&$\{1,8,16\}$
		&$\{1,9,18\}$&$\{1,9,15\}$
		&$\{1,10,18\}$
		\\$\{1,11,15\}$
		&$\{1,11,13\}$
		&$\{1,13,17\}$&$\{1,16,17\}$
		&$\{2,6,16\}$
		&$\{2,6,15\}$
		&$\underline{2}\{2,7,13\}$
		\\$\{2,8,15\}$
		&$\{2,8,12\}$
		&$\{2,9,17\}$
		&$\{2,9,14\}$
		&$\{2,10,17\}$&$\{2,10,18\}$
		&$\{2,11,16\}$
		
		\\$\{2,11,14\}$
		&$\{2,12,18\}$
		&$\{3,6,18\}$&$\{3,6,11\}$
		&$\{3,7,16\}$
		&$\{3,7,18\}$
		
		&$\{3,8,13\}$
		\\$\{3,8,15\}$&$\{3,9,14\}$
		&$\{3,9,17\}$
		&$\{3,10,17\}$
		&$\{3,10,12\}$
		
		&$\{3,11,14\}$&$\{3,12,15\}$
		\\$\{3,13,16\}$
		&$\{4,6,12\}$
		&$\{4,6,15\}$
		&$\underline{2}\{4,7,10\}$
		&$\{4,8,18\}$
		&$\{4,8,11\}$
		&$\{4,9,12\}$
		\\$\{4,9,18\}$&$\{4,11,16\}$
		&$\{4,13,17\}$
		
		&$\{4,13,16\}$
		&$\{4,14,17\}$
		&$\{4,14,15\}$&$\underline{2}\{5,6,17\}$
		\\$\{5,7,14\}$
		&$\{5,7,15\}$
		
		&$\underline{2}\{5,8,10\}$&$\{5,9,11\}$
		&$\{5,9,13\}$
		&$\{5,11,12\}$
		&$\{5,12,13\}$\\$\{5,14,18\}$
		
		&$\{5,15,16\}$
		&$\{5,16,18\}$
		&$\underline{2}\{6,13,14\}$
		&$\{6,16,18\}$
		&$\{7,11,12\}$
		&$\{7,11,17\}$
		\\$\{7,15,17\}$
		&$\{8,11,18\}$&$\{8,12,13\}$
		&$\{8,14,16\}$
		&$\{9,10,15\}$
		
		&$\{9,10,16\}$
		&$\{9,12,16\}$\\$\{10,11,13\}$
		&$\{10,11,15\}$
		&$\{10,12,14\}$
		&$\{10,14,16\}$
		&$\{11,17,18\}$&$\{12,14,17\}$
		&$\{12,17,18\}$
		\\$\underline{2}\{13,15,18\}$
		&$\{15,16,17\}$
	\end{longtable}
\end{center}
A $ \Diamond $-ITS$_2 (21;7,4) $ on $I_{21}$ with holes  $ \{0,1,\ldots,6\} $ and $ \{7,8,9,10\} $:

\begin{center}\tabcolsep 0.05in
	\begin{longtable}{lllllll}
		$\{0,7,12\}$
		&$\{0,7,19\}$
		&$\{0,8,20\}$
		&$\{0,8,14\}$
		&$\{0,9,15\}$&
		$\{0,9,19\}$
		&$\{0,10,20\}$
		\\$\{0,10,18\}$
		&$\{0,11,17\}$
		&$\{0,11,18\}$&
		$\{0,12,15\}$&$\{0,13,17\}$
		&$\{0,13,16\}$
		&$\{0,14,16\}$
		\\$\{1,7,20\}$&
		$\{1,7,18\}$
		&$\{1,8,14\}$
		&$\{1,8,16\}$
		
		&$\{1,9,20\}$
		&$\{1,9,17\}$&
		$\{1,10,13\}$
		\\$\{1,10,18\}$
		&$\{1,11,16\}$
		&$\{1,11,19\}$
		
		&$\{1,12,15\}$&
		$\{1,12,14\}$
		&$\{1,13,19\}$
		&$\{1,15,17\}$
		\\$\{2,7,14\}$
		&$\{2,7,18\}$
		&
		$\{2,8,11\}$
		&$\{2,8,12\}$
		&$\{2,9,17\}$
		&$\{2,9,18\}$
		&$\{2,10,14\}$\\
		$\{2,10,20\}$
		
		&$\{2,11,15\}$&$\{2,12,20\}$
		&$\{2,13,15\}$
		&$\{2,13,17\}$&
		$\underline{2}\{2,16,19\}$&$\{3,7,11\}$
		
		\\$\{3,7,15\}$
		&$\{3,8,18\}$&
		$\{3,8,13\}$
		&$\{3,9,14\}$
		&$\{3,9,11\}$
		&$\{3,10,15\}$
		
		&$\{3,10,17\}$\\$\{3,12,18\}$&$\{3,12,19\}$&$\{3,13,16\}$&$\{3,14,19\}$&$\{3,16,20\}$
		&$\{3,17,20\}$
		&$\{4,7,11\}$
		\\
\end{longtable}
\end{center}
A $ \Diamond $-ITS$_2 (23;8,4) $ on $I_{23}$ with holes  $ \{0,1,\ldots,7\} $ and $ \{8,9,10,11\} $:
\begin{center}\tabcolsep 0.05in
	\begin{longtable}{lllllll}
		
		$\{0,8,14\}$
		&$\{0,8,17\}$
		&$\{0,9,13\}$
		&$\{0,9,20\}$
		&$\{0,10,13\}$&
		$\{0,10,18\}$
		
		&$\{0,11,19\}$\\$\{0,11,16\}$&$\{0,12,17\}$&$\{0,12,21\}$&
		$\{0,14,19\}$
		&$\{0,15,20\}$
		
		&$\{0,15,22\}$&$\{0,16,21\}$
		\\$\{0,18,22\}$&
		$\{1,8,16\}$
		&$\{1,8,15\}$
		&$\{1,9,19\}$
		
		&$\{1,9,21\}$&$\{1,10,18\}$&
		$\{1,10,21\}$
		\\$\{1,11,13\}$
		&$\{1,11,15\}$
		&$\{1,12,20\}$
		
		&$\{1,12,16\}$&$\{1,13,22\}$
		&$\{1,14,17\}$
		&$\{1,14,18\}$
		\\$\{1,17,19\}$
		&$\{1,20,22\}$
		&$\{2,8,20\}$
		&$\{2,8,17\}$
		&$\{2,9,16\}$
		&$\{2,9,12\}$
		&$\{2,10,20\}$\\
		$\{2,10,21\}$
		
		&$\{2,11,15\}$
		&$\{2,11,22\}$&$\{2,12,15\}$&$\{2,13,21\}$&$\{2,13,19\}$&$\{2,14,18\}$
		
		\\$\{2,14,19\}$&$\{2,16,22\}$
		&$\{2,17,18\}$&$\{3,8,19\}$
		&$\{3,8,21\}$
		&$\{3,9,21\}$
		
		&$\{3,9,20\}$
		\\$\{3,10,17\}$&
		$\{3,10,14\}$
		&$\{3,11,12\}$
		&$\{3,11,13\}$
		&$\{3,12,18\}$
		
		&$\{3,13,18\}$&$\{3,14,22\}$
		\\
\end{longtable}
\end{center}
A $ \Diamond $-ITS$_2 (26;8,7) $ on $I_{26}$ with holes  $ \{0,1,\ldots,7\} $ and $ \{8,9,\ldots,14\} $:
\begin{center}\tabcolsep 0.05in
	\begin{longtable}{lllllll}
		$\{0,8,15\}$
		&$\{0,8,23\}$
		&$\{0,9,16\}$
		&$\{0,9,25\}$
		&$\{0,10,16\}$&
		$\{0,10,18\}$
		&$\{0,11,23\}$\\$\{0,11,15\}$
		&$\{0,12,19\}$
		&$\{0,12,20\}$&
		$\{0,13,21\}$&$\{0,13,19\}$
		
		&$\{0,14,24\}$&$\{0,14,20\}$\\$\{0,17,21\}$&$\{0,17,25\}$&$\{0,18,22\}$&$\{0,22,24\}$
		
		&$\{1,8,22\}$&$\{1,8,15\}$&$\{1,9,20\}$\\$\{1,9,25\}$&$\{1,10,15\}$
		&$\{1,10,17\}$
		
		&$\{1,11,18\}$&$\{1,11,25\}$&$\{1,12,19\}$&$\{1,12,21\}$\\$\{1,13,24\}$&$\{1,13,23\}$
		&$\{1,14,23\}$&$\{1,14,24\}$&$\{1,16,20\}$&$\{1,16,22\}$&$\{1,17,18\}$\\$\{1,19,21\}$
		
		&$\{2,8,20\}$&$\{2,8,16\}$&$\{2,9,17\}$&$\{2,9,16\}$&$\{2,10,17\}$
		&$\{2,10,21\}$
		
		\\$\{2,11,24\}$&$\{2,11,21\}$&$\{2,12,23\}$&$\{2,12,25\}$&$\{2,13,20\}$&$\{2,13,18\}$
		
		&$\{2,14,15\}$\\$\{2,14,22\}$&$\{2,15,18\}$&$\{2,19,24\}$&$\{2,19,25\}$&$\{2,22,23\}$
		
		&$\{3,8,16\}$&$\{3,8,17\}$\\
\end{longtable}
\end{center}
A $ \Diamond $-ITS$_2 (29;9,8) $ on $I_{29}$ with holes  $ \{0,1,\ldots,8\} $ and $ \{9,10,\ldots,16\} $:
\begin{center}\tabcolsep 0.05in
	\begin{longtable}{lllllll}
		$\{0,9,17\}$
		&$\{0,9,22\}$
		&$\{0,10,23\}$
		&$\{0,10,17\}$
		&$\{0,11,19\}$&
		$\{0,11,24\}$
		
		&$\{0,12,28\}$\\$\{0,12,19\}$&$\{0,13,28\}$&$\{0,13,25\}$&
		$\{0,14,24\}$&$\{0,14,21\}$
		
		&$\{0,15,20\}$&$\{0,15,23\}$\\$\{0,16,22\}$&$\{0,16,20\}$&$\{0,18,27\}$&$\{0,18,25\}$
		
		&$\{0,21,26\}$&$\{0,26,27\}$&$\{1,9,26\}$
		\\$\{1,9,18\}$
		&$\{1,10,17\}$
		&$\{1,10,18\}$
		
		&$\{1,11,28\}$&
		$\{1,11,20\}$
		&$\{1,12,20\}$&$\{1,12,19\}$\\$\{1,13,17\}$
		&$\{1,13,19\}$
		&$\{1,14,24\}$&$\{1,14,25\}$
		&$\{1,15,27\}$&$\{1,15,26\}$&$\{1,16,23\}$\\$\{1,16,22\}$
		
		&$\{1,21,23\}$&$\{1,21,27\}$
		&$\{1,22,28\}$&$\{1,24,25\}$&$\{2,9,21\}$&$\{2,9,27\}$
		
		\\$\{2,10,21\}$&$\{2,10,18\}$&$\{2,11,19\}$&$\{2,11,22\}$
		&$\{2,12,23\}$&$\{2,12,17\}$
		
		&$\{2,13,18\}$\\$\{2,13,19\}$&$\{2,14,22\}$&$\{2,14,27\}$&$\{2,15,26\}$&$\{2,15,23\}$
		
		&$\{2,16,25\}$&$\{2,16,17\}$\\
\end{longtable}
\end{center}
A $ \Diamond $-ITS$_2 (27;10,4) $ on $I_{27}$ with holes  $ \{0,1,\ldots,9\} $ and $ \{10,11,12,13\} $:
\begin{center}\tabcolsep 0.05in
	\begin{longtable}{lllllll}
		
		$\{0,10,14\}$&$\{0,10,26\}$&$\{0,11,25\}$&$\{0,11,21\}$&$\{0,12,15\}$&$\{0,12,22\}$
		
		&$\{0,13,15\}$\\$\{0,13,23\}$&$\{0,14,21\}$&$\{0,16,25\}$&$\{0,16,18\}$&$\{0,17,23\}$
		
		&$\{0,17,18\}$&$\underline{2}\{0,19,24\}$\\$\{0,20,22\}$&$\{0,20,26\}$&$\{1,10,15\}$&$\{1,10,24\}$
		
		&$\{1,11,14\}$&$\{1,11,19\}$&$\{1,12,16\}$\\$\{1,12,24\}$
		&$\{1,13,21\}$
		&$\{1,13,19\}$
		&$\{1,14,22\}$&$\{1,15,22\}$&$\{1,16,23\}$&$\underline{2}\{1,17,20\}$\\$\underline{2}\{1,18,25\}$&$\{1,21,26\}$
		
		&$\{1,23,26\}$&$\{2,10,23\}$&$\{2,10,20\}$&$\{2,11,15\}$
		&$\{2,11,17\}$
		\\$\{2,12,22\}$
		
		&$\{2,12,23\}$&$\{2,13,24\}$&$\{2,13,19\}$&$\{2,14,26\}$&$\{2,14,21\}$&$\{2,15,22\}$\\$\{2,16,19\}$&$\{2,16,25\}$&$\{2,17,24\}$&$\{2,18,20\}$&$\{2,18,26\}$&$\{2,21,25\}$
		
		&$\{3,10,14\}$\\$\{3,10,21\}$&$\{3,11,14\}$&$\{3,11,25\}$&$\{3,12,24\}$&$\{3,12,18\}$
		
		&$\{3,13,26\}$&$\{3,13,23\}$\\
\end{longtable}
\end{center}
A $ \Diamond $-ITS$_2 (30;10,7) $ on $I_{30}$ with holes  $ \{0,1,\ldots,9\} $ and $ \{10,11,\ldots,16\} $:
\begin{center}\tabcolsep 0.05in
	\begin{longtable}{lllllll}
		
		$\{0,10,17\}$&$\{0,10,21\}$&$\{0,11,28\}$&$\{0,11,17\}$
		&$\{0,12,23\}$&
		$\{0,12,21\}$
		
		&$\{0,13,27\}$\\$\{0,13,18\}$&$\underline{2}\{0,14,19\}$&$\{0,15,20\}$&$\{0,15,23\}$&$\{0,16,25\}$
		&$\{0,16,22\}$&$\{0,18,28\}$\\$\{0,20,29\}$&$\{0,22,27\}$&$\{0,24,25\}$&$\{0,24,26\}$
		&$\{0,26,29\}$&$\{1,10,23\}$&$\{1,10,17\}$\\$\{1,11,18\}$&$\{1,11,21\}$&$\{1,12,18\}$
		&$\{1,12,22\}$&$\{1,13,19\}$&$\{1,13,29\}$&$\{1,14,21\}$\\$\{1,14,17\}$&$\{1,15,23\}$
		
		&$\{1,15,24\}$&$\{1,16,25\}$&$\{1,16,27\}$&$\{1,19,26\}$&$\{1,20,27\}$\\$\{1,20,25\}$
		
		&$\{1,22,29\}$&$\{1,24,28\}$&$\{1,26,28\}$&$\{2,10,18\}$&$\{2,10,25\}$&$\{2,11,23\}$
		
		\\$\{2,11,17\}$&$\{2,12,26\}$&$\{2,12,21\}$&$\{2,13,22\}$&$\{2,13,17\}$&$\{2,14,18\}$
		
		&$\{2,14,21\}$\\$\{2,15,25\}$&$\{2,15,24\}$&$\{2,16,27\}$&$\{2,16,23\}$&$\{2,19,28\}$
		&$\{2,19,26\}$&$\{2,20,24\}$\\
\end{longtable}
\end{center}
A $ \Diamond $-ITS$_2 (29;11,4) $ on $I_{29}$ with holes  $ \{0,1,\ldots,10\} $ and $ \{11,12,13,14\} $:
\begin{center}\tabcolsep 0.05in
	\begin{longtable}{lllllll}
		$\{0,11,15\}$
		&$\{0,11,24\}$&$\{0,12,25\}$
		&$\{0,12,15\}$
		&$\{0,13,20\}$&$\{0,13,28\}$
		
		&$\{0,14,22\}$\\$\{0,14,28\}$&$\{0,16,25\}$
		&$\{0,16,22\}$&$\{0,17,21\}$&$\{0,17,26\}$
		
		&$\{0,18,24\}$&$\{0,18,23\}$\\$\{0,19,21\}$&$\{0,19,26\}$&$\{0,20,27\}$&$\{0,23,27\}$
		
		&$\{1,11,28\}$&$\{1,11,20\}$&$\{1,12,17\}$\\$\{1,12,16\}$&$\underline{2}\{1,13,15\}$&$\{1,14,24\}$
		&$\{1,14,22\}$&$\{1,16,25\}$&$\{1,17,23\}$&$\{1,18,26\}$\\$\{1,18,20\}$&$\{1,19,24\}$
		
		&$\{1,19,22\}$&$\{1,21,25\}$&$\{1,21,27\}$&$\{1,23,27\}$&$\{1,26,28\}$\\$\{2,11,24\}$
		
		&$\{2,11,22\}$&$\{2,12,15\}$&$\{2,12,20\}$&$\{2,13,21\}$&$\{2,13,24\}$&$\{2,14,16\}$
		
		\\$\{2,14,19\}$&$\{2,15,20\}$&$\{2,16,26\}$&$\{2,17,25\}$&$\{2,17,28\}$&$\{2,18,27\}$
		
		&$\{2,18,23\}$\\$\{2,19,21\}$&$\{2,22,27\}$&$\{2,23,28\}$&$\{2,25,26\}$&$\{3,11,15\}$&$\{3,11,22\}$&$\{3,12,27\}$\\
\end{longtable}
\end{center}
A $ \Diamond $-ITS$_2 (35;11,10) $ on $I_{35}$ with holes  $ \{0,1,\ldots,10\} $ and $ \{11,12,\ldots,20\} $:
\begin{center}\tabcolsep 0.05in
	\begin{longtable}{lllllll}
		$\{0,11,21\}$
		&$\{0,11,22\}$
		&$\{0,12,22\}$
		&$\{0,12,24\}$
		&$\{0,13,26\}$&$\{0,13,21\}$
		&$\{0,14,23\}$\\$\{0,14,25\}$&$\{0,15,24\}$&$\{0,15,25\}$&$\{0,16,27\}$&$\{0,16,32\}$
		
		&$\{0,17,26\}$&$\{0,17,34\}$\\$\{0,18,32\}$&$\{0,18,31\}$
		&$\{0,19,30\}$&$\{0,19,29\}$
		
		&$\{0,20,30\}$&$\{0,20,29\}$&$\{0,23,31\}$\\$\{0,27,34\}$&$\underline{2}\{0,28,33\}$&$\{1,11,28\}$
		&$\{1,11,21\}$&$\{1,12,27\}$&$\{1,12,22\}$&$\{1,13,21\}$\\$\{1,13,22\}$&$\{1,14,25\}$
		
		&$\{1,14,34\}$&$\{1,15,30\}$&$\{1,15,29\}$&$\underline{2}\{1,16,26\}$&$\{1,17,28\}$\\$\{1,17,23\}$
		
		&$\{1,18,31\}$&$\{1,18,34\}$&$\{1,19,32\}$&$\{1,19,33\}$&$\{1,20,31\}$&$\{1,20,29\}$
		
		\\$\{1,23,27\}$&$\{1,24,33\}$&$\{1,24,30\}$&$\{1,25,32\}$&$\{2,11,22\}$&$\{2,11,23\}$
		&$\{2,12,23\}$\\$\{2,12,25\}$&$\{2,13,22\}$&$\{2,13,26\}$&$\{2,14,21\}$&$\{2,14,26\}$
		
		&$\{2,15,27\}$&$\{2,15,24\}$\\
\end{longtable}
\end{center}

\subsection{$\Diamond$-ITS$_6(2w+u+1;w,u)$s in  Lemma 4.32}
A $ \Diamond $-ITS$_6 (7;2,2) $ on $I_{7}$ with holes  $ \{0,1\} $ and $ \{2,3\} $:
\begin{center}\tabcolsep 0.05in
	\begin{longtable}{llllllll}
		\{0, 2, 5\}&\underline{5}\{0, 2, 6\}&\underline{3}\{0, 3, 4\}&\underline{3}\{0, 3, 5\}&\underline{2}\{0, 4, 5\}&\{0, 4, 6\} &
		\underline{3}\{1, 2, 4\}&\underline{3}\{1, 2, 5\}\\\underline{2}\{2, 4, 5\}&\{2, 4, 6\}&\underline{2}\{1, 3, 4\}&\{1, 3, 5\}&
		\underline{3}\{1, 3, 6\}&\{1, 4, 6\}&\underline{2}\{1, 5, 6\}&\{3, 4, 6\}\\\underline{2}\{3, 5, 6\}&\underline{2}\{4, 5, 6\}
	\end{longtable}
\end{center}
A $ \Diamond $-ITS$_6 (9;3,2) $ on $I_{9}$ with holes  $ \{0,1\} $ and $ \{2,3,4\} $:
\begin{center}\tabcolsep 0.05in
	\begin{longtable}{llllllll}
		\underline{6}\{0, 2, 5\}&\underline{4}\{0, 3, 6\}&\{0, 3, 7\}&\{0, 3, 8\}&\underline{4}\{0, 4, 7\}&\underline{2}\{0, 4, 8\}&\underline{2}\{0, 6, 8\}&\{0, 7, 8\}\\\{1, 2, 6\}&\underline{2}\{1, 2, 7\}&\underline{3}\{1, 2, 8\}&\underline{3}\{2, 6, 7\}&\underline{2}\{2, 6, 8\}&\{2, 7, 8\}&\underline{3}\{1, 3, 5\}&\{1, 3, 6\}\\\underline{2}\{1, 3, 8\}&\{3, 5, 6\}&\underline{2}\{3, 5, 7\}&\underline{3}\{3, 7, 8\}&\underline{4}\{1, 4, 6\}&\underline{2}\{1, 4, 7\}&\underline{2}\{1, 5, 7\}&\{1, 5, 8\}\\\underline{2}\{4, 5, 6\}&\underline{4}\{4, 5, 8\}&\underline{2}\{5, 6, 7\}&\{5, 6, 8\}&\{6, 7, 8\}
	\end{longtable}
\end{center}
A $ \Diamond $-ITS$_6 (13;5,2) $ on $I_{13}$ with holes  $ \{0,1,2,3,4\} $ and $ \{5,6\} $:
\begin{center}\tabcolsep 0.05in
	\begin{longtable}{lllllll}
		\underline{6}\{0, 5, 7\}&\underline{6}\{0, 6, 8\}&\underline{6}\{0, 9, 10\}&\underline{6}\{0, 11, 12\}&\underline{6}\{1, 5, 8\} &
		\underline{6}\{1, 6, 9\}&\underline{6}\{1, 7, 11\}\\\underline{6}\{1, 10, 12\}&\underline{6}\{2, 5, 9\}&\underline{6}\{3, 5, 12\} &
		\underline{3}\{4, 5, 10\}&\underline{3}\{4, 5, 11\}&\underline{3}\{5, 10, 11\}&\underline{3}\{2, 6, 11\} \\
		\underline{2}\{2, 8, 11\}&\{2, 10, 11\}&\{2, 6, 7\}&\{2, 6, 10\}&\{2, 6, 12\}&\underline{3}\{2, 7, 10\}&\{2, 8, 10\}\\\underline{2}\{2, 7, 12\}&\underline{3}\{2, 8, 12\}&\underline{4}\{4, 6, 12\}&\{6, 7, 12\}&\underline{4}\{3, 6, 7\}&\underline{2}\{3, 6, 10\}&\{4, 6, 10\}\\\{4, 6, 11\} &\underline{2}\{6, 10, 11\}&\underline{2}\{3, 7, 10\}&\underline{2}\{3, 8, 10\} &
		\underline{2}\{3, 8, 9\}&\underline{2}\{3, 8, 11\}&\underline{4}\{3, 9, 11\}\\\underline{2}\{4, 8, 10\}&\{7, 8, 10\}&\{4, 8, 11\}&\{4, 9, 11\}&\{8, 9, 11\}&\underline{3}\{4, 7, 8\}&\underline{3}\{4, 7, 9\}\\\underline{2}\{4, 9, 12\}&\{7, 8, 9\}&\{7, 8, 12\} &
		\underline{2}\{7, 9, 12\}&\underline{2}\{8, 9, 12\}
	\end{longtable}
\end{center}
A $ \Diamond $-ITS$_6 (15;5,4) $ on $I_{15}$ with holes  $ \{0,1,2,3,4\} $ and $ \{5,6,7,8\} $:
\begin{center}\tabcolsep 0.05in
	\begin{longtable}{lllllll}
		$\{0,5,10\}$
		&$\underline{4}\{0,5,9\}$
		&$\{0,5,11\}$
		
		&$\{0,6,14\}$
		&$\underline{3}\{0,6,11\}$
		&$\underline{2}\{0,6,10\}$
		&$\{0,7,13\}$
		\\$\{0,7,14\}$
		&$\{0,7,9\}$&$\underline{2}\{0,7,12\}$
		&$\{0,7,13\}$
		&$\underline{3}\{0,8,14\}$
		
		&$\{0,8,11\}$&$\{0,8,13\}$
		\\$\{0,8,12\}$
		&$\{0,9,12\}$&$\{0,10,13\}$
		&$\underline{2}\{0,10,12\}$
		
		&$\{0,11,13\}$
		&$\{0,13,14\}$&$\underline{2}\{1,5,14\}$
		\\$\underline{2}\{1,5,10\}$
		&$\{1,5,12\}$
		&$\{1,5,11\}$
		&$\{1,6,10\}$
		&$\underline{2}\{1,6,14\}$
		&$\underline{2}\{1,6,13\}$
		&$\{1,6,9\}$
		\\$\{1,7,10\}$
		&$\{1,7,9\}$
		
		&$\underline{2}\{1,7,12\}$&$\{1,7,14\}$
		&$\{1,7,9\}$
		&$\underline{2}\{1,8,13\}$
		&$\{1,8,11\}$
		\\$\underline{2}\{1,8,12\}$
		
		&$\{1,8,10\}$
		&$\{1,9,13\}$&$\underline{2}\{1,9,11\}$
		&$\{1,10,12\}$
		&$\{1,11,13\}$
		&$\{1,11,14\}$
		\\$\underline{2}\{2,5,9\}$
		&$\{2,5,13\}$
		&$\{2,5,14\}$&
		$\underline{2}\{2,5,11\}$
		&$\{2,6,10\}$
		&$\underline{3}\{2,6,13\}$
		
		&$\{2,6,11\}$
		\\$\{2,6,14\}$&$\{2,7,12\}$
		&$\{2,7,11\}$
		&$\{2,7,9\}$&$\{2,7,14\}$
		
		&$\{2,7,13\}$
		&$\{2,7,11\}$
		\\
\end{longtable}
\end{center}
A $ \Diamond $-ITS$_6 (15;6,2) $ on $I_{15}$ with holes  $ \{0,1,\ldots,5\} $ and $ \{6,7\} $:
\begin{center}\tabcolsep 0.05in
	\begin{longtable}{lllllll}
		$\underline{2}\{0,6,9\}$
		&$\{0,6,8\}$
		&$\{0,6,12\}$
		&$\underline{2}\{0,6,10\}$
		&$\underline{2}\{0,7,12\}$
		&$\underline{2}\{0,7,9\}$
		
		&$\{0,7,10\}$\\$\{0,7,11\}$
		&$\underline{3}\{0,8,13\}$
		&$\{0,8,14\}$
		&$\{0,8,11\}$&$\{0,9,14\}$
		
		&$\{0,9,12\}$
		&$\{0,10,11\}$\\$\underline{2}\{0,10,14\}$
		&$\{0,11,14\}$
		&$\underline{2}\{0,11,13\}$&$\{0,12,13\}$
		
		&$\{0,12,14\}$
		&$\underline{5}\{1,6,8\}$
		&$\{1,6,10\}$
		
		\\$\underline{2}\{1,7,12\}$
		&$\underline{2}\{1,7,13\}$&$\{1,7,11\}$
		
		&$\{1,7,14\}$
		&$\{1,8,11\}$&$\{1,9,10\}$
		&$\underline{2}\{1,9,13\}$
		\\$\underline{2}\{1,9,12\}$
		&$\{1,9,14\}$
		&$\{1,10,14\}$
		&$\{1,10,11\}$
		&$\{1,10,13\}$
		&$\{1,10,12\}$&$\{1,11,13\}$
		\\$\underline{2}\{1,11,14\}$
		
		&$\{1,12,14\}$
		&$\underline{3}\{2,6,10\}$&$\{2,6,9\}$
		&$\{2,6,12\}$
		&$\{2,6,11\}$&$\{2,7,11\}$
		
		\\$\underline{3}\{2,7,13\}$
		&$\{2,7,12\}$
		&$\{2,7,14\}$&$\{2,8,10\}$
		&$\{2,8,13\}$
		&$\{2,8,12\}$
		
		&$\underline{2}\{2,8,14\}$\\$\{2,8,9\}$
		&$\underline{2}\{2,9,10\}$
		&$\{2,9,13\}$&$\{2,9,12\}$
		&$\underline{3}\{2,11,14\}$
		
		&$\{2,11,12\}$
		&$\{2,12,13\}$
		\\
\end{longtable}
\end{center}
A $ \Diamond $-ITS$_6 (17;6,4) $ on $I_{17}$ with holes  $ \{0,1,\ldots,5\} $ and $ \{6,7,8,9\} $:
\begin{center}\tabcolsep 0.05in
	\begin{longtable}{lllllll}
$\underline{2}\{0,6,10\}$
&$\underline{2}\{0,6,13\}$
&$\{0,6,11\}$
&$\{0,6,14\}$
&$\underline{2}\{0,7,11\}$
&$\{0,7,13\}$
&$\underline{2}\{0,7,10\}$
\\$\{0,7,16\}$&
$\underline{2}\{0,8,13\}$
&$\underline{2}\{0,8,15\}$
&$\{0,8,14\}$
&$\{0,8,10\}$
&$\underline{2}\{0,9,16\}$
&$\{0,9,13\}$\\
$\underline{2}\{0,9,14\}$
&$\{0,9,10\}$
&$\{0,11,16\}$&
$\underline{2}\{0,11,15\}$
&$\underline{2}\{0,12,16\}$
&$\underline{2}\{0,12,14\}$
&$\underline{2}\{0,12,15\}$
\\$\{1,6,15\}$
&$\{1,6,10\}$&
$\{1,6,11\}$
&$\{1,6,12\}$
&$\{1,6,13\}$
&$\{1,6,14\}$
&$\underline{2}\{1,7,13\}$
\\$\{1,7,11\}$
&$\underline{2}\{1,7,10\}$
&$\{1,7,14\}$
&$\underline{2}\{1,8,12\}$
&$\{1,8,10\}$
&$\{1,8,15\}$
&$\{1,8,14\}$
\\$\{1,8,13\}$
&$\underline{2}\{1,9,10\}$
&$\underline{3}\{1,9,16\}$
&$\{1,9,15\}$
&$\underline{3}\{1,11,15\}$
&$\{1,11,14\}$
&$\{1,12,13\}$
\\$\underline{2}\{1,12,16\}$
&$\{1,13,14\}$&
$\{1,14,16\}$
&$\{2,6,10\}$
&$\underline{3}\{2,6,13\}$
&$\{2,6,14\}$
&$\{2,6,16\}$
\\$\underline{2}\{2,7,10\}$
&$\{2,7,13\}$
&$\{2,7,11\}$
&$\underline{2}\{2,7,15\}$
&$\{2,8,10\}$
&$\{2,8,16\}$&
$\underline{2}\{2,8,11\}$
\\
\end{longtable}
\end{center}
A $ \Diamond $-ITS$_6 (18;6,5) $ on $I_{18}$ with holes  $ \{0,1,\ldots,5\} $ and $ \{6,7,8,9,10\} $:
\begin{center}\tabcolsep 0.05in
	\begin{longtable}{lllllll}
		$\{0,6,11\}$
		&$\underline{2}\{0,6,13\}$
		&$\underline{2}\{0,6,12\}$
		&$\{0,6,17\}$
		&$\underline{2}\{0,7,15\}$
		&$\underline{2}\{0,7,14\}$
		
		&$\{0,7,12\}$
		\\$\{0,7,11\}$
		&$\{0,8,11\}$
		&$\underline{2}\{0,8,17\}$
		&$\underline{2}\{0,8,15\}$&$\{0,8,16\}$
		
		&$\{0,9,15\}$&$\underline{2}\{0,9,12\}$\\$\{0,9,14\}$
		&$\{0,9,13\}$
		&$\{0,9,11\}$
		&$\{0,10,17\}$
		&$\underline{2}\{0,10,16\}$&$\underline{2}\{0,10,14\}$&$\{0,10,13\}$\\$\{0,11,15\}$
		&$\{0,11,17\}$&$\{0,12,16\}$
		
		&$\underline{2}\{0,13,16\}$&$\{0,14,17\}$
		&$\{1,6,15\}$&$\underline{3}\{1,6,11\}$\\$\{1,6,12\}$&$\{1,6,13\}$
		
		&$\underline{2}\{1,7,16\}$
		&$\underline{2}\{1,7,14\}$
		&$\{1,7,17\}$
		&$\{1,7,13\}$
		&$\{1,8,15\}$
		\\$\underline{2}\{1,8,14\}$
		&$\{1,8,16\}$&$\underline{2}\{1,8,12\}$&$\underline{3}\{1,9,17\}$&$\underline{2}\{1,9,13\}$
		&$\{1,9,11\}$
		&$\{1,10,11\}$
		
		\\$\underline{2}\{1,10,13\}$
		&$\{1,10,17\}$&$\{1,10,15\}$&$\{1,10,16\}$&$\{1,11,16\}$&$\underline{3}\{1,12,15\}$
		
		&$\{1,14,17\}$
		\\$\{1,14,16\}$
		&$\underline{2}\{2,6,11\}$
		&$\{2,6,13\}$&$\underline{3}\{2,6,12\}$&$\underline{2}\{2,7,17\}$
		
		&$\underline{3}\{2,7,12\}$&$\{2,7,11\}$\\
\end{longtable}
\end{center}
A $ \Diamond $-ITS$_6 (21;8,4) $ on $I_{21}$ with holes  $ \{0,1,\dots,7\} $ and $ \{8,9,10,11\} $:
\begin{center}\tabcolsep 0.05in
	\begin{longtable}{lllllll}
		$\{0,8,13\}$
		&$\underline{2}\{0,8,20\}$
		&$\underline{3}\{0,8,14\}$
		&$\underline{2}\{0,9,17\}$
		&$\underline{2}\{0,9,18\}$
		&$\{0,9,19\}$
		&$\{0,9,15\}$
		\\$\underline{2}\{0,10,19\}$
		&$\{0,10,18\}$
		&$\{0,10,12\}$
		&$\{0,10,14\}$
		&$\{0,10,15\}$
		&$\underline{3}\{0,11,18\}$
		&$\{0,11,15\}$
		\\$\{0,11,13\}$
		&$\{0,11,14\}$
		&$\{0,12,17\}$&
		$\underline{3}\{0,12,16\}$
		&$\{0,12,19\}$
		&$\{0,13,14\}$&
		$\{0,13,15\}$
		\\$\{0,13,17\}$
		&$\{0,13,20\}$
		&$\{0,15,17\}$
		&$\{0,15,19\}$&
		$\underline{3}\{0,16,20\}$
		&$\{0,17,19\}$
		&$\{1,8,15\}$\\
		$\{1,8,14\}$
		&$\{1,8,16\}$
		&$\{1,8,20\}$
		&$\{1,8,19\}$
		&$\{1,8,13\}$&
		$\{1,9,14\}$
		&$\{1,9,16\}$
		\\$\{1,9,20\}$
		&$\{1,9,13\}$
		&$\{1,9,12\}$&
		$\{1,9,15\}$
		&$\underline{2}\{1,10,14\}$
		&$\underline{3}\{1,10,18\}$
		&$\{1,10,17\}$
		\\$\underline{2}\{1,11,12\}$
		&$\underline{2}\{1,11,20\}$
		&$\{1,11,17\}$
		&$\{1,11,19\}$
		&$\{1,12,14\}$
		&$\{1,12,17\}$&
		$\{1,12,19\}$
		\\$\{1,13,16\}$
		&$\{1,13,20\}$
		&$\underline{2}\{1,13,19\}$
		&$\{1,14,17\}$
		&$\{1,15,16\}$
		&$\{1,15,18\}$
		&$\underline{2}\{1,15,17\}$
		\\
\end{longtable}
\end{center}
A $ \Diamond $-ITS$_6 (22;8,5) $ on $I_{22}$ with holes  $ \{0,1,\ldots,7\} $ and $ \{8,9,10,11,12\} $:
\begin{center}\tabcolsep 0.03in
	\begin{longtable}{lllllll}
		$\{0,8,14\}$&$\underline{2}\{0,8,18\}$&$\underline{2}\{0,8,20\}$
		&$\{0,8,17\}$&$\{0,9,13\}$&$\underline{4}\{0,9,16\}$
		
		&$\{0,9,14\}$\\$\underline{3}\{0,10,19\}$&$\{0,10,13\}$&$\{0,10,18\}$&$\{0,10,16\}$&$\underline{4}\{0,11,20\}$
		
		&$\{0,11,21\}$&$\{0,11,14\}$\\$\underline{2}\{0,12,19\}$&$\underline{2}\{0,12,17\}$&$\{0,12,13\}$&$\{0,12,18\}$&$\{0,13,21\}$&$\underline{2}\{0,13,18\}$&$\{0,14,15\}$\\$\underline{2}\{0,14,17\}$&$\underline{4}\{0,15,21\}$&$\{0,15,17\}$
		
		&$\{0,16,19\}$&$\underline{2}\{1,8,17\}$&$\{1,8,14\}$&$\underline{2}\{1,8,16\}$\\$\{1,8,13\}$
		&$\underline{2}\{1,9,21\}$
		&$\{1,9,16\}$&$\{1,9,15\}$&$\{1,9,13\}$&$\{1,9,18\}$&$\underline{2}\{1,10,18\}$\\$\{1,10,13\}$
		
		&$\underline{2}\{1,10,14\}$&$\{1,10,19\}$&$\{1,11,16\}$&$\{1,11,21\}$
		&$\underline{2}\{1,11,14\}$&$\underline{2}\{1,11,19\}$
		
		\\$\{1,12,21\}$&$\underline{4}\{1,12,15\}$&$\{1,12,18\}$&$\underline{3}\{1,13,20\}$&$\{1,14,18\}$&$\{1,15,17\}$
		
		&$\{1,16,20\}$\\$\{1,16,19\}$&$\{1,17,19\}$&$\{1,17,21\}$&$\{1,17,20\}$&$\{1,18,20\}$
		
		&$\{1,19,21\}$&$\{2,8,17\}$\\
\end{longtable}
\end{center}
A $ \Diamond $-ITS$_6 (24;8,7) $ on $I_{24}$ with holes  $ \{0,1,\ldots,7\} $ and $ \{8,9,\ldots,14\} $:
\begin{center}\tabcolsep 0.03in
	\begin{longtable}{lllllll}
		$\{0,8,21\}$&$\underline{3}\{0,8,15\}$&$\{0,8,18\}$&$\{0,8,16\}$&$\underline{2}\{0,9,15\}$&$\{0,9,22\}$
		
		&$\{0,9,17\}$\\$\{0,9,18\}$&$\{0,9,21\}$&$\underline{2}\{0,10,18\}$&$\{0,10,21\}$&$\{0,10,17\}$
		&$\{0,10,23\}$&$\{0,10,22\}$\\$\underline{2}\{0,11,21\}$&$\{0,11,17\}$&$\{0,11,15\}$&$\{0,11,20\}$
		
		&$\{0,11,16\}$&$\underline{2}\{0,12,20\}$&$\{0,12,16\}$\\$\{0,12,23\}$&$\underline{2}\{0,12,22\}$&$\{0,13,19\}$
		
		&$\{0,13,18\}$&$\{0,13,16\}$&$\{0,13,16\}$&$\{0,13,20\}$\\$\{0,13,23\}$&$\underline{2}\{0,14,23\}$
		
		&$\underline{2}\{0,14,17\}$&$\underline{2}\{0,14,22\}$&$\{0,16,20\}$&$\{0,17,19\}$&$\{0,18,19\}$\\$\{0,19,20\}$
		
		&$\{0,19,23\}$&$\{0,19,21\}$&$\{1,8,16\}$&$\underline{2}\{1,8,19\}$&$\{1,8,22\}$&$\underline{2}\{1,8,21\}$
		
		\\$\underline{4}\{1,9,15\}$&$\{1,9,17\}$&$\{1,9,16\}$&$\underline{2}\{1,10,15\}$&$\{1,10,22\}$&$\underline{2}\{1,10,16\}$
		
		&$\{1,10,20\}$\\$\{1,11,19\}$&$\underline{2}\{1,11,20\}$&$\underline{2}\{1,11,16\}$&$\{1,11,17\}$&$\{1,12,22\}$
		
		&$\{1,12,23\}$&$\{1,12,20\}$\\
\end{longtable}
\end{center}

\subsection{A case in  Lemma 6.6}
A $ \Diamond $-ITS$_3 (11;3,3) $ on $I_{11}$ with holes  $ \{0,1,2\} $ and $ \{3,4,5\} $:
\begin{center}\tabcolsep 0.05in
	\begin{longtable}{llllllll}
		\underline{3}\{0, 3, 6\}&\underline{3}\{0, 4, 7\}&\underline{3}\{0, 5, 8\}&\underline{3}\{0, 9, 10\} &\underline{3}\{1, 3, 7\}
		&\underline{3}\{2, 3, 9\} &\underline{3}\{3, 8, 10\}&\{1, 5, 6\} \\\underline{2}\{1, 5, 9\} &\underline{3}\{1, 4, 10\}
		&\underline{2}\{1, 6, 8\} &\{1, 8, 9\} &\{2, 4, 6\}&\underline{2}\{2, 4, 8\} &\{2, 7, 8\} &\underline{2}\{4, 6, 9\}
		\\\{4, 8, 9\} &\{5, 7, 9\} &\{6, 7, 8\} &\{6, 7, 9\}
		&\{7, 8, 9\} &\{2, 5, 6\} &\{2, 5, 7\} &\{2, 5, 10\} \\\{2, 6, 10\}
		&\{2, 7, 10\} &\{5, 6, 10\} &\{5, 7, 10\}& \{6, 7, 10\}
	\end{longtable}
\end{center}

\end{document}